\documentclass{sigma}
\usepackage{graphicx}
\usepackage[nobysame]{myamsrefs}
\PassOptionsToPackage{pdfauthor={Vladimir V. Kisil},%
    pdftitle={Erlangen Program at Large--1: Geometry of Invariants},%
    pdfsubject={mathematics},%
  pdfkeywords={symmetries, cycles, analytic function theory, semisimple groups, elliptic,
  parabolic, hyperbolic, Clifford algebras, complex numbers, dual
  numbers, double numbers, split-complex numbers, Moebius
  transformations},
 }{hyperref}
\usepackage{amssymb}
\usepackage{xr-hyper}
\usepackage{hyperref}

\providecommand{\GiNaC}{\textsf{GiNaC}}

\providecommand{\comment}[1]{}

\providecommand{\such}{\,\mid\,}
\providecommand{\norm}[2][\relax]{\left\|#2\right\|\ifx#1\relax\else_{#1}\fi}
\providecommand{\modulus}[2][\relax]{\left| #2 \right|\ifx#1\relax\else_{#1}\fi}
\providecommand{\lvec}[1]{\overrightarrow{#1}}
\newcommand{\zcycle}[3][]{#1 Z^{#2}_{#3}}
\newcommand{\realline}[3][]{#1 R^{#2}_{#3}}

\providecommand{\spec}[1][]{\ensuremath{\mathbf{sp}}\,}
\providecommand{\GiNaC}{\textsf{GiNaC}}
\providecommand{\bs}{\breve{\sigma}}
\providecommand{\SL}[1][2]{\ensuremath{\FSpace{SL}{#1}(\Space{R}{})}}
\providecommand{\scalar}[3][\relax]{\left\langle #2,#3 
        \right\rangle\ifx#1\relax\else_{#1}\fi}
\providecommand{\Space}[3][]{\ensuremath{\mathbb{#2}^{#3}_{#1}{}}}
  \providecommand{\FSpace}[3][]{\ensuremath{\ifx#2l \ell_{#3}^{#1}{}\else
  #2_{#3}^{#1}{}\fi}} 
\providecommand{\rmi}{\mathrm{i}}
\providecommand{\rme}{\varepsilon}

\providecommand{\tr}{\mathop{tr}}
\providecommand{\algebra}[1]{\ensuremath{\mathfrak{#1}}}

\providecommand{\modulus}[2][\relax]{\left| #2 \right|\ifx#1\relax\else_{#1}\fi}

\providecommand{\href}[2]{#2}

\newcommand{\CPP}{\texttt{C++}}
\newcommand{\NoWEB}{\texttt{noweb}}

\newcommand{\sperp}{\dashv}

\makeatletter
\newcommand{\ephname}[1]{\ifcase\csname c@#1\endcsname a\or e\or p\or h\fi}
\makeatother
\input{cyr.fd}

\providecommand{\tr}{\mathop{tr}}
\providecommand{\GiNaC}{\textsf{GiNaC}}

\providecommand{\se}[1]{\breve{e}_#1}
\providecommand{\DSpace}[2]{\ensuremath{ { \dot{\mathbb{#1}}^{#2}} }}
\providecommand{\TSpace}[2]{\ensuremath{ { \widetilde{\mathbb{#1}}^{#2}} }}
\providecommand{\cycle}[3][]{{#1 C^{#2}_{#3}}}
\providecommand{\realline}[3][]{#1 R^{#2}_{#3}}
\providecommand{\zcycle}[3][]{#1 Z^{#2}_{#3}}
\providecommand{\bs}{\breve{\sigma}}
\providecommand{\rs}{\mathring{\sigma}}
\providecommand{\half}{{\textstyle\frac{1}{2}}}
\providecommand{\Cliff}[2][\comment]{{\ensuremath{%
\mathcal{C}\kern-0.12em\ell(#1,#2)}}}
\usepackage[all]{xy}

\newtheorem{problem}{Problem}
\newtheorem{conv}{Convention}

 \numberwithin{equation}{section}
 \numberwithin{theorem}{section}
 \numberwithin{proposition}{section}
 \numberwithin{lemma}{section}
 \numberwithin{corollary}{section}
 \numberwithin{definition}{section}
 \numberwithin{example}{section}
 \numberwithin{remark}{section}
 \numberwithin{note}{section}
\numberwithin{problem}{section}
\numberwithin{conv}{section}
\newcommand{\epigraph}[3]{\par
\hfill\parbox{0.6\textwidth}{\footnotesize #1 \par \hfil #2 
\textit{#3}}\par\medskip}%
\providecommand{\person}[1]{#1}

\begin{document}

\renewcommand{\PaperNumber}{076}
\renewcommand{\Volume}{{\bf 6}}
\renewcommand{\PublicationYear}{2010}

\FirstPageHeading

\ShortArticleName{EPAL1: Geometry of Invariants}

\ArticleName{Erlangen Program at Large--1: Geometry of Invariants}

\Author{Vladimir V. Kisil}

\AuthorNameForHeading{V. V. Kisil}

\Address{
School of Mathematics,
University of Leeds, 
Leeds LS2\,9JT, 
UK\\
On  leave from the Odessa University}

\Email{kisilv@maths.leeds.ac.uk}

\URLaddress{http://www.maths.leeds.ac.uk/\~{}kisilv/}

\ArticleDates{Received April 20, 2010, in f\/inal form September 10, 2010;  Published online September 26, 2010}

\Abstract{
  This paper presents geometrical foundation for a systematic
  treatment of three main (elliptic, parabolic and hyperbolic) types
  of analytic function theories based on the representation theory of
  \(\SL\) group. We describe here geometries of corresponding domains.
  The principal r\^ole is played by Clifford algebras of matching
  types. In this paper we also generalise the
  Schwerdtfeger--Fillmore--Springer--Cnops construction which describes cycles as
  points in the extended space. This allows to consider many algebraic
  and geometric invariants of cycles within the Erlangen program
  approach.
}
\Classification{30G35, 22E46, 30F45, 32F45}


\epigraph{\noindent\textcyr{I pryamiznu tetivy lomayu,\\
    i luk sgibaet{}sya v krug\ldots}}{\textcyr{A.V.~Makarevich}}{}
\section{Introduction}
\label{sec:introduction}

This paper describes geometry of simples two-dimensional domains in
the spirit of the Erlangen program of \person{F.~Klein} influenced by
works of \person{S.~Lie}, for its development see
books~\citelist{\cite{BarkerHowe07} \cite{Sharpe97}
  \cite{RozenfeldZamakhovski03}} and their references. Further works
in this series will use the Erlangen approach for analytic function
theories and spectral theory of operators~\cite{Kisil10b}. In the
present paper we are focused on the geometry and study objects in a
plane and their properties which are \emph{invariant} under
linear-fractional transformations associated to the \(\SL\) group. The
basic observation is that geometries obtained in this way are
naturally classified as \emph{elliptic}, \emph{parabolic} and
\emph{hyperbolic}.

\subsection{Background and History}
\label{sec:background}
We repeatedly meet such a division of various mathematical objects
into three main classes.  They are named by the historically first
example---the classification of conic sections: elliptic, parabolic and
hyperbolic---however the pattern persistently reproduces itself in many
very different areas (equations, quadratic forms, metrics, manifolds,
operators, etc.). We will abbreviate this separation as
\emph{EPH-classification}. The \emph{common origin} of this
fundamental division can be seen from the simple picture of a
coordinate line split by the zero into negative and positive
half-axes:
\begin{equation}
  \label{eq:eph-class}
    \raisebox{-15pt}{\includegraphics[scale=1]{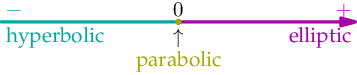}}
\end{equation}

Connections between different objects admitting EPH-classification are
not limited to this common source. There are many deep results
linking, for example, ellipticity of quadratic forms, metrics and
operators. On the other hand there are still a lot of white spots and
obscure gaps between some subjects as well.

For example, it is well known that elliptic operators are effectively
treated through complex analysis, which can be naturally identified as
the \emph{elliptic analytic function
  theory}~\cites{Kisil97c,Kisil01a}.  Thus there is a natural quest for
\emph{hyperbolic} and \emph{parabolic} analytic function theories,
which will be of similar importance for corresponding types of
operators. A search for hyperbolic function theory was attempted
several times starting from 1930's, see for
example~\cites{VignauxDuranona35a,LavrentShabat77,MotterRosa98}.
Despite of some important advances the obtained hyperbolic theory does
not look as natural and complete as complex analysis is.  
Parabolic geometry was considered in an excellent
book~\cite{Yaglom79}, which is still a source of valuable
inspirations.  However the corresponding ``parabolic calculus''
described in various places~\citelist{\cite{CatoniCannataNichelatti04}
  \cite{Gromov90a} \cite{Zejliger34}} is rather trivial.  

There is also a recent interest to this topic in different areas:
differential geometry \cites{BekkaraFrancesZeghib06,%
  CatoniCannataNichelatti04,CatoniCannataZampetti05,
  FjelstadGal01,Benz07a,BairdWood09a}, modal
logic~\cite{KuruczWolterZakharyaschev05}, quantum
mechanics~\cites{Khrennikov05a,KhrennikovSegre07a,McRae07,Kisil10a}, space-time
geometry~\cites{BocCatoniCannataNichZamp07,HerranzSantander02b,%
  HerranzSantander02a,Mirman01a,GromovKuratov05a,Garasko09a,Pavlov06a},
hypercomplex
analysis~\cites{CerejeirasKahlerSommen05a,EelbodeSommen04a,EelbodeSommen05a}.
A brief history of the topic can be found
in~\cite{CatoniCannataZampetti05} and further references are provided
in the above papers.

Most of previous research had an algebraic flavour. An alternative
approach to analytic function theories based on the representation
theory of semisimple Lie groups was developed in the series of
papers~\cites{Kisil95d,Kisil96d,Kisil96c,Kisil97c,Kisil94e,Kisil98a,Kisil01a}.
Particularly, some elements of hyperbolic function theory were built
in~\cites{Kisil96d,Kisil97c} along the same lines as the elliptic
one---standard complex analysis. Covariant functional calculus of
operators and respective covariant spectra were considered
in~\cites{Kisil95i,Kisil02a}. 

This paper continues this line of research and significantly expands
results of the earlier paper~\cite{Kisil04b}, see also~\cite{Kisil06a}
for an easy-reading introduction. A brief outline of the Erlangen
Programme at Large, which includes geometry, analytic functions and
functional calculus, is written in~\cite{Kisil10b}. 

\subsection{Highlights of Obtained Results}
\label{sec:novelty-highlights}

In the previous
paper~\cite{Kisil04b} we identify geometric objects called
\emph{cycles}~\cite{Yaglom79}, which are circles, parabolas and
hyperbolas in the corresponding EPH cases. They are invariants of the
M\"obius transformations, i.e. the natural geometric objects in the
sense of the Erlangen program. Note also that cycles are algebraically
defined through the quadratic expressions~\eqref{eq:cycle-def-real}
which may lead to interesting connections with the innovative approach
to the geometry presented in~\cite{Wildberger05}.

In this paper we systematically study those cycles through
an essential extension of the Schwerdtfeger--Fillmore--Springer--Cnops
construction\footnote{In
  the case of circles this technique was already spectacularly
  developed by H.~Schwerdtfeger in 1960-ies, see
  \cite{Schwerdtfeger79a}. Unfortunately, that beautiful book was not
  known to the present author until he accomplished his own
  works~\cites{Kisil05a,Kisil05b,Kisil12a}.}~\cites{Schwerdtfeger79a,Cnops02a,Porteous95,FillmoreSpringer90a,Kirillov06,Kisil12a} abbreviated in this paper as
SFSCc. The idea behind SFSCc is to consider cycles not as loci of points
from the initial \emph{point space} but rather as points of the new
\emph{cycle space}, see \S~\ref{sec:fillm-spring-cnops}. Then many
geometrical properties of the point space may be better expressed
through properties of the cycle space. Notably M\"obius
linear-fractional transformations of the point space are linearised in
the cycle space, see Prop.~\ref{pr:transform-cycles}.

An interesting feature of the correspondence between the point and
cycle spaces is that many relations between cycles, which are of local
nature in the cycle space, looks like non-local if translated back to
the point space, see for example non-local character of cycle
orthogonality in Fig.~\ref{fig:orthogonality1}
and~\ref{fig:orthogonality2}. Such a non-point behaviour is oftenly
thought to be a characteristic property of \emph{non-commutative
  geometry} but appears here within the Erlangen program
approach~\cites{Kisil97a,Kisil02c}. This also demonstrates that our
results cannot be reduced to the ordinary differential geometry or
nine Cayley--Klein geometries of the
plane~\citelist{\cite{Yaglom79}*{App.~C} \cite{Pimenov65a}}. 
 
\begin{remark}
  \label{rem:prejudice}
  Introducing parabolic objects on a common ground with elliptic
  and hyperbolic ones we should warn against some common prejudices
  suggested by picture~\eqref{eq:eph-class}:
  \begin{enumerate}
  \item \label{it:par-is-zero} The parabolic case is unimportant (has
    ``zero measure'') in comparison to the elliptic and hyperbolic ones.
    As we shall see (e.g. Remark~\ref{re:par-more-cayley}
    and~\ref{re:focal-lenght}.\ref{it:more-lengths}) some geometrical features are richer in
    parabolic case.
  \item 
    \label{it:par-is-limit} The parabolic case is a limiting situation
    (a contraction) or an intermediate position between the elliptic
    and hyperbolic ones: all properties of the former can be guessed
    or obtained as a limit or an average from the latter two.
    Particularly this point of view is implicitly supposed
    in~\cite{LavrentShabat77}.
    
    Although there are few confirmations of this (e.g.
    Fig.~\ref{fig:unit-disks}(\(E\))--(\(H\))), we shall see (e.g.
    Remark~\ref{re:par-is-not-limit}) that some properties of the
    parabolic case cannot be straightforwardly guessed from a
    combination of the elliptic and hyperbolic cases.
  \item All three EPH cases are even less disjoint than it is usually
    thought. For example, there are meaningful notions of \emph{centre of
      a parabola}~\eqref{de:center-first} or \emph{focus of
      a circle}~\eqref{eq:def-focus}.
  \item A (co-)invariant geometry is believed to be ``coordinate free''
    which sometimes is pushed to an absolute mantra. However our
    study within the Erlangen program framework reveals two useful notions
    (Defn.~\ref{de:center-first} and ~\eqref{eq:def-focus}) mentioned
    above which are defined by coordinate expressions and look very
    ``non-invariant'' on the first glance. 
  \end{enumerate}
\end{remark}

An amazing aspect of this topic is a transparent similarity between
all three EPH cases which is combined with some non-trivial exceptions
like \emph{non-invariance} of the upper half-plane in the hyperbolic
case (subsection~\ref{sec:invar-upper-half}) or \emph{non-symmetric}
length and orthogonality in the parabolic case
(Lemma~\ref{lem:orthogonal1}.\ref{item:par-orthogon}). The elliptic case seems to be free
from any such irregularities only because it is the standard model by
which the others are judged.

\begin{remark}
  \label{re:klein-descartes}
  We should say a word or two on proofs in this paper.  Majority of
  them are done through symbolic computations performed in the
  paper~\cite{Kisil05b} on the base of \GiNaC~\cite{GiNaC} computer
  algebra system. As a result we can reduce many proofs just to a
  one-line reference to the paper~\cite{Kisil05b}. In a sense this is
  the complete fulfilment of the \emph{Cartesian program} of reducing
  geometry to algebra with the latter to be done by straightforward
  \emph{mechanical calculations}.  Therefore the Erlangen program is
  nicely compatible with the Cartesian approach: the former defines
  the set of geometrical object with invariant properties and the
  latter provides a toolbox for their study. Another example of their
  unification in the field of non-commutative geometry was discussed
  in~\cite{Kisil02c}. 
\end{remark}

However the lack of intelligent proofs based on
smart arguments is undoubtedly a deficiency. An enlightening
reasoning (e.g. the proof of Lem.~\ref{le:invariance-of-cycles}) besides
establishing the correctness of a mathematical statement  gives
valuable insights about deep relations between objects. Thus it will
be worth to reestablish key results of this paper in a more synthetic
way. 

\subsection{The Paper Outline}
\label{sec:paper-outline}

Section~\ref{sec:ellipt-parab-hyperb}
describes the \(\SL\) group, its one-dimensional subgroups and
corresponding homogeneous spaces. Here corresponding Clifford algebras
show their relevance and cycles naturally appear as \(\SL\)-invariant
objects. 

To study cycles we extend in Section~\ref{sec:space-cycles} the
Schwerdtfeger--Fillmore--Springer--Cnops construction (SFSCc) to include parabolic
case. We also refine SFSCc from a traditional severe restriction that
space of cycles posses the same metric as the initial point space.
Cycles became points in a bigger space and got their presentation by
matrix. We derive first \(\SL\)-invariants of cycles from the classic
matrix invariants. 

Mutual disposition of two cycles may be also characterised through an
invariant notions of (normal and focal) orthogonalities, see
Section~\ref{sec:orth-invers-1}, both are defined in matrix terms of
SFSCc.  Orthogonality in generalised SFSCc is not anymore a
local property defined by tangents in the intersection point of cycles.
Moreover, the focal orthogonality is not even symmetric. The
corresponding notion of inversion (in a cycle) is considered as well.

Section~\ref{sec:metric-properties} describes distances and lengths
defined by cycles. Although they share some strange properties (e.g.
non-local character or non-symmetry) with the orthogonalities they are
legitimate objects in Erlangen approach since they are conformal under
the M\"obius maps. We also consider the corresponding perpendicularity
and its relation to orthogonality. Invariance of ``infinitesimal'' cycles and
corresponding version of conformality is considered in
Section~\ref{sec:invar-infin-scale}. 

Section~\ref{sec:global-properties} deals with the global properties
of the plane, e.g. its proper compactification by a zero-radius cycle at
infinity. Finally, Section~\ref{sec:unit-circles} considers some
aspects of the Cayley transform, with nicely interplays  with other notions
(e.g. focal orthogonality, lengths, etc.) considered in the previous sections.

To finish this introduction we point out the following natural question.
\begin{problem}
  \label{prob:higher-dim-gener}
  To which extend the subject presented here can be generalised to
  higher dimensions?
\end{problem}
 
\section{Elliptic, Parabolic and Hyperbolic Homogeneous Spaces}
\label{sec:ellipt-parab-hyperb}
We begin from representations of the \(\SL\) group in Clifford
algebras with two generators. They naturally introduce circles,
parabolas and hyperbolas as invariant objects of corresponding
geometries. 

\subsection[Special linear group and Clifford algebras]{$\SL$ group
  and Clifford Algebras}  
\label{sec:mobi-transf-spr}
We consider Clifford algebras defined by elliptic, parabolic and
hyperbolic bilinear forms. Then representations of \(\SL\) defined by
the same formula~\eqref{eq:moebius-def} will inherit this division. 
\begin{conv}
  \label{co:eph-and-a}
  There will be three different Clifford algebras \(\Cliff{e}\),
  \(\Cliff{p}\), \(\Cliff{h}\) corresponding to \emph{elliptic},
  \emph{parabolic}, and \emph{hyperbolic} cases respectively. The
  notation \(\Cliff{\sigma}\), with assumed values \(\sigma=-1\),
  \(0\), \(1\), refers to \emph{any} of these three algebras.
\end{conv}

A Clifford algebra \(\Cliff{\sigma }\) as a \(4\)-dimensional linear
space is spanned\footnote{We label generators of our Clifford algebra
  by \(e_0\) and \(e_1\) following the \texttt{C/C++} indexing
  agreement which is used by computer algebra calculations
  in~\cite{Kisil05b}.} by \(1\), \(e_0\), \(e_1\), \(e_0e_1\) with
\emph{non-commutative} multiplication defined by the following
identities\footnote{In
        light of usefulness of infinitesimal
        numbers~\cites{Devis77,Uspenskii88} in the parabolic spaces
        (see \S~\ref{sec:zero-length-cycles}) it may be worth to
        consider the parabolic Clifford algebra
        \(\Cliff{\varepsilon}\) with a generator
        \(e_1^2=\varepsilon\), where \(\varepsilon\) is an
        infinitesimal number.}:
\begin{equation}
  \label{eq:alg-mult}
 e_0^2=-1, \qquad
e_1^2=\sigma = \left\{
  \begin{array}{cl}
    -1, & \textrm{for \(\Cliff{e}\)---{elliptic} case}\\
    0, & \textrm{for \(\Cliff{p}\)---{parabolic} case}\\
    1, & \textrm{for \(\Cliff{h}\)---{hyperbolic} case}\\
  \end{array}
\right., \qquad 
e_0e_1=-e_1e_0.
\end{equation}
The two-dimensional subalgebra of \(\Cliff{e}\) spanned by \(1\) and
\(\rmi=e_1e_0=-e_0e_1\) is \emph{isomorphic} (and can actually replace
in all calculations!) the field of complex numbers \(\Space{C}{}\).
For example, from~\eqref{eq:alg-mult} follows that
\(\rmi^2=(e_1e_0)^2=-1\).  For any \(\Cliff{\sigma}\) we identify
\(\Space{R}{2}\) with the set of vectors \(w=ue_0+ve_1\), where
\((u,v)\in\Space{R}{2}\).  In the elliptic case of \(\Cliff{e}\) this
maps
\begin{equation}
  \label{eq:complexification}
  (u,v)\mapsto e_0(u+\rmi v)=e_0z, \qquad \textrm{ with } z=u+\rmi v,
\end{equation} 
in  the standard form of complex numbers. Similarly,
see~\cite{Kisil06a} and~\cite{Yaglom79}*{Supl.~C}  
\begin{enumerate}
\renewcommand{\theenumi}{(\ephname{enumi})}
\stepcounter{enumi}
\item in the parabolic case \(\varepsilon=e_1 e_0\) (such that
  \(\varepsilon^2=0\)) is known as \emph{dual unit} and all expressions
  \(u+\varepsilon v\), \(u, v\in\Space{R}{}\) form \emph{dual numbers}.
\item in the hyperbolic case \(e=e_1 e_0\) (such that \(e^2=1\)) is
  known as \emph{double unit} and all expressions \(u+e v\), \(u,
  v\in\Space{R}{}\) constitute \emph{double numbers}.
\end{enumerate}
\begin{remark}
  \label{re:dual-double}
  A part of this paper can be rewritten in terms of complex, dual
  and double numbers and it will have some common points with
  Supplement~C of the book~\cite{Yaglom79}. However the usage of
  Clifford algebras provides some facilities which do not have natural
  equivalent in complex numbers, see Rem.~\ref{re:dual-limit}. Moreover
  the language of Clifford algebras is more uniform and also allows
  straightforward generalisations to higher
  dimensions~\cite{Kisil96d}.
\end{remark}
We denote the space \(\Space{R}{2}\) of vectors \(u e_0+v e_1\) by
\(\Space{R}{e}\), \(\Space{R}{p}\) or \(\Space{R}{h}\) to highlight
which of Clifford algebras is used in the present context. The
notation \(\Space{R}{\sigma}\) assumes \(\Cliff{\sigma}\).

The \(\SL\) group~\cites{HoweTan92,Lang85,MTaylor86} consists of 
\(2\times 2\) matrices
\begin{displaymath}
  \begin{pmatrix}
    a&b\\c&d
  \end{pmatrix},\qquad  \textrm{ with } a, b, c, d\in\Space{R}{}
  \textrm{ and the determinant } ad-bc=1.
\end{displaymath}
An isomorphic realisation of \(\SL\) with the same multiplication is
obtained if we replace a matrix 
\(\begin{pmatrix} 
  a&b\\c&d
\end{pmatrix}
\) by \(\begin{pmatrix}
  a&b e_0\\-c e_0&d
\end{pmatrix}\) within any \(\Cliff{\sigma}\). 
The advantage of the latter form is that we can define the \emph{M\"obius 
transformation}  of \(\Space{R}{\sigma}\rightarrow \Space{R}{\sigma}\) for \emph{all}
three algebras \(\Cliff{\sigma}\) by the same expression:
\begin{equation}
  \label{eq:moebius-def}
  \begin{pmatrix}
    a&b e_0\\-c e_0&d
  \end{pmatrix}:\  ue_0+ve_1 \  \mapsto \ 
  \frac {
    a(ue_0+ve_1)+be_0
  } {
    -ce_0(ue_0+ve_1)+d
  },
\end{equation}
where the expression \(\frac{a}{b}\) in a non-commutative algebra is
always understood as \(ab^{-1}\), see~\cites{Cnops94a,Cnops02a}. Therefore
\(\frac{ac}{bc}=\frac{a}{b}\) but \(\frac{ca}{cb}\neq\frac{a}{b}\) in
general.

Again in the elliptic case the transformation~\eqref{eq:moebius-def}
is equivalent to, cf.~\citelist{\cite{Beardon05a}*{Ch.~13}
\cite{Beardon95}*{Ch.~3}}:
\begin{displaymath}
  \begin{pmatrix}
    a&b e_0\\-c e_0&d
  \end{pmatrix}:\  e_0z \  \mapsto \ 
  \frac {
    e_0(a(u+e_1e_0 v)+b)
  } {
    c(u+e_1e_0 v)+d
  }=   e_0\frac {
   az+b
  } {
    cz+d
  }, \textrm{ where } z=u+\rmi v,
\end{displaymath}
which is the standard form of a M\"obius transformation.
One can straightforwardly verify that the map~\eqref{eq:moebius-def}
is a left action of \(\SL\) on \(\Space{R}{\sigma}\), i.e. \(g_1(g_2
w)=(g_1g_2)w\).

To study finer structure of M\"obius transformations it is useful to
decompose an element \(g\) of \(\SL\) into the product \(g=g_a g_n g_k\): 
\begin{equation}
  \label{eq:iwasawa-decomp}
  \begin{pmatrix}
    a&b e_0\\-c e_0&d
  \end{pmatrix}=
  {\begin{pmatrix}
      \alpha^{-1} & 0\\0&\alpha
    \end{pmatrix}}
  {\begin{pmatrix}
      1&\nu e_0\\0&1
    \end{pmatrix}}
  {\begin{pmatrix}
      \cos\phi &  e_0\sin\phi\\ 
      e_0\sin\phi & \cos\phi 
    \end{pmatrix}},
\end{equation}
where the values of parameters are as follows:
\begin{equation}
  \label{eq:iwasawa-param}
  \alpha=\sqrt{c^2+d^2}, \qquad
  \nu=ac+bd,\qquad
  \phi = -\arctan \frac{c}{d}.
\end{equation} 
Consequently \(\cos \phi=\frac{d}{\sqrt{c^2+d^2}}\) and
\(\sin \phi=\frac{-c}{\sqrt{c^2+d^2}}\).  The
product~\eqref{eq:iwasawa-decomp} gives a realisation of the \emph{Iwasawa
  decomposition }~\cite{Lang85}*{\S~III.1} in the form \(\SL=ANK\),
where \(K\) is the maximal compact group, \(N\) is nilpotent and \(A\)
normalises \(N\).

\subsection{Actions of Subgroups} 
\label{sec:mobi-transf}
\begin{figure}[htbp]
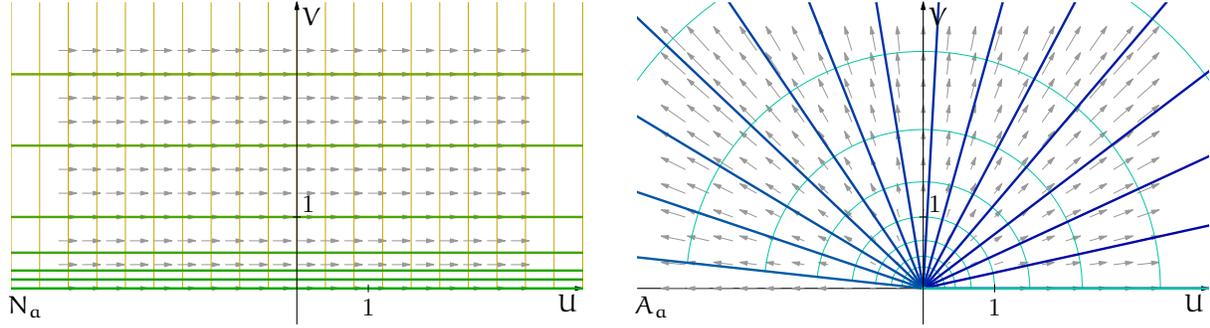

  \includegraphics[scale=.95]{parabolic0.7}\hfill
   \includegraphics[scale=.95]{parabolic0.8}
  \caption[Actions of the subgroups $A$ and $N$ by M\"obius
    transformations]{Actions of the subgroups \(A\) and \(N\) by M\"obius
    transformations.}
  \label{fig:a-n-action}
\end{figure}

We describe here orbits of the \emph{three} subgroups from the Iwasawa
decomposition~\eqref{eq:iwasawa-decomp} for all \emph{three} types of
Clifford algebras.  However there are \emph{less} than nine
(\(=3\times 3\)) different orbits since in all three EPH cases the
subgroups \(A\) and \(N\) act through M\"obius transformation
uniformly:
\begin{lemma} 
\label{le:an-action}
For any type of the Clifford algebra \(\Cliff{\sigma}\):
  \begin{enumerate}
  \item The subgroup \(N\) defines shifts \(ue_0+ve_1 \mapsto
    (u+\nu)e_0+ve_1 \) along the ``real'' axis \(U\) by
    \(\nu\).\\
    The vector field of the derived representation is \(dN_a{(u,v)}=(1,0)\).
  \item The subgroup \(A\) defines dilations \(ue_0+ve_1 \mapsto
    \alpha^{-2}(ue_0+ve_1)\) by the factor \(\alpha^{-2}\) which fixes origin
    \((0,0)\).\\
    The vector field of the derived representation is \(dA_a{(u,v)}=(2u,2v)\).
  \end{enumerate}
\end{lemma}
Orbits and vector fields corresponding to the \emph{derived
representation}~\citelist{\cite{Kirillov76}*{\S~6.3}
\cite{Lang85}*{Chap.~VI}} 
of the Lie algebra \(\algebra{sl}_2\) for subgroups \(A\) and \(N\)
are shown in Fig.~\ref{fig:a-n-action}. Thin transverse lines
join points of orbits corresponding to the same values of the
parameter along the subgroup. 

\begin{figure}[htbp]
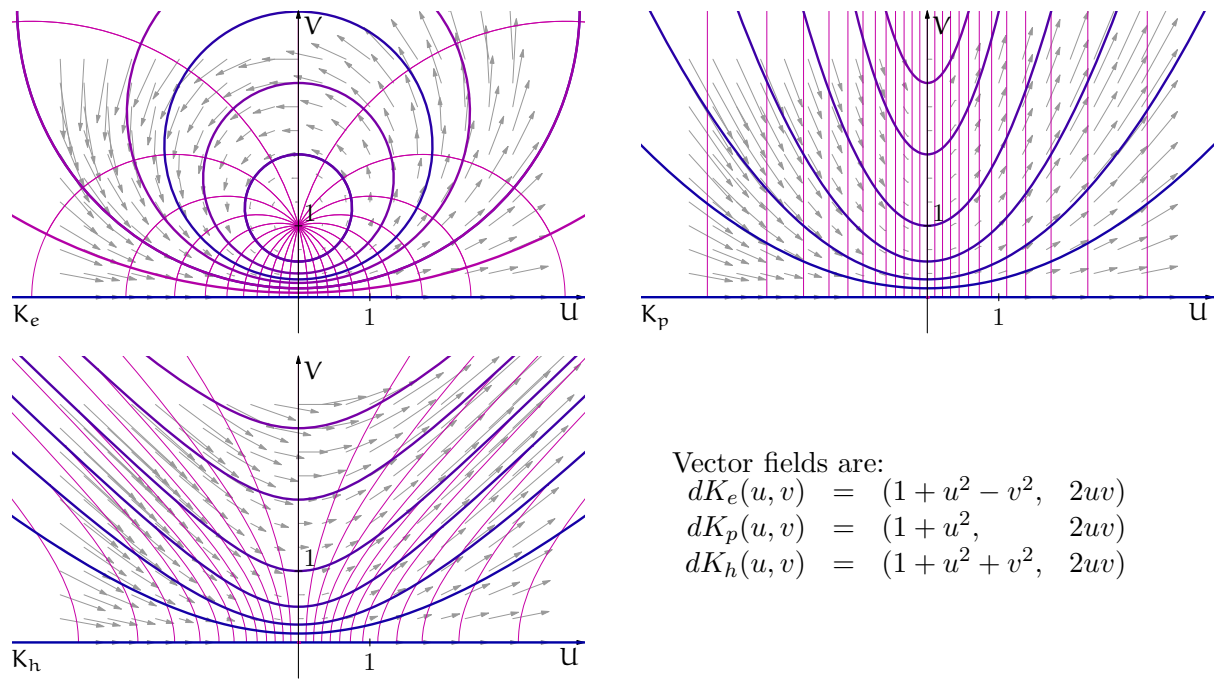

  \centering
\includegraphics[scale=.95
]{parabolic0.4}\hfill
\includegraphics[scale=.95
]{parabolic0.5}\\[2mm]
\parbox{.45\textwidth}{
\includegraphics[scale=.95
]{parabolic0.6}
}
\hfill \parbox{.45\textwidth}{
  Vector fields are:\\
  \(\begin{array}{rcll}
    dK_e(u,v)&=&(1+u^2-v^2,&2uv)\\
    dK_p(u,v)&=&(1+u^2,&2uv)\\
    dK_h(u,v)&=&(1+u^2+v^2,&2uv)\\
  \end{array}\)
}
  \caption[Action of the $K$ subgroup]{Action of the \(K\) subgroup. The corresponding orbits are
    circles, parabolas and hyperbolas.}
  \label{fig:k-subgroup}
\end{figure}

By contrast the actions of the subgroup  \(K\) look 
differently between the EPH cases, see Fig.~\ref{fig:k-subgroup}.
They obviously correlate with names chosen for
\(\Cliff{e}\), \(\Cliff{p}\), \(\Cliff{h}\). However algebraic
expressions for these orbits are uniform. 
\begin{lemma}
  \label{le:k-orbit-gen}
  A \(K\)-orbit in \(\Space{R}{\sigma}\) passing the point
  \((0,t)\) has the following equation: 
  \begin{equation}
    \label{eq:k-orbit-eq}
    (u^2-\sigma v^2)-2v\frac{t^{-1}-\sigma t }{2}+1=0, \qquad
    \text{where } \sigma = e_1^2 \ (\textrm{i.e } -1, 0 \textrm{ or } 1).
  \end{equation}
  The curvature of a \(K\)-orbit at point \((0,t)\) is equal to 
  \begin{displaymath}
    \kappa=\frac{2t}{1+\sigma t^2}.
  \end{displaymath}
\end{lemma}
A proof will be given later (see Ex.~\ref{ex:orbits}.\ref{it:k-orbit-similarity}),
when a more suitable tool will be in our disposal. Meanwhile these
formulae allows to produce geometric characterisation of \(K\)-orbits.
\begin{lemma} 
  \label{le:k-action}
  \begin{enumerate}
    \renewcommand{\theenumi}{(\ephname{enumi})}
  \item \label{item:circle-desc} 
    For \(\Cliff{e}\) the orbits of \(K\) are circles, they are
    coaxal~\cite{CoxeterGreitzer}*{\S~2.3} with the real line being
    the radical axis. A circle with centre
     at \((0, (v+v^{-1})/2)\) passing through two points \((0,v)\)
    and \((0,v^{-1})\).\\
    The vector field of the derived representation is
    \(dK_e{(u,v)}=(u^2-v^2+1, 2uv)\). 
  \item \label{item:parab-desc}
    For \(\Cliff{p}\) the orbits of \(K\) are parabolas with the
    vertical axis \(V\). A parabola passing through \((0,v/2)\) has horizontal
    directrix passing through \((0, v-v^{-1}/2)\) 
    and focus at \((0,(v+v^{-1})/2)\). \\
    The vector field of the derived representation is
    \(dK_p{(u,v)}=(u^2+1, 2uv)\).
  \item \label{item:hyperb-desc}
    For \(\Cliff{h}\) the orbits of \(K\) are hyperbolas with
    asymptotes parallel to lines \(u=\pm v\). A hyperbola passing
    through the point \((0,v)\) has the focal distance
    \(2p\), where \(p=\frac{v^2+1}{\sqrt{2}v}\) and the upper focus
    is located at \((0,f)\) with:
    \begin{displaymath}
    f=\left\{\begin{array}{ll}
      p-\sqrt{\frac{p^2}{2}-1}, &\textrm{ for } 0<v<1; \textrm{
        and }\\
      p+\sqrt{\frac{p^2}{2}-1}, & \textrm{ for } v\geq 1.
    \end{array}\right.
  \end{displaymath}
  The vector field of the derived representation is
    \(dK_h{(u,v)}=(u^2+v^2+1, 2uv)\).
  \end{enumerate}
\end{lemma} 
\begin{figure}[htbp]
  \centering
  (a)\includegraphics[scale=1.2]{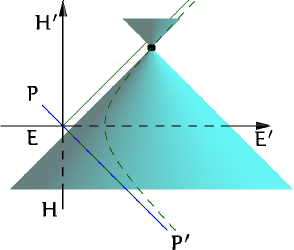}\hfill
  (b)\includegraphics[scale=1.2]{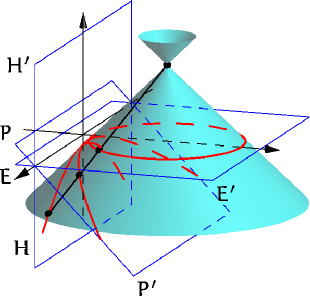}\hfill
  \caption[$K$-orbits as conic sections]{\(K\)-orbits as conic
    sections:\\
    (a) a flat projection along \(U\) axis;\\ 
    (b) same values of \(\phi\) on different orbits belong to the same
    generator of the cone.}
  \label{fig:k-orbit-sect}
\end{figure}
Since all \(K\)-orbits are conic sections it is tempting to obtain
them as sections of some cones. To this end we define the family of
double-sided right-angle cones be parametrised by \(t>0\): 
\begin{equation}
  \label{eq:cones-fam}
  x^2+(y-\half(t+t^{-1}))^2-(z-\half(t-t^{-1}))^2=0.
\end{equation}
The vertices of cones belong to the hyperbola \(\{x=0,
y^2-z^2=1\}\), see Fig.~\ref{fig:k-orbit-sect} for illustration.
\begin{lemma} \(K\)-orbits may be obtained cases as follows: 
  \begin{enumerate}
    \renewcommand{\theenumi}{(\ephname{enumi})}
  \item elliptic \(K\)-orbits are sections of cones~\eqref{eq:cones-fam} by the plane
    \(z=0\) (\(EE'\) on Fig.~\ref{fig:k-orbit-sect}); 
  \item parabolic \(K\)-orbits are sections of~\eqref{eq:cones-fam} by the plane
    \(y=\pm z\) (\(PP'\) on Fig.~\ref{fig:k-orbit-sect});
  \item hyperbolic   \(K\)-orbits are sections of~\eqref{eq:cones-fam} by the plane \(y=0\) 
    (\(HH'\) on Fig.~\ref{fig:k-orbit-sect});
  \end{enumerate}

  Moreover, each straight line generating a cone from the
  family~\eqref{eq:cones-fam} is crossing corresponding elliptic,
  parabolic and hyperbolic \(K\)-orbits at points with the same value of
  parameter \(\phi\)~\eqref{eq:iwasawa-param} of the subgroup \(K\).
\end{lemma}
From the above algebraic and geometric descriptions of the orbits we
can make several observations.
\begin{remark}  
  \label{rem:orbits2}
  \begin{enumerate}
  \item The values of all three vector fields \(dK_e\), \(dK_p\) and \(dK_h\)
    coincide on the ``real'' \(U\)-axis \(v=0\), i.e. they are three
    different extensions into the domain of the same boundary
    condition. Another source of this: the axis \(U\) is the
    intersection of planes \(EE'\), \(PP'\) and \(HH'\) on
    Fig.~\ref{fig:k-orbit-sect}. 
  \item \label{it:hyp-orbit} The hyperbola passing through the point
    \((0,1)\) has the shortest focal length \(\sqrt{2}\) among all
    other hyperbolic orbits since it is the section of the cone
    \(x^2+(y-1)^2+z^2=0\) closest from the family to the plane
    \(HH'\).   
  \item Two hyperbolas passing through \((0,v)\) and \((0,v^{-1})\)
    have the same focal length since they are sections of two cones
    with the same distance from \(HH'\). Moreover, two such hyperbolas
    in the lower- and upper half-planes passing the points
    \((0,v)\) and \((0,-v^{-1})\) are sections of the same
    double-sided cone.  They are related to each other as explained in
    Remark~\ref{rem:two-cver}.\ref{it:hyp-object}.
\end{enumerate}
\end{remark}
One can see from the first picture in Fig.~\ref{fig:k-subgroup} that the
elliptic action of subgroup \(K\) fixes the point \(e_1\). More
generally we have:
\begin{figure}[htbp]
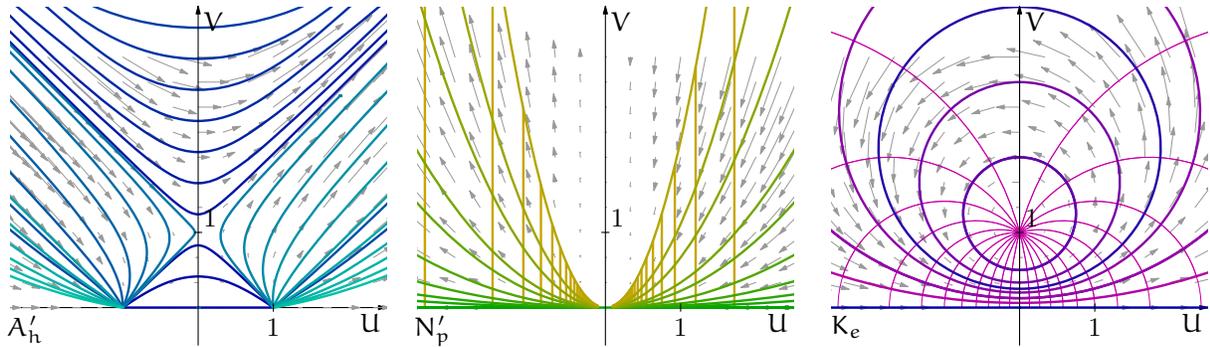

  \includegraphics[scale=1]{parabolic0.17}\hfill
   \includegraphics[scale=1]{parabolic0.18}
   \hfill
   \includegraphics[scale=1]{parabolic0.19}
  \caption[Actions of fix-subgroups]{Actions of the subgroups which
    fix point \(e_1\) in three cases.} 
  \label{fig:fix-sbroups}
\end{figure}
\begin{lemma} 
  \label{le:fix-subgroups}
  The fix group of the point \(e_1\) is
  \begin{enumerate}
    \renewcommand{\theenumi}{(\ephname{enumi})}
  \item\label{it:fix-group-ell} 
    the subgroup \(K'_e=K\) in the elliptic case. Thus the
    elliptic upper half-plane is a model for the homogeneous space
    \(\SL/K\); 
  \item \label{it:fix-group-par}
    the subgroup \(N'_p\) of matrices
    \begin{equation}
      \label{eq:par-fix-subgroup}
      \begin{pmatrix}
        1&0\\
         \nu e_0 & 1
      \end{pmatrix}
      =
      \begin{pmatrix}
        0&e_0\\
        e_0 & 0
      \end{pmatrix}
      \begin{pmatrix}
        1&\nu e_0\\
         0& 1
      \end{pmatrix}
      \begin{pmatrix}
        0&-e_0\\
        -e_0 & 0
      \end{pmatrix}
    \end{equation}
    in the parabolic case.  It also fixes any point \(ve1\).  It is
    conjugate to subgroup \(N\), thus the parabolic upper half-plane is
    a model for the homogeneous space \(\SL/N\);
  \item \label{it:fix-group-hyp}
    the subgroup \(A'_h\) of matrices
    \begin{equation}
      \label{eq:hyp-fix-subgroup}
      \begin{pmatrix}
        \cosh(\tau) & \sinh(\tau)e_0\\
        -\sinh(\tau)e_0 & \cosh(\tau) 
      \end{pmatrix} =
      \frac{1}{2}
      \begin{pmatrix}
        1 & -e_0\\
        -e_0 & 1
      \end{pmatrix}
      \begin{pmatrix}
        e^\tau & 0\\
        0 & e^{-\tau}
      \end{pmatrix}
      \begin{pmatrix}
        1 & e_0\\
        e_0 & 1
      \end{pmatrix},
    \end{equation}
    in the hyperbolic case. It is conjugate to subgroup \(A\), thus
    two copies of the upper halfplane (see
    Section~\ref{sec:invar-upper-half}) is a model for \(\SL/A\).
  \end{enumerate}

  Moreover, vectors fields of these actions are \((u^2+\sigma(v^2-1), 2uv)\) for
  the corresponding values of \(\sigma\).
  Orbits of the fix groups satisfy to the equation:
  \begin{displaymath}
        (u^2-\sigma v^2)-2lv-\sigma=0, \qquad
    \text{where } l\in\Space{R}{}.
  \end{displaymath}
\end{lemma}
\begin{remark}
  \begin{enumerate}
  \item   Note that we can uniformly express the fix-subgroups of \(e_1\) in
    all EPH cases by
    matrices of the form:
    \begin{displaymath}
      \begin{pmatrix}
        a & -\sigma e_0b\\
        -e_0b & a
      \end{pmatrix},\qquad \text{ where } a^2-\sigma b^2=1.
    \end{displaymath}
  \item In the hyperbolic case the subgroup \(A'_h\) may be extended to a
    subgroup \(A''_h\) by the element \(
    \begin{pmatrix}
      0&e_0\\e_0&0
    \end{pmatrix}\), which flips upper and lower half-planes (see
    Section~\ref{sec:invar-upper-half}).  The subgroup \(A''_h\) fixes
    the set \(\{e_1,-e_1\}\).
  \end{enumerate}
\end{remark}

\begin{lemma}
  \label{le:fix-ax+b-gen-sl2}
  M\"obius action of \(\SL\) in each EPH case is generated by
  action  the corresponding fix-subgroup (\(A''_h\) in the hyperbolic case) and
  actions of the \(ax+b\) group, e.g. subgroups \(A\) and \(N\).
\end{lemma}
\begin{proof}
  The \(ax+b\) group transitively acts on the upper or lower half-plane. Thus
  for any \(g\in\SL\) there is \(h\) in \(ax+b\) group such that
  \(h^{-1}g\) either fixes \(e_1\) or sends it to \(-e_1\). Thus
  \(h^{-1}g\) is in the corresponding fix-group.
\end{proof}

\subsection{Invariance of Cycles}
\label{sec:invariance-cycles}

As we will see soon the three types of \(K\)-orbits are principal
invariants of the constructed geometries, thus we will unify them in
the following definition.
\begin{definition}
  We use the word \emph{cycle} to denote loci in \(\Space{R}{\sigma}\)
  defined by the equation:
\begin{subequations}
  \label{eq:cycle-def}
  \begin{equation}
    \label{eq:cycle-def-cliff}
    -k(u e_0+v e_1)^2-2\scalar{(l,n)}{(u,v)}+m=0
  \end{equation}
  or equivalently (avoiding any reference to Clifford algebra generators): 
  \begin{equation}
    \label{eq:cycle-def-real}
    k(u^2-\sigma v^2)-2lu-2nv+m=0, \qquad \text{ where } \sigma = e_1^2,
  \end{equation}
  or equivalently (using \emph{only} Clifford algebra operations,
  cf.~\cite{Yaglom79}*{Supl.~C(42a)}):  
  \begin{equation}
    \label{eq:cycle-def-compl}
    K w^2+Lw-wL+M=0,
  \end{equation}
  where \(w= ue_0+ve_1\), \(K=-ke_{01}\), \(L=-n e_0+l e_1\),
  \(M=me_{01}\).
\end{subequations}
\end{definition}
Such cycles obviously mean for certain \(k\), \(l\), \(n\), \(m\)
straight lines \textbf{and one of the following}:
\begin{enumerate}
  \renewcommand{\theenumi}{(\ephname{enumi})}
\item in the elliptic case: circles with centre
  \(
  \left(\frac{l}{k},\frac{n}{k}\right)\) and squared radius
  \(
  m-\frac{l^2+n^2}{k}\);
\item in the parabolic case: parabolas with horizontal directrix and
  focus at \(
  \left(\frac{l}{k},
    \frac{m}{2n}-\frac{l^2}{2nk}+\frac{n}{2k}\right)\); 
\item  in the hyperbolic case: rectangular hyperbolas with centre
  \(
  \left(\frac{l}{k},-\frac{n}{k}\right)\) and a
  vertical axis of symmetry.  
\end{enumerate}
Moreover words \emph{parabola} and \emph{hyperbola} in this paper
always assume only the above described types. Straight lines
are also called \emph{flat cycles}.

All three EPH types of cycles are enjoying many common properties,
sometimes even beyond that we normally expect. For example, the
following definition is quite intelligible even when extended from the
above elliptic and hyperbolic cases to the parabolic one.
\begin{definition}
  \label{de:center-first}
  \emph{\(\rs\)-Centre} of the \(\sigma\)-cycle~\eqref{eq:cycle-def}
  for any EPH case is the point \(\left(\frac{l}{k},
    -\rs\frac{n}{k}\right)\in\Space{R}{\sigma}\). Notions of e-centre,
  p-centre, h-centre are used along the adopted EPH notations.
  
  Centres of straight lines are at infinity, see  subsection~\ref{sec:comp-spra}.
\end{definition}
\begin{remark}
  \label{re:different-signatures}
  Here we use a signature \(\rs=-1\), \(0\) or \(1\) of a  Clifford
  algebra which is not related to the signature \(\sigma\) of the
  space \(\Space{R}{\sigma}\). We will need also a third signature
  \(\bs\) to describe the geometry of cycles in
  Defn.~\ref{de:cycle-2-matrix}. 
\end{remark}
The meaningfulness of this definition even in the parabolic case is
justified, for example, by:
\begin{itemize}
\item the uniformity of description of relations between centres of
  orthogonal cycles, see the next subsection and
  Fig.~\ref{fig:orthogonality1}. 
\item the appearance of \emph{concentric parabolas} in
  Fig.~\ref{fig:unit-disks}(\(N_{P_{e}}\)) and~(\(N_{P_{h}}\)).
\end{itemize}

Using the Lemmas~\ref{le:an-action} and~\ref{le:k-action} we can give
an easy (and virtually calculation-free!) proof of invariance for
corresponding cycles.  

\begin{lemma} \label{le:invariance-of-cycles}
  M\"obius transformations preserve the cycles in the
  upper   half-plane, i.e.:
  \begin{enumerate}
    \renewcommand{\theenumi}{(\ephname{enumi})}
  \item For \(\Cliff{e}\) M\"obius transformations map circles to circles.
  \item For \(\Cliff{p}\) M\"obius transformations map parabolas to parabolas.
  \item For \(\Cliff{h}\) M\"obius transformations map hyperbolas to hyperbolas. 
  \end{enumerate}
\end{lemma}
\begin{proof}
\begin{figure}[htbp]
  \centering
  \includegraphics[scale=1.1]{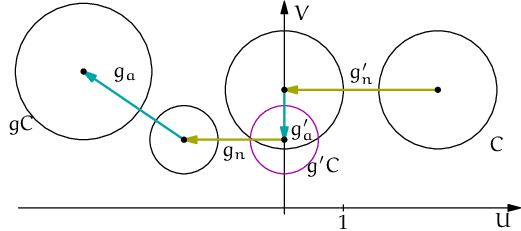}
  \caption[The decomposition of an arbitrary M\"obius
  transformation]{Decomposition of an arbitrary M\"obius
    transformation \(g\) into a product \(g=g_a g_n g_k g_a' g_n'\).}
  \label{fig:moeb-decomp}
\end{figure}
  Our first observation is that the subgroups \(A\) and \(N\) obviously
  preserve all circles, parabolas, hyperbolas and straight lines in all
  \(\Cliff{\sigma}\). Thus we use subgroups \(A\) and \(N\) to  fit a given
  cycle exactly on a particular orbit of subgroup \(K\) shown on
  Fig.~\ref{fig:k-subgroup} of the corresponding type. 

  To this end for an arbitrary cycle \(S\) we can find \(g_n'\in N\)
  which puts centre of \(S\) on the \(V\)-axis, see
  Fig.~\ref{fig:moeb-decomp}. Then there is a unique \(g_a'\in A\)
  which scales it exactly to an orbit of \(K\), e.g. for a
  circle passing through points \((0,v_1)\) and \((0,v_2)\) the scaling
  factor is \(\frac{1}{\sqrt{v_1v_2}}\) according to
  Lemma~\ref{le:k-action}.\ref{item:circle-desc}. Let \(g'=g_a' g_n'\), then for any
  element \(g\in \SL\) using the Iwasawa decomposition of \(g
  g'^{-1}=g_a g_n g_k\) we get the presentation \(g=g_a g_n g_k g_a'
  g_n'\) with \(g_a, g_a'\in A\), \(g_n, g_n'\in N\) and \(g_k \in
  K\).

  Then the image \(g'S\) of the cycle \(S\) under \(g'=g_a' g_n'\)
  is a cycle itself in the obvious way, then \(g_k (g'S)\) is again a
  cycle since \(g'S\) was arranged to coincide with a \(K\)-orbit, and
  finally  \(gS=g_a g_n (g_k (g'S))\) is a cycle due to the obvious
  action of \(g_a g_n\), see Fig.~\ref{fig:moeb-decomp} for an
  illustration. 
\end{proof}
One can naturally wish that all other proofs in this paper will be of
the same sort. This is likely to be possible, however we use a lot of
computer algebra calculations as well.

\section{Space of Cycles}
\label{sec:space-cycles}

We saw in the previous sections that cycles are M\"obius invariant,
thus they are natural objects of the corresponding geometries in the
sense of F.~Klein. An efficient tool of their study is to represent
all cycles in \(\Space{R}{\sigma}\) by points of a new bigger space.

\subsection{Schwerdtfeger--Fillmore--Springer--Cnops Construction (SFSCc)}
\label{sec:fillm-spring-cnops}

It is well known that linear-fractional transformations can be
linearised by a transition into a suitable projective
space~\cite{Olver99}*{Cha.~1}.  The fundamental idea of the
Schwerdtfeger--Fill\-more--Springer--Cnops construction
(SFSCc)~\cites{Schwerdtfeger79a,Cnops02a,Porteous95,FillmoreSpringer90a,Kirillov06,Kisil12a} is that for linearisation of
M\"obius transformation in \(\Space{R}{\sigma}\) the required
projective space can be identified with the \emph{space of all cycles}
in \(\Space{R}{\sigma}\).  The latter can be associated with certain
subset of \(2\times 2\) matrices.  SFSCc can be adopted
from~\cites{Cnops02a,Porteous95} to serve all three EPH cases with some
interesting modifications.
\begin{definition}
  \label{de:cycle-2-matrix}  
  Let \(\Space{P}{3}\) be the projective space, i.e. collection of
  the rays passing through points in \(\Space{R}{4}\). We define the
  following two identifications (depending from some additional
  parameters \(\sigma\), \(\bs\) and \(s\) described below) which map a
  point \((k, l, n, m)\in \Space{P}{3}\) to:
  \begin{itemize}
  \item[\(Q\):]
    the cycle (quadric) \(\cycle{}{}\) on
    \(\Space{R}{\sigma}\)
    defined by the equations~\eqref{eq:cycle-def} with constant
    parameters \(k\), \(l\), \(n\), \(m\):
    \begin{equation}
      \label{eq:cycle-def-1}
          -k(e_0u+e_1v)^2-2\scalar{(l,n)}{(u,v)}+m=0,
    \end{equation}
     for some \(\Cliff{\sigma}\) with generators \(e_0^2=-1\), \(e_1^2=\sigma\).
  \item[\(M\):]
    the ray of \(2\times 2\) matrices 
    passing through 
    \begin{equation}
      \label{eq:matrix-for-cycle}
      \cycle{s}{\bs} = \begin{pmatrix}
      l\se{0}+sn\se{1} & m\\
      k & -l\se{0}-sn\se{1}
    \end{pmatrix} \in M_2(\Cliff{\bs}), \text{ with }
     \se{0}^2=-1, \se{1}^2=\bs,
    \end{equation}
    i.e. generators \(\se{0}\) and \(\se{1}\) of \(\Cliff{\bs}\)
    can be \emph{of any type}: elliptic, parabolic or hyperbolic
    regardless of the \(\Cliff{\sigma}\) in~\eqref{eq:cycle-def-1}.
  \end{itemize}
   The meaningful values of parameters \(\sigma\), \(\bs\) and \(s\)
   are \(-1\), \(0\) or \(1\), and in many cases \(s\) is equal to
   \(\sigma\).
\end{definition}
\begin{remark}
  A hint for the composition of the matrix~\eqref{eq:matrix-for-cycle}
  is provided by the following identity:
  \begin{displaymath}
    \begin{pmatrix} 
      1&w
    \end{pmatrix}
    \begin{pmatrix} 
      L&M\\K&-L
    \end{pmatrix}
    \begin{pmatrix} 
      w\\1
    \end{pmatrix}=     w K w+Lw-wL+M,
  \end{displaymath}
  which realises the equation~\eqref{eq:cycle-def-compl} of a cycle.
\end{remark}

The both identifications \(Q\) and \(M\) are straightforward. Indeed,
a point \((k, l, n, m)\in\Space{P}{3}\) equally well represents (as
soon as \(\sigma \), \(\bs\) and \(s\) are already fixed) both the
equation~\eqref{eq:cycle-def-1} and the ray of
matrix~\eqref{eq:matrix-for-cycle}. Thus for fixed \(\sigma \),
\(\bs\) and \(s\) one can introduce the correspondence between
quadrics and matrices shown by the horizontal arrow on the following
diagram:
\begin{equation}
  \label{eq:cycles-diagram}
\xymatrix{ 
    &\Space{P}{3} \ar@{<->}[dl]_{Q} \ar@{<->}[dr]^{M}&\\
    \text{Quadrics on }\Space{R}{\sigma} \ar@{<->}[rr]^(.55){Q\circ M} &&{M_2(\Cliff{\bs})}}
\end{equation}
which combines \(Q\) and \(M\). On the first glance the dotted arrow
seems to be of a little practical interest since it depends from too
many different parameters (\(\sigma\), \(\bs\) and
\(s\)). However the following result demonstrates that it is 
compatible with easy calculations of images of cycles under the
M\"obius transformations.
\begin{proposition} 
  \label{pr:transform-cycles} 
  A cycle \(-k(e_0u+e_1v)^2-2\scalar{(l,n)}{(u,v)}+m=0\) is
  transformed by \(g\in\SL\) into the cycle
  \(-\tilde{k}(e_0u+e_1v)^2-2\scalar{(\tilde{l},
    \tilde{n})}{(u,v)}+\tilde{m}=0\) such that
  \begin{equation}
    \label{eq:cycle-transform-short}
    \cycle[\tilde]{s}{\bs}=  g\cycle{s}{\bs}g^{-1}
  \end{equation}
  for any Clifford algebras \(\Cliff{\sigma}\) and
  \(\Cliff{\bs}\). Explicitly this means:
  \begin{eqnarray}
    \label{eq:cycle-transform}
    \lefteqn{
    \begin{pmatrix}
      \tilde{l}\se{0}+s\tilde{n}\se{1} & \tilde{m}\\
      \tilde{k} & -\tilde{l}\se{0}-s\tilde{n}\se{1}
    \end{pmatrix}
  }\qquad\qquad \\
    &=&
    \begin{pmatrix}
      a & b\se{0}\\
      -c\se{0} & d
    \end{pmatrix}
    \begin{pmatrix}
      l\se{0}+sn\se{1} & m\\
      k & -l\se{0}-sn\se{1}
    \end{pmatrix}
    \begin{pmatrix}
      d & -b\se{0}\\
      c\se{0} & a
    \end{pmatrix}. \nonumber 
  \end{eqnarray}
\end{proposition}
\begin{proof}
  It is already established in the elliptic and hyperbolic cases for
  \(\sigma=\bs\), see \cite{Cnops02a}. For all EPH cases (including
  parabolic) it can be done by the direct calculation in
  \GiNaC~\cite{Kisil05b}*{\S~\ref{G-sec:mobi-invar-cycl}}. An
  alternative idea of an elegant proof based on the zero-radius cycles
  and orthogonality (see below) may be borrowed from~\cite{Cnops02a}.
\end{proof} 
\begin{example}
  \label{ex:orbits}
  \begin{enumerate}
  \item   
    \label{it:real-line-inv}
    The real axis \(v=0\) is represented by the ray coming through
  \((0,0,1,0)\) and a matrix \(
  \begin{pmatrix}
    s\se{1}&0\\
    0&-s\se{1}
  \end{pmatrix}\). For any \(
  \begin{pmatrix}
    a&b\se{0}\\
    -c\se{0}& d
  \end{pmatrix}\in \SL\) we have:
  \begin{displaymath}
    \begin{pmatrix}
      a&b\se{0}\\
      -c\se{0}& d
    \end{pmatrix}
    \begin{pmatrix}
      s\se{1}&0\\
      0&-s\se{1}
    \end{pmatrix}
    \begin{pmatrix}
      d & -b\se{0}\\
      c\se{0} & a
    \end{pmatrix} =
    \begin{pmatrix}
      s\se{1}&0\\
      0&-s\se{1}
    \end{pmatrix},
  \end{displaymath}
  i.e. the real line is \(\SL\)-invariant.
\item 
  \label{it:k-orbit-similarity}
  A direct calculation in
  GiNaC~\cite{Kisil05b}*{\S~\ref{G-sec:transf-k-orbits}} shows that
  matrices representing cycles from~\eqref{eq:k-orbit-eq} are
  invariant under the similarity with elements of \(K\), thus they are
  indeed \(K\)-orbits.
  \end{enumerate}
\end{example}

It is surprising on the first glance that the \(\cycle{s}{\bs}\) is
defined through a Clifford algebra \(\Cliff{\bs}\) with an arbitrary
sign of \(\se{1}^2\). However a moment of reflections reveals that
transformation~\eqref{eq:cycle-transform} depends only from the sign
of \(\se{0}^2\) but does not involve any quadratic (or higher) terms
of \(\se{1}\). 

\begin{remark}
  \label{re:real-line-inv-0}
  Such a variety of choices is a consequence of the usage of \(\SL\)---a
  smaller group of symmetries in comparison to the all M\"obius maps
  of \(\Space{R}{2}\). The \(\SL\) group fixes the real line and consequently a
  decomposition of vectors into ``real'' (\(e_0\)) and ``imaginary''
  (\(e_1\)) parts is obvious. This permits to assign an arbitrary
  value to the square of the ``imaginary unit'' \(e_1e_0\).

  Geometric invariants defined below, e.g. orthogonalities in
  sections~\ref{sec:orthogonality} and~\ref{sec:focal-orthogonality},
  demonstrate ``awareness'' of the real line invariance in one way or
  another. We will call this the \emph{boundary effect} in the
  upper half-plane geometry. The famous question
  \href{http://en.wikipedia.org/wiki/Hearing_the_shape_of_a_drum}{on
    hearing drum's shape} has a sister:
  \begin{quote}
  \emph{Can we see/feel the boundary from inside a domain?}
  \end{quote}
  Rems.~\ref{re:real-line-inv-1},~\ref{re:real-line-inv-2}
  and~\ref{re:real-line-inv-3} provide hints for positive answers.
\end{remark} 

\begin{figure}[htbp]
  \centering
  (a) \includegraphics[scale=1.3]{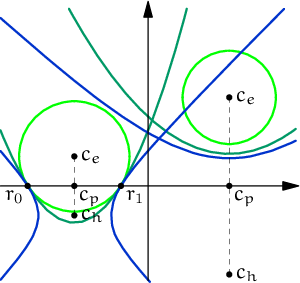}\hfill
  (b) \includegraphics[scale=1.3]{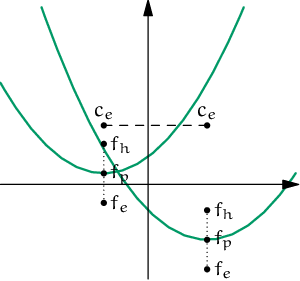}
  \caption[Cycle implementations, centres and foci]{(a) Different
    EPH implementations of the same cycles defined by quadruples of
    numbers.\\
  (b) Centres and foci of two parabolas with the same focal length.} 
  \label{fig:eph-cycle}
\end{figure}
To encompass all aspects from~\eqref{eq:cycles-diagram} we think
a cycle \(\cycle{s}{\bs}\) defined by a quadruple \((k, l, n, m)\) as an
``imageless'' object which have distinct implementations (a circle, a
parabola or a hyperbola) in the corresponding space \(\Space{R}{\sigma}\).
These implementations may look very different, see
Fig.~\ref{fig:eph-cycle}(a), but still have some properties in common.
For example,
\begin{itemize}
\item All implementations have the same vertical axis of symmetries;
\item Intersections with the real axis (if exist) coincide, see
  \(r_1\) and \(r_2\) for the left cycle in
  Fig.~\ref{fig:eph-cycle}(a).
\item Centres of circle \(c_e\) and corresponding hyperbolas \(c_h\)
  are mirror reflections of each other in the real axis with the
  parabolic centre be in the middle point.
\end{itemize}
Lemma~\ref{le:k-action} gives another example of similarities between 
different implementations of the same cycles defined by the
equation~\eqref{eq:k-orbit-eq}. 

Finally, we may restate the Prop.~\ref{pr:transform-cycles} as an
intertwining property.
\begin{corollary}
  Any implementation of cycles shown on~\eqref{eq:cycles-diagram} by
  the dotted arrow for any combination of \(\sigma\),
  \(\bs\) and \(s\) intertwines two actions of \(\SL\): by
  matrix conjugation~\eqref{eq:cycle-transform-short} and M\"obius
  transformations~\eqref{eq:moebius-def}.
\end{corollary}
\begin{remark}
  A similar representation of circles by \(2\times 2\) complex
  matrices which intertwines M\"obius transformations and matrix
  conjugations was used recently by A.A.~Kirillov~\cite{Kirillov06} in
  the study of the Apollonian gasket. Kirillov's matrix
  realisation~\cite{Kirillov06} of a cycle has an attractive
  ``self-adjoint'' form:
  \begin{equation}
    \label{eq:kirrilov-matrix}
    \cycle{s}{\bs} = \begin{pmatrix}
      m & l\se{0}+sn\se{1} \\
      -l\se{0}-sn\se{1} & k  
    \end{pmatrix} \text{ (in notations of this paper)}.
  \end{equation}
  Note that the matrix inverse to~\eqref{eq:kirrilov-matrix}  is intertwined with the SFSCc
  presentation~\eqref{eq:matrix-for-cycle} by the matrix \(
  \begin{pmatrix}
    0&1\\1&0
  \end{pmatrix}\).
\end{remark}
\subsection{First Invariants of Cycles}
\label{sec:invariants-cycles}

Using implementations from Definition~\ref{de:cycle-2-matrix} and
relation~\eqref{eq:cycle-transform-short} we can derive some
invariants of cycles (under the M\"obius transformations) from
well-known invariants of matrix (under similarities). First we use
trace to define an \emph{invariant} inner product in the space of
cycles. 

\begin{definition}
  \label{de:inner-product}
  \emph{Inner \(\bs\)-product} of two cycles is given by the trace of their
  product as matrices:
  \begin{equation}
    \label{eq:inner-product}
    \scalar{\cycle{s}{\bs}}{\cycle[\tilde]{s}{\bs}} = \tr(\cycle{s}{\bs}\cycle[\tilde]{s}{\bs}).
  \end{equation}
\end{definition}
The above definition is very similar to an inner product defined in
operator algebras~\cite{Arveson76}. This is not a coincidence: cycles
act on points of \(\Space{R}{\sigma}\) by inversions, see
subsection~\ref{sec:invers-in-cycl}, and this action is linearised
by SFSCc, thus cycles can be viewed as linear operators as well. 

Geo\emph{metrical} interpretation of the inner product will be given
in Cor.~\ref{co:inner-product-dist}.

An obvious but interesting observation is that for matrices
representing cycles we obtain the second classical invariant
(determinant) under similarities~\eqref{eq:cycle-transform-short} from
the first (trace) as follows: 
  \begin{equation}
    \label{eq:trace-det}
    \scalar{\cycle{s}{\bs}}{\cycle{s}{\bs}} = -2 \det\cycle{s}{\bs}.
  \end{equation}
The explicit expression for the determinant is:
\begin{equation}
  \label{eq:determinant-def}
  \det\cycle{s}{\bs} = l^2-\bs s^2 n^2 -mk
\end{equation}

We recall that the same cycle is defined by any matrix \(\lambda
\cycle{s}{\bs}\), \(\lambda \in \Space[+]{R}{}\), thus the
determinant, even being M\"obius-invariant, is useful only in the
identities of the sort \(\det\cycle{s}{\bs}=0\).  Note also that
\(\tr(\cycle{s}{\bs})=0\) for any matrix of the
form~\eqref{eq:matrix-for-cycle}. Since it may be convenient to have a
predefined representative of a cycle out of the ray of equivalent
SFSCc matrices we introduce the following normalisation.
\begin{definition}
  \label{de:normalisation}
  A SFSCc matrix representing a cycle is said to be \(k\)-normalised if its
  \((2,1)\)-element  is \(1\) and it is \(\det\)-normalised if its
  determinant is equal \(1\).
\end{definition}
Each normalisation has its own advantages: element \((1,1)\) of
\(k\)-normalised matrix immediately tell us the centre of the cycle, meanwhile
\(\det\)-normalisation is preserved by matrix conjugation with \(\SL\)
element (which is important in view of
Prop.~\ref{pr:transform-cycles}). 
The later normalisation is used, for example, in~\cite{Kirillov06}  

Taking into account its invariance it is not surprising that the
determinant of a cycle enters the following definition~\ref{de:focus}
of the focus and the invariant zero-radius
cycles from Def.~\ref{de:zero-radius-cycle}.
\begin{definition}
  \label{de:focus}
  \emph{\(\rs\)-Focus} of a cycle \(\cycle{s}{\bs}\)
  is the point in \(\Space{R}{\sigma}\)
  \begin{equation}
    \label{eq:def-focus}
    f_{\rs}=\left(\frac{l}{k}, -\frac{\det\cycle{s}{\rs}}{2nk}\right)
    \quad \text{ or  explicitly } \quad
    f_{\rs} =\left(\frac{l}{k}, \frac{mk-l^2+\rs n^2}{2nk}\right).
  \end{equation}
  We also use \emph{e-focus}, \emph{p-focus}, \emph{h-focus} and
  \emph{\(\rs\)-focus}, in line with
  Convention~\ref{co:eph-and-a} to take into account of the type of
  \(\Cliff{\rs}\).

  \emph{Focal length} of a cycle is \(\frac{n}{2k}\).
\end{definition}
\begin{remark}
  Note that focus of \(\cycle{s}{\bs}\) is independent of the sign of
  \(s\). Geometrical meaning of focus is as follows. 
  If a cycle is realised in the parabolic space \(\Space{R}{p}\)
  \emph{h-focus}, \emph{p-focus}, \emph{e-focus} are correspondingly
  geometrical focus of the parabola, its vertex and the point on
  directrix nearest to the vertex, see
  Fig.~\ref{fig:eph-cycle}(b). Thus the traditional focus is h-focus
  in our notations.
\end{remark}

We may describe a finer structure of the cycle space through invariant
subclasses of them. Two such families are \emph{zero-radius} and
\emph{self-adjoint} cycles which are naturally appearing from
expressions~\eqref{eq:trace-det} and~\eqref{eq:inner-product}
correspondingly. 
\begin{definition}
\label{de:zero-radius-cycle}
  \emph{\(\bs\)-Zero-radius cycles} are defined by the condition
  \(\det(\cycle{s}{\bs})=0\), i.e. are explicitly given by matrices
  \begin{equation}
    \label{eq:zero-cycle-matrix}
    \begin{pmatrix}
      y & -y^2\\
      1 & -y
    \end{pmatrix} =\frac{1}{2}
    \begin{pmatrix}
      y & y\\
      1 & 1
    \end{pmatrix}
    \begin{pmatrix}
      1 & -y\\
      1 & -y
    \end{pmatrix} = 
    \begin{pmatrix}
      \se{0} u+\se{1} v & u^2-\bs v^2\\
      1 & -\se{0} u-\se{1} v
    \end{pmatrix}, 
  \end{equation}
  where  \(y=\se{0} u+\se{1} v\).
  We denote such a \(\bs\)-zero-radius cycle by \(\zcycle{s}{\bs}(y)\).
\end{definition}

\begin{figure}[htbp]
  \centering
  \includegraphics[scale=1.1]{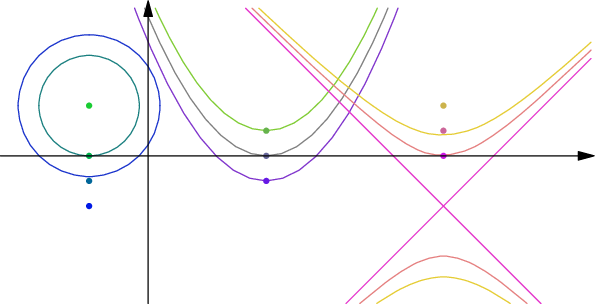}
  \caption[Different implementations of the same
    zero-radius cycles]{Different \(\sigma\)-implementations of the same
    \(\bs\)-zero-radius cycles and corresponding foci.}
  \label{fig:zero-radius}
\end{figure}
Geometrically \(\bs\)-zero-radius cycles are \(\sigma\)-implemented by \(Q\)
from Defn.~\ref{de:cycle-2-matrix} rather differently, see
Fig.~\ref{fig:zero-radius}. Some notable rules are:
\begin{itemize}
\item[(\(\sigma\bs=1\))] Implementations are zero-radius
  cycles in the standard sense: the point \(ue_0-ve_1\) in
  elliptic case and the \emph{light cone} with the centre at
  \(ue_0+ ve_1\) in hyperbolic space~\textup{\cite{Cnops02a}}.
\item[(\(\sigma=0\))] Implementations are parabolas with focal
  length \(v/2\) and the real axis 
  passing through the \(\bs\)-focus.  In other words, for
  \(\bs=-1\) focus at \((u,v)\) (the real axis is
  directrix), for \(\bs=0\) focus at \((u, v/2)\) (the real
  axis passes through the vertex), for \(\bs=1\) focus at
  \((u, 0)\) (the real axis passes through the focus). Such parabolas as
  well have ``zero-radius'' for a suitable parabolic metric, see
  Lemma~\ref{le:distance-first}.
\item[(\(\bs=0\))] \(\sigma\)-Implementations are corresponding conic
  sections which touch the real axis.
\end{itemize}
\begin{remark}
  \label{re:real-line-inv-1}
  The above ``touching'' property of zero-radius cycles for \(\bs=0\)
  is an example of boundary effect inside the domain mentioned in
  Rem.~\ref{re:real-line-inv-0}. It is not surprising after all since
  \(\SL\) action on the upper half-plane may be considered as an
  extension of its action on the real axis.
\end{remark}

\(\bs\)-Zero-radius cycles are significant since they are completely
determined by their centres and thus ``encode'' points into the ``cycle
language''. The following result states that this encoding is M\"obius
invariant as well.
\begin{lemma}
  \label{le:moeb-conj-z-cycle}
  The conjugate \(g^{-1}\zcycle{s}{\bs}(y)g\) of a \(\bs\)-zero-radius cycle
  \(\zcycle{s}{\bs}(y)\) with \(g\in\SL\) is a \(\bs\)-zero-radius cycle
  \(\zcycle{s}{\bs}(g\cdot y)\) with centre at \(g\cdot y\)---the M\"obius
  transform of the centre of \(\zcycle{s}{\bs}(y)\).
\end{lemma}
\begin{proof}
  This may be calculated in \GiNaC~\cite{Kisil05b}*{\S~\ref{G-sec:mobi-invar-cycl}}.
\end{proof}
Another important class of cycles is given by next definition based on
the invariant inner product~\eqref{eq:inner-product} and the
invariance of the real line.
\begin{definition}
   \label{de:self-adj-cycle}
   \emph{Self-adjoint cycle} for \(\bs\neq 0\) are defined by the
   condition \(\Re \scalar{\cycle{s}{\bs}}{\realline{s}{\bs}}=0\),
   where \(\realline{s}{\bs}\) corresponds to the ``real'' axis
   \(v=0\) and \(\Re\) denotes the real part of a Clifford number.
\end{definition}
Explicitly a self-adjoint cycle \(\cycle{s}{\bs}\) is defined
by \(n=0\) in~\eqref{eq:cycle-def-1}. Geometrically they are:
\begin{itemize}
\item[(e, h)] circles or hyperbolas with centres on the real line;
\item[(p)] vertical lines, which are also ``parabolic
  circles''~\textup{\cite{Yaglom79}}, i.e. are given by
  \(\norm{x-y}=r^2\) in the parabolic metric defined below
  in~\eqref{eq:distance-first-par}.
\end{itemize}

\begin{lemma} 
  \label{le:invariant-cycles}
  Self-adjoint cycles form a family, which is invariant under the M\"obius
  transformations.  
\end{lemma}
\begin{proof}
  The proof is either geometrically obvious from the transformations
  described in Section~\ref{sec:mobi-transf}, or follows analytically
  from the action described in Proposition~\ref{pr:transform-cycles}.
\end{proof}
\begin{remark}
  Geometric objects, which are invariant under infinitesimal action of \(\SL\), were
  studied recently in papers~\cites{KonovenkoLychagin08a,Konovenko09a}. 
\end{remark}

\section{Joint invariants: Orthogonality and Inversions}
\label{sec:orth-invers-1}

\subsection{Invariant Orthogonality Type Conditions}
\label{sec:orthogonality}
We already use the matrix invariants of a single cycle in
Definition~\ref{de:focus}, \ref{de:zero-radius-cycle} and
\ref{de:self-adj-cycle}. Now we will consider joint invariants of
several cycles.  Obviously, the relation
\(\tr(\cycle{s}{\bs}\cycle[\tilde]{s}{\bs})=0\) between two
cycles is invariant under M\"obius transforms and characterises the
mutual disposition of two cycles \(\cycle{s}{\bs}\) and
\(\cycle[\tilde]{s}{\bs}\).  More generally the relations
\begin{equation}
  \label{eq:general-relation}
  \tr(p({\cycle{s}{\bs}}^{(1)}, \ldots, {\cycle{s}{\bs}}^{(n)}))=0 
  \quad\text{or}\quad
  \det(p({\cycle{s}{\bs}}^{(1)}, \ldots, {\cycle{s}{\bs}}^{(n)}))=0 
\end{equation}
between \(n\) cycles \({\cycle{s}{\bs}}^{(1)}\), \ldots,
\({\cycle{s}{\bs}}^{(n)}\) based on a polynomial
\(p(x^{(1)},\ldots,x^{(n)})\) of \(n\) non-commuting variables \(x'\),
\ldots, \(x^{(n)}\) is M\"obius invariant if
\(p(x^{(1)},\ldots,x^{(n)})\) is homogeneous in every \(x^{(i)}\). 
Non-homogeneous polynomials will also create M\"obius invariants if we
substitute cycles' \(\det\)-normalised matrices only.  Let us consider some
lower order realisations of~\eqref{eq:general-relation}.

\begin{figure}[htbp]
  \includegraphics[scale=1.2]{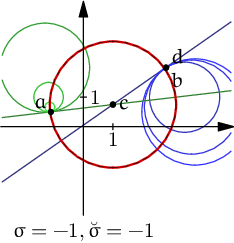}\hfill
  \includegraphics[scale=1.2]{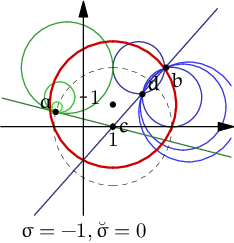}\hfill
  \includegraphics[scale=1.2]{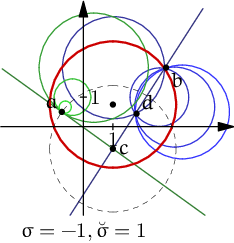}\\[4mm]
  \includegraphics[scale=1.2]{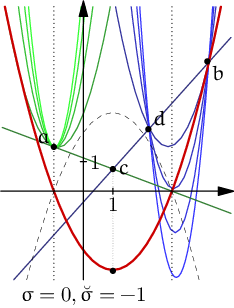}\hfill
  \includegraphics[scale=1.2]{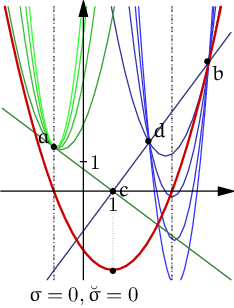}\hfill
  \includegraphics[scale=1.2]{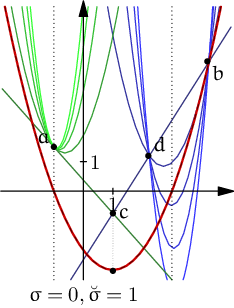}\\[4mm]
  \includegraphics[scale=1.2]{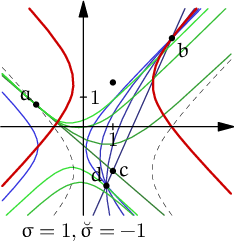}\hfill
  \includegraphics[scale=1.2]{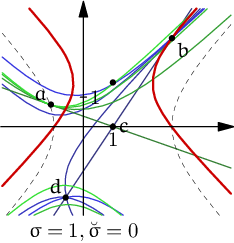}\hfill
  \includegraphics[scale=1.2]{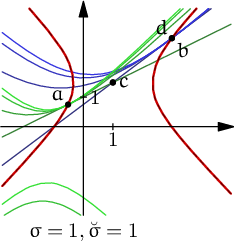}
  \caption[Orthogonality of the first kind]{Orthogonality of the first
    kind in nine combinations.\\
    Each picture presents two groups (green and blue) of cycles which
    are orthogonal to the red cycle \(\cycle{s}{\bs}\).  Point \(b\)
    belongs to \(\cycle{s}{\bs}\) and the family of blue cycles
    passing through \(b\) is orthogonal to \(\cycle{s}{\bs}\). They
    all also intersect in the point \(d\) which is the inverse of
    \(b\) in \(\cycle{s}{\bs}\). Any orthogonality is reduced to the usual
    orthogonality with a new (``ghost'') cycle (shown by the dashed
    line), which may or may not coincide with \(\cycle{s}{\bs}\). For
    any point \(a\) on the ``ghost'' cycle the orthogonality is
    reduced to the local notion in the terms of tangent lines at the
    intersection point. Consequently such a point \(a\) is always the
    inverse of itself.}
  \label{fig:orthogonality1}
\end{figure}

\begin{definition}
  \label{de:orthogonality-first}
  Two cycles \(\cycle{s}{\bs}\) and \(\cycle[\tilde]{s}{\bs}\) are
  \emph{\(\bs\)-orthogonal} if the real part of their inner
  product~\eqref{eq:inner-product} vanishes:
  \begin{equation}
    \label{eq:orthogonality-first}
    \Re\scalar{\cycle{s}{\bs}}{\cycle[\tilde]{s}{\bs}}
    = 0\quad \text{ or, equivalently, }\quad
    \scalar{\cycle{s}{\bs}}{\cycle[\tilde]{s}{\bs}}+\scalar{\cycle[\tilde]{s}{\bs}}{\cycle{s}{\bs}}
    = 0
  \end{equation}
  In light of~\eqref{eq:trace-det} the zero-radius
  cycles (Defn.~\ref{de:zero-radius-cycle}) are also called
  \emph{self-orthogonal} or \emph{isotropic}.
\end{definition}
\begin{lemma}
  The \(\bs\)-orthogonality condition~\eqref{eq:orthogonality-first} is invariant
  under M\"obius transformations.
\end{lemma}
\begin{proof}
  It immediately follows from Definition~\ref{de:orthogonality-first},
  formula~\eqref{eq:cycle-transform-short} and the invariance of trace under
  similarity.
\end{proof}
We also get by the straightforward
calculation~\cite{Kisil05b}*{\S~\ref{G-sec:vari-orth-cond}}: 
\begin{subequations}
  \renewcommand{\theequation}{\theparentequation\ephname{equation}}
  \setcounter{equation}{-1}
  \begin{lemma}
    \label{le:first-orthogonality}
    The
    \(\bs\)-orthogonality~\eqref{eq:orthogonality-first} of
    cycles \(\cycle[\tilde]{s}{\bs}\) and
    \(\cycle{s}{\bs}\) is given through their defining
    equation~\eqref{eq:cycle-def-1} coefficients by
    \begin{eqnarray}
      \label{eq:first-ortho-gen}
      2\bs\tilde{n}n-2\tilde{l}l+\tilde{k}m+\tilde{m}k = 0.
    \end{eqnarray}
    or specifically by
    \begin{eqnarray}
      \label{eq:first-ortho-e}
      -2\tilde{n}n-2\tilde{l}l+\tilde{k}m+\tilde{m}k &=& 0,\\
      \label{eq:first-ortho-p}
      -2\tilde{l}l+\tilde{k}m+\tilde{m}k &=& 0,\\
      \label{eq:first-ortho-h}
      2\tilde{n}n-2\tilde{l}l+\tilde{k}m+\tilde{m}k &=& 0
    \end{eqnarray}
    in the elliptic, parabolic and hyperbolic cases of \(\Cliff{\sigma}\) correspondingly.
  \end{lemma}
\end{subequations}
Note that the orthogonality identity~\eqref{eq:first-ortho-gen} is
linear for coefficients of one cycle if the other cycle is fixed. Thus
we obtain several simple conclusions.
\begin{corollary}
  \begin{enumerate}
  \item A \(\bs\)-self-orthogonal cycle is
    \(\bs\)-zero-radius one~\eqref{eq:zero-cycle-matrix}.
  \item For \(\bs=\pm 1\) there is no non-trivial cycle
    orthogonal to all other non-trivial cycles.\\
    For \(\bs=0\) only the real axis \(v=0\) is orthogonal
    to all other non-trivial cycles.
  \item For \(\bs=\pm 1\) any cycle is uniquely defined
    by the family of cycles orthogonal to it,
    i.e. \((\cycle{s}{\bs}^\perp)^\perp=\{\cycle{s}{\bs}\}\) .\\ 
    For \(\bs= 0\) the set
    \((\cycle{s}{\bs}^\perp)^\perp\) consists of all cycles
    which have the same roots as \(\cycle{s}{\bs}\),
    see middle column of pictures in Fig.~\ref{fig:orthogonality1}.
  \end{enumerate}
\end{corollary}
We can visualise the orthogonality with a zero-radius cycle as follow: 
\begin{lemma}
  \label{le:zero-radius-otho}
  A cycle \(\cycle{s}{\bs}\) is
  \(\bs\)-orthogonal to \(\sigma\)-zero-radius cycle
  \(\zcycle{s}{\sigma }(u,v)\) if 
  \begin{equation}
    \label{eq:ortho-zero-rad}
    k(u^2 k- \sigma v^2) - 2\scalar{(l,n)}{(u, \bs v)}+m =0,
  \end{equation}
  i.e. \(\sigma\)-implementation of \(\cycle{s}{\bs}\) is
  passing through the point \((u, \bs v)\), which \(\bs\)-centre of
  \(\zcycle{s}{\sigma}(u,v)\).
\end{lemma}
The important consequence of the above observations is the possibility
to extrapolate results from zero-radius cycles to the entire space.
\begin{proposition}
  \label{pr:cycle-map-ext}
  Let \(T: \Space{P}{3} \rightarrow \Space{P}{3}\) is
  an orthogonality preserving map of the cycles space,
  i.e. \(\scalar{\cycle{s}{\bs}}{\cycle[\tilde]{s}{\bs}}=0 \
    \Leftrightarrow  \scalar{T\cycle{s}{\bs}}{T\cycle[\tilde]{s}{\bs}}=0\). Then
    for \(\sigma\neq 0\) there is a map \(T_\sigma : \Space{R}{\sigma}
    \rightarrow\Space{R}{\sigma}\), such that \(Q\) intertwines \(T\) and
    \(T_\sigma \):
    \begin{equation}
      \label{eq:Q-intertwine}
       Q T_\sigma  = T Q. 
    \end{equation}
\end{proposition}
\begin{proof}
  If \(T\) preserves the orthogonality (i.e. the inner
  product~\eqref{eq:inner-product} and consequently the determinant
  from~\eqref{eq:trace-det}) then by the image
  \(T\zcycle{s}{\bs}(u,v)\) of a zero-radius cycle
  \(\zcycle{s}{\bs}(u,v)\) is again a zero-radius cycle
  \(\zcycle{s}{\bs}(u_1,v_1)\) and we can define \(T_\sigma \) by the identity
  \(T_\sigma : (u,v)\mapsto (u_1,v_1)\).  

  To prove the intertwining property~\eqref{eq:Q-intertwine} we need
  to show that if a cycle \(\cycle{s}{\bs}\) passes through \((u,v)\)
  then  the image \(T\cycle{s}{\bs}\) passes through \(T_\sigma
  (u,v)\). However for \(\sigma\neq 0\) this is a consequence of the
  \(T\)-invariance of orthogonality and the expression of the
  point-to-cycle incidence through the orthogonality from
  Lemma~\ref{le:zero-radius-otho}.  
\end{proof}
\begin{corollary}
  Let \(T_i: \Space{P}{3} \rightarrow \Space{P}{3}\), \(i=1,2\) are
  two orthogonality preserving maps of the cycles space. If they
  coincide on the subspace of \(\bs\)-zero-radius cycles, \(\bs\neq
  0\), then they are identical in the whole \(\Space{P}{3}\).
\end{corollary}

\begin{remark}
  \label{re:local-ortho}
  Note, that the orthogonality is reduced to local notion in terms of
  tangent lines to cycles in their intersection points \emph{only for
    \(\sigma\bs=1\)}, i.e. this happens only in NW and SE
  corners of Fig.~\ref{fig:orthogonality1}.  In other cases the local
  condition can be formulated in term of ``ghost'' cycle defined below. 
\end{remark}
 We denote by  \(\chi(\sigma)\) the \emph{Heaviside function}:
  \begin{equation}
    \label{eq:heaviside-function}
    \chi(t)=\left\{
      \begin{array}{ll}
        1,& t\geq 0;\\
        -1,& t<0.
      \end{array}\right.
  \end{equation}
\begin{proposition}
  \label{pr:ghost-cycle}
  Let cycles \(\cycle{}{\bs}\) and \(\cycle[\tilde]{}{\bs}\) be
  \(\bs\)-orthogonal. For their \(\sigma\)-implementations we define
  the \emph{ghost cycle} \(\cycle[\hat]{}{\bs}\) by the following two
  conditions:
  \begin{enumerate}
  \item \label{item:centre-centre-rel}
    \(\chi(\sigma)\)-centre of \(\cycle[\hat]{}{\bs}\) coincides
    with  \(\bs\)-centre of \(\cycle{}{\bs}\).
  \item Determinant of  \(\cycle[\hat]{1}{\sigma}\) is equal to
    determinant of \(\cycle{\chi(\bs)}{\sigma}\).
  \end{enumerate}
  Then:
  \begin{enumerate}
  \item   \(\cycle[\hat]{}{\sigma}\) coincides with
    \(\cycle{}{\sigma}\) if \(\sigma \bs=1\);
  \item   \(\cycle[\hat]{}{\sigma}\) has common roots (real or imaginary)
    with \(\cycle{}{\sigma}\);
  \item In the \(\sigma\)-implementation the tangent line to
    \(\cycle[\tilde]{}{\bs}\) at points of its intersections with the
    ghost cycle \(\cycle[\hat]{}{\sigma}\) are passing the
    \(\sigma\)-centre of  \(\cycle[\hat]{}{\sigma}\).
  \end{enumerate}
\end{proposition}
\begin{proof}
  The calculations are done in \GiNaC, see~\cite{Kisil05b}*{\S~\ref{G-sec:ghost-cycle}}.
  For illustration see Fig.~\ref{fig:orthogonality1}, where the ghost
  cycle is shown by the black dashed line.
\end{proof}
Consideration of the ghost cycle does present the orthogonality in the
local terms however it hides the symmetry of this relation.
\begin{remark}
  \label{re:real-line-inv-2}
  Elliptic and hyperbolic ghost cycles are symmetric in the real line,
  the parabolic ghost cycle has its centre on it, see
  Fig.~\ref{fig:orthogonality1}. This is an illustration to the
  boundary effect from Rem.~\ref{re:real-line-inv-0}. 
\end{remark}

\subsection{Inversions in Cycles}
\label{sec:invers-in-cycl}

Definition~\ref{de:cycle-2-matrix} associates a \(2\times 2\)-matrix to any
cycle.  Similarly to \(\SL\) action~\eqref{eq:moebius-def} we can consider
a fraction-linear transformation on \(\Space{R}{\sigma}\) defined by such a
matrix: 
\begin{equation}
  \label{eq:cycle-frac-linear}
  \cycle{s}{\sigma}: u e_0 + v e_1 \mapsto 
  \cycle{s}{\sigma}(u e_0 + v e_1)
  = \frac{(le_0 + ne_1)(u e_0 + v e_1)+m}{k(u e_0 + v e_1)-(le_0 + ne_1)} 
\end{equation}
where \(\cycle{s}{\sigma}\) is as usual~\eqref{eq:matrix-for-cycle}
\begin{displaymath}
  \cycle{s}{\sigma} =
  \begin{pmatrix}
    le_0 + ne_1 & m\\
    k & -le_0 - ne_1
  \end{pmatrix}.
\end{displaymath}
Another natural action of cycles in the matrix form is given by
the conjugation on other cycles:
\begin{equation}
  \label{eq:cycle-conjugation}
  \cycle{s}{\bs}: \cycle[\tilde]{s}{\bs} \mapsto \cycle{s}{\bs}\cycle[\tilde]{s}{\bs}\cycle{s}{\bs}.
\end{equation}
Note that \(\cycle{s}{\bs}\cycle{s}{\bs}= -\det(\cycle{s}{\bs}) I\),
where \(I\) is the identity matrix. Thus the
definition~\eqref{eq:cycle-conjugation} is equivalent to expressions
\(\cycle{s}{\bs}\cycle[\tilde]{s}{\bs}\cycle{s}{\bs}^{-1}\) for
\(\det\cycle{s}{\bs} \neq 0\) since cycles form a projective
space. 
There is a connection between two actions~\eqref{eq:cycle-frac-linear}
and~\eqref{eq:cycle-conjugation} of cycles, which is similar to
\(\SL\) action in  Lemma~\ref{le:moeb-conj-z-cycle}.
\begin{lemma} Let \(\det \cycle{s}{\bs} \neq 0\), then: 
  \label{le:cycle-conj-is-moeb}
  \begin{enumerate}
  \item The conjugation~\eqref{eq:cycle-conjugation} preserves the
    orthogonality relation~\eqref{eq:orthogonality-first}. 
  \item The image \(\cycle{s_2}{\sigma}\zcycle[\tilde]{s_1}{\sigma}(u,
    v)\cycle{s_2}{\sigma}\) of a \(\sigma\)-zero-radius cycle
    \(\zcycle[\tilde]{s_1}{\sigma}\) under the
    conjugation~\eqref{eq:cycle-conjugation} is a
    \(\sigma\)-zero-radius cycle \(\zcycle[\tilde]{s_1}{\sigma}(u',
    v')\), where \((u',v')\) is calculated by the linear-fractional
    transformation~\eqref{eq:cycle-frac-linear}
    \((u',v')=\cycle{s_1s_2}{\sigma}(u,v)\) associated to the cycle
    \(\cycle{s_1s_2}{\sigma}\).
  \item Both formulae~\eqref{eq:cycle-frac-linear}
    and~\eqref{eq:cycle-conjugation} define the same transformation of
    the point space. 
  \end{enumerate}
\end{lemma}
\begin{proof}
  The first part is obvious, the second is calculated in
  \GiNaC~\cite{Kisil05b}*{\S~\ref{G-sec:cycles-conjugation}}.  The last
  part follows from the first two and Prop.~\ref{pr:cycle-map-ext}.
\end{proof}
\begin{figure}[htbp]
  \centering
  (a)\includegraphics[scale=1.5]{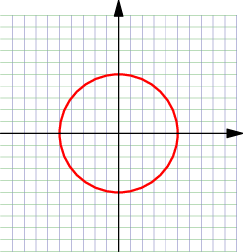}\qquad
  (b)\includegraphics[scale=1.5]{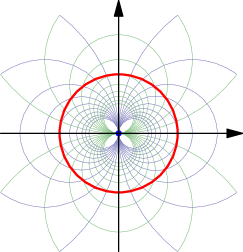}\\[5pt]
  (c)\includegraphics[scale=1.5]{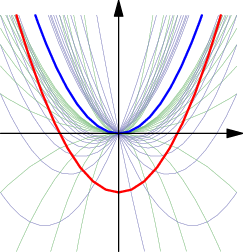}\qquad
  (d)\includegraphics[scale=1.5]{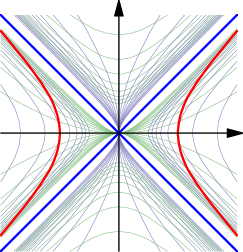}
  \caption[Three types of inversions of the rectangular grid]{Three
    types of inversions of the rectangular grid. The initial
    rectangular grid (a) is inverted elliptically in the unit circle
    (shown in red) on (b), parabolically on (c) and hyperbolically on
    (d). The blue cycle (collapsed to a point at the origin on (b))
    represent the image of the cycle at infinity under inversion.}  
  \label{fig:inversions}
\end{figure}

There are at least two natural ways to define inversions in cycles.
One of them use the orthogonality condition, another define them as
``reflections in cycles''.
\begin{definition}
  \label{def:inversion}
  \begin{enumerate}
  \item\label{item:inversion-orthogonality} \emph{Inversion} in a
    cycle \(\cycle{s}{\sigma}\) sends a point \(p\) to the second point
    \(p'\) of intersection of all cycles orthogonal to
    \(\cycle{s}{\sigma}\) and passing through \(p\).
  \item\label{item:inversion-reflection} \emph{Reflection} in a cycle
     \(\cycle{s}{\sigma}\) is given by \(M^{-1}RM\) where \(M\) sends the
    cycle  \(\cycle{s}{\sigma}\) into the horizontal axis and \(R\) is the mirror
    reflection in that axis.
  \end{enumerate}
\end{definition}
We are going to see that inversions are given
by~\eqref{eq:cycle-frac-linear} and reflections are expressed
through~\eqref{eq:cycle-conjugation}, thus they are essentially the
same in light of Lemma~\ref{le:cycle-conj-is-moeb}.
\begin{remark}
  \label{re:dual-limit} Here is a simple example where usage of
  complex (dual or double) numbers is weaker then Clifford algebras,
  see Rem.~\ref{re:dual-double}. A reflection of a cycle in the axis \(v=0\) is
  represented by the conjugation~\eqref{eq:cycle-conjugation} with
  the corresponding matrix \(
  \begin{pmatrix}
    \se{1}&0\\0&-\se{1}
  \end{pmatrix}\). The same transformation in term of complex numbers
  should involve a complex conjugation and thus cannot be expressed by
  multiplication.
\end{remark}
Since we have three different EPH orthogonality between cycles
there are also three different inversions:
\begin{proposition}
  \label{pr:orthogonality-to-cycle}
  A cycle \(\cycle[\tilde]{s}{\bs}
  \) is
  orthogonal to a cycle \(\cycle{s}{\bs}
  \) if for any point \(u_1e_0+v_1e_1\in \cycle[\tilde]{s}{\bs}\)
  the cycle \(\cycle[\tilde]{s}{\bs}\) is also passing through its
  image
  \begin{equation}
    \label{eq:inversion-ortho}
    u_2e_0+v_2e_1 
    =\begin{pmatrix}
      le_0+sre_1 & m\\
      k & le_0+sre_1
    \end{pmatrix}
    (u_1e_0+v_1e_1)
  \end{equation}
  under the M\"obius transform defined by the matrix
  \(\cycle{s}{\bs}\). Thus the point
  \(u_2e_0+v_2e_1=\cycle{s}{\bs}(u_1e_0+v_1e_1)\) is the inversion of
  \(u_1e_0+v_1e_1\) in \(\cycle{s}{\bs}\).
\end{proposition}
\begin{proof}
  The symbolic calculations done by \GiNaC\cite{Kisil05b}*{\S~\ref{G-sec:orth-invers}}.
\end{proof}
\begin{proposition}
  The reflection~\ref{def:inversion}.\ref{item:inversion-reflection} of
  a zero-radius cycle \(\zcycle{s}{\bs}\) in a cycle \(\cycle{s}{\bs}\) is
  given by the conjugation: \(\cycle{s}{\bs}\zcycle{s}{\bs}\cycle{s}{\bs}\).
\end{proposition}
\begin{proof}
  Let \(\cycle[\tilde]{s}{\bs}\) has the property \(\cycle[\tilde]{s}{\bs}
  \cycle{s}{\bs} \cycle[\tilde]{s}{\bs} = \Space{R}{}\). Then
  \(\cycle[\tilde]{s}{\bs} \Space{R}{} \cycle[\tilde]{s}{\bs} =
  \cycle{s}{\bs}\). Mirror reflection in the real line is given by the
  conjugation with \(\Space{R}{}\), thus the transformation described
  in~\ref{def:inversion}.\ref{item:inversion-reflection} is a conjugation with the cycle
  \(\cycle[\tilde]{s}{\bs} \Space{R}{} \cycle[\tilde]{s}{\bs} =
  \cycle{s}{\bs}\) and thus coincide with~\eqref{eq:inversion-ortho}.
\end{proof}
The cycle \(\cycle[\tilde]{s}{\bs}\) from the above proof can be
characterised as follows. 
\begin{lemma}
  \label{le:cycle-refl-real-line}
  Let \(\cycle{s}{\bs}=(k, l, n,m)\) be a cycle and for \(\bs\neq 0\)
  the \(\cycle[\tilde]{s}{\bs}\) be given by \((k, l,
    n\pm\sqrt{\det\cycle{\sigma}{\bs}},m)\). Then 
  \begin{enumerate}
  \item \(\cycle[\tilde]{s}{\bs} \cycle{s}{\bs} \cycle[\tilde]{s}{\bs}
    = \Space{R}{}\) and \(\cycle[\tilde]{s}{\bs} \Space{R}{}
    \cycle[\tilde]{s}{\bs} =  \cycle{s}{\bs}\)
  \item\label{it:inv-cycle-root}
    \(\cycle[\tilde]{s}{\bs}\) and \(\cycle{s}{\bs}\) have common
    roots.
  \item \label{it:inv-cycle-center}
    In the \(\bs\)-implementation the cycle \(\cycle{s}{\bs}\) passes
    the centre of \(\cycle[\tilde]{s}{\bs}\).
  \end{enumerate}
\end{lemma}
\begin{proof}
  This is calculated by
  \GiNaC~\cite{Kisil05b}*{\S~\ref{G-sec:reflection-cycles}}. Also one
  can direct observe~\ref{le:cycle-refl-real-line}.\ref{it:inv-cycle-root} for real roots, since
  they are fixed points of the inversion. Also the transformation of
  \(\cycle{s}{\bs}\) to a flat cycle implies that \(\cycle{s}{\bs}\)
  is passing the centre of inversion, hence~\ref{le:cycle-refl-real-line}.\ref{it:inv-cycle-center}. 
\end{proof}

In~\cite{Yaglom79}*{\S~10} the \emph{inversion of second kind}
related to a parabola \(v=k(u-l)^2+m\) was defined by the map:
\begin{equation}
  \label{eq:invers-sec-kind}
  (u,v) \mapsto (u, 2(k(u-l)^2+m)-v),
\end{equation}
i.e. the parabola bisects the vertical line joining a point and its
image. Here is the result expression this transformation through the
usual inversion in parabolas:
\begin{proposition}
  \label{pr:inv-2nd-kind-decomp}
  The inversion of second kind~\eqref{eq:invers-sec-kind} is a
  composition of three inversions:  in parabolas \(u^2- 2lu
  -4mv-m/k=0\), \(u^2- 2lu-m/k=0\), and the real line.
\end{proposition}
\begin{proof}
   See symbolic calculation in~\cite{Kisil05b}*{\S~\ref{G-sec:invers-second-kind}}.
\end{proof}
\begin{remark}
  Yaglom in~\cite{Yaglom79}*{\S~10} considers the usual inversion (``of the
  first kind'') only in degenerated parabolas (``parabolic circles'')
  of the form \(u^2-2lu+m=0\). However the inversion of the second
  kind requires for its decomposition like in
  Prop.~\ref{pr:inv-2nd-kind-decomp} at least one inversion in a proper
  parabolic cycle \(u^2-2lu-2nv+m=0\). Thus such inversions are indeed
  of \emph{another kind} within Yaglom's framework~\cite{Yaglom79},
  but are not in our.

  Another important difference between inversions
  from~\cite{Yaglom79} and our wider set of 
  transformations~\eqref{eq:cycle-frac-linear} is what ``special''
  (vertical) lines does not form an invariant set, as can be
  seen from Fig.~\ref{fig:inversions}(c), and thus they are not
  ``special'' lines anymore.
\end{remark}

\subsection{Focal Orthogonality}
\label{sec:focal-orthogonality}
It is natural to consider invariants of higher orders which are
generated by~\eqref{eq:general-relation}. Such invariants shall have
at least one of the following properties
\begin{itemize}
\item contains a non-linear power of the same cycle;
\item accommodate more than two cycles.
\end{itemize}
The consideration of higher order invariants is similar to a
transition from Riemannian geometry to Finsler
one~\cites{Chern96a,Garasko09a,Pavlov06a}. 

It is interesting that higher order invariants 
\begin{enumerate}
\item can be built on top of the already defined ones;
\item can produce lower order invariants.
\end{enumerate}
For each of the two above transitions we consider an example. We already
know that a similarity of a cycle with another cycle is a new
cycle~\eqref{eq:cycle-conjugation}. The inner product of later with a
third given cycle form a joint invariant of those three cycles:
\begin{equation}
  \label{eq:trio}
  \scalar{\cycle{s}{\bs}_1 \cycle{s}{\bs}_2\cycle{s}{\bs}_1}{\cycle{s}{\bs}_3},
\end{equation}
which is build from the second-order invariant
\(\scalar{\cdot}{\cdot}\). Now we can reduce the order of this
invariant by fixing \(\cycle{s}{\bs}_3\) be the real line (which is
itself invariant).  The obtained invariant of two cycles deserves a
special consideration.  Alternatively it emerges from
Definitions~\ref{de:orthogonality-first} and~\ref{de:self-adj-cycle}.

\begin{definition}
  \label{de:f-ortho}
  The \emph{focal orthogonality} (f-orthogonality) of
  a cycle \(\cycle{s}{\bs}\) \emph{to} a cycle  \(\cycle[\tilde]{s}{\bs}\) is
  defined by the condition that the cycle \(\cycle{s}{\bs}
  \cycle[\tilde]{s}{\bs} \cycle{s}{\bs}\) is orthogonal (in the
  sense of Definition~\ref{de:orthogonality-first}) to the real line,
  i.e is a self-adjoint cycle in the sense of
  Definition~\ref{de:self-adj-cycle}. Analytically this is defined by
  \begin{equation}
    \label{eq:f-orthog-def}
    \Re\tr(\cycle{s}{\bs} \cycle[\tilde]{s}{\bs}\cycle{s}{\bs}\realline{s}{\bs})=0
  \end{equation}
  and  we denote it by \(\cycle{s}{\bs} \sperp   \cycle[\tilde]{s}{\bs}\).
\end{definition}
\begin{remark}
  \label{re:real-line-inv-3}
  This definition is explicitly based on the invariance of the real
  line and is an illustration to the boundary value effect from
  Rem.~\ref{re:real-line-inv-0}.
\end{remark}
\begin{figure}[htbp]
  \includegraphics[scale=1.2]{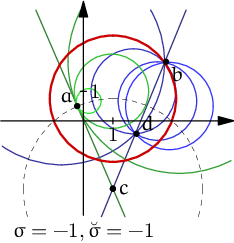}\hfill
  \includegraphics[scale=1.2]{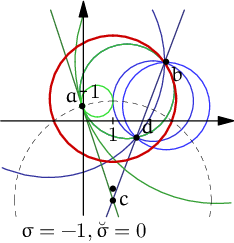}\hfill
  \includegraphics[scale=1.2]{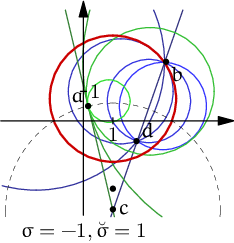}\\[4mm]
  \includegraphics[scale=1.2]{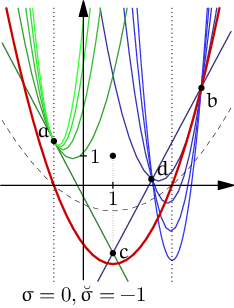}\hfill
  \includegraphics[scale=1.2]{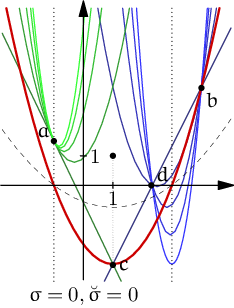}\hfill
  \includegraphics[scale=1.2]{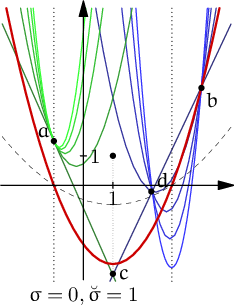}\\[4mm]
  \includegraphics[scale=1.2]{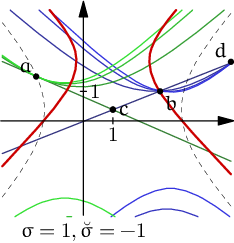}\hfill
  \includegraphics[scale=1.2]{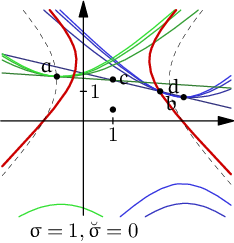}\hfill
  \includegraphics[scale=1.2]{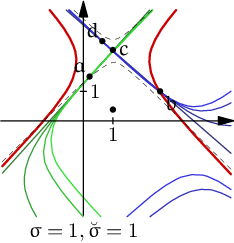}
  \caption[Focal Orthogonality]{Focal orthogonality in all nine
    combinations. To highlight both similarities and distinctions with
    the ordinary orthogonality we use the same notations as that in
    Fig.~\ref{fig:orthogonality1}. The cycles
    \(\cycle[\tilde]{\bs}{\bs}\) from
    Proposition~\ref{pr:focal-ortho} are drawn by dashed lines.}
  \label{fig:orthogonality2}
\end{figure}
\begin{remark}
  It is easy to observe the following
  \begin{enumerate}
  \item f-orthogonality is not a symmetric: \(\cycle{s}{\bs} \sperp
  \cycle[\tilde]{s}{\bs}\) does not implies \(\cycle[\tilde]{s}{\bs} \sperp
  \cycle{s}{\bs}\);
  \item 
    \label{it:s-orth-sl2-inv}
    Since the real axis \(\Space{R}{}\) and
    orthogonality~\eqref{eq:orthogonality-first} are \(\SL\)-invariant
    objects f-orthogonality is also \(\SL\)-invariant.
  \end{enumerate}
\end{remark}
However an invariance of
f-orthogonality under inversion of cycles required some study since
the real line is not an invariant of such transformations in general.
\begin{lemma}
  The image \(\cycle{s_1}{\bs} \realline{s}{\bs} \cycle{s_1}{\bs}\) of the
  real line under inversion in \(\cycle{s_1}{\bs}=(k,l,n,m)\) is the cycle:
  \begin{displaymath}
    (2 s s_1 \bs k n,\ 2 s s_1 \bs l n,\ s^2 (l^2+\bs n^2-m k),\ 2 s s_1 \bs m n).
  \end{displaymath}
  It is the real line if \(s\cdot\det(\cycle{s_1}{\bs})\neq0\) and   either
  \begin{enumerate}
  \item  \(s_1 n=0\), in this case it is a composition of \(\SL\)-action by \(
    \begin{pmatrix}
      l&-me_0\\ke_0 &-l
    \end{pmatrix}\) and the reflection in the real line; or
  \item  \(\bs=0\), i.e. the parabolic case of the cycle space.
  \end{enumerate}
  If this condition is satisfied than f-orthogonality
  preserved by the inversion \(\cycle[\tilde]{s}{\bs}\rightarrow \cycle{s_1}{\bs} \cycle[\tilde]{s}{\bs}
  \cycle{s_1}{\bs}\) in \(\cycle[\tilde]{s}{\bs}\).
\end{lemma}
The following explicit expressions of f-orthogonality reveal further
connections with cycles' invariants.
\begin{proposition} f-orthogonality of \({\cycle{s}{p}}\) to
  \({\cycle[\tilde]{s}{p}}\) is given by either of the following
  equivalent identities
  \begin{eqnarray*}
     \tilde{n} (l^{2}-\se{1}^2 n^{2}- m k) + \tilde{m} n k-2\tilde{l}
     n l+ \tilde{k} m n & =& 0, \qquad \textrm{ or} \\
     \tilde{n} \det(\cycle{s}{\bs}) + n\scalar{\cycle{s}{p}}{\cycle[\tilde]{s}{p}} & =& 0. 
  \end{eqnarray*}
\end{proposition}
\begin{proof}
  This is another \GiNaC\ calculation~\cite{Kisil05b}*{\S~\ref{G-sec:expr-orth-s.k}}
\end{proof}
The f-orthogonality may be again related to the usual orthogonality
through an appropriately chosen \emph{f-ghost cycle}, compare the
next Proposition with Prop.~\ref{pr:ghost-cycle}:
\begin{proposition}
  \label{pr:focal-ortho}
  Let \(\cycle{s}{\bs}\) be a cycle, then its \emph{f-ghost cycle}
  \(\cycle[\tilde]{\bs}{\bs} = \cycle{\chi(\sigma)}{\bs}
  \Space[\bs]{R}{\bs} \cycle{\chi(\sigma)}{\bs}\) is the reflection of
  the real line in \(\cycle{\chi(\sigma)}{\bs}\), where
  \(\chi(\sigma)\) is the \emph{Heaviside function}~\ref{eq:heaviside-function}.
  Then
  \begin{enumerate}
  \item Cycles \(\cycle{s}{\bs}\) and \(\cycle[\tilde]{\bs}{\bs}\)
    have the same roots.
  \item \label{item:focal-centre-rel}
    \(\chi(\sigma)\)-Centre of \(\cycle[\tilde]{\bs}{\bs}\) coincides with the \(\bs\)-focus of
    \(\cycle{s}{\bs}\), consequently all lines f-orthogonal to
    \(\cycle{s}{\bs}\) are passing one of its foci.
  \item s-Reflection in \(\cycle{s}{\bs}\) defined from
    f-orthogonality (see
    Definition~\ref{def:inversion}.\ref{item:inversion-orthogonality}) coincides with
    usual inversion in \(\cycle[\tilde]{\bs}{\bs}\).
  \end{enumerate}
\end{proposition}
\begin{proof}
  This again is calculated in \GiNaC, see~\cite{Kisil05b}*{\S~\ref{G-sec:invers-from-orth}}.
\end{proof}
For the reason~\ref{pr:focal-ortho}.\ref{item:focal-centre-rel} this relation between
cycles may be labelled as \emph{focal orthogonality},
cf. with~\ref{pr:ghost-cycle}.\ref{item:centre-centre-rel}. It can generates the
corresponding inversion similar to
Defn.~\ref{def:inversion}.\ref{item:inversion-orthogonality} which obviously reduces to
the usual inversion in the f-ghost cycle. 
The extravagant f-orthogonality will unexpectedly appear again from
consideration of length and distances in the next section and is
useful for infinitesimal cycles \S~\ref{sec:zero-length-cycles}.

\section{Metric Properties from Cycle Invariants}
\label{sec:metric-properties}

So far we discussed only invariants like orthogonality, which are
related to angles. Now we turn to metric properties similar to
distance. 

\subsection{Distances and Lengths}
\label{sec:lengths-orth}

The covariance of cycles (see Lemma~\ref{le:invariance-of-cycles})
suggests them as ``circles'' in each of the EPH cases. 
Thus we play \emph{the standard mathematical game}: turn some
properties of classical objects into definitions of new ones.
\begin{definition}
  \label{de:radius}
  The \(\bs\)-\emph{radius} of a cycle \(\cycle{s}{\bs}\) if squared is equal
  to the \(\bs\)-determinant of cycle's \(k\)-normalised (see
  Defn.~\ref{de:normalisation}) matrix, i.e.
  \begin{equation}
    \label{eq:radius}
    r^2= \frac{\det\cycle{s}{\bs}}{k^2} =\frac{l^2- \bs n^{2}-km}{k^2}.
  \end{equation}
  As usual, the \(\bs\)-\emph{diameter} of a cycles is two times its radius.
\end{definition}
\begin{lemma}
  The \(\bs\)-\emph{radius} of a cycle \(\cycle{s}{\bs}\) is equal to
  \(1/k\), where \(k\) is \((2,1)\)-entry of \(\det\)-normalised
  matrix (see Defn.~\ref{de:normalisation}) of the cycle.
\end{lemma}

\begin{figure}[htbp]
  \centering
  (a) \includegraphics[scale=1.2]{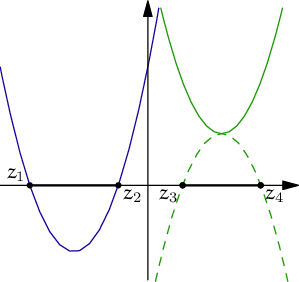}\hfill
  (b) \includegraphics[scale=1.2]{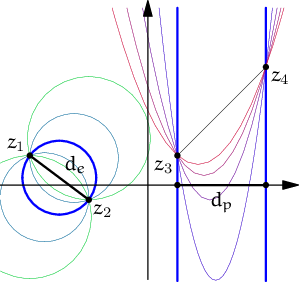}
  \caption[Radius and distance]{(a) The square of the parabolic diameter
    is the square of the distance between roots if they are real (\(z_1\)
    and \(z_2\)), otherwise the negative square of the distance between the
    adjoint roots (\(z_3\) and \(z_4\)).\\
    (b) Distance as extremum of diameters in elliptic (\(z_1\) and
    \(z_2\)) and parabolic (\(z_3\) and \(z_4\)) cases.}
  \label{fig:distances}
\end{figure}

Geometrically in various EPH cases this corresponds to the following
\begin{enumerate}
\item[(e, h)] The value of~\eqref{eq:radius} is the usual radius of
  a circle or hyperbola; 
\item[(p)] The diameter of a parabola is the
  (Euclidean) distance between its (real) roots, i.e.  solutions of
  \(ku^2-2lu+m=0\), or roots of its ``adjoint'' parabola
  \(-ku^2+2lu+m-\frac{2l^2}{k}=0\) (see
  Fig.~\ref{fig:distances}(a)).
\end{enumerate}
\begin{remark}
  Note that
  \begin{displaymath}
    r^2_{\bs}=-4*f*u_{\bs},    
  \end{displaymath}
  where \(r^2_{\bs}\) is the square of cycle's \(\bs\)-radius,
  \(u_{\bs}\) is the second coordinate of its \(\bs\)-focus and \(f\)
  its focal length.
\end{remark}
An intuitive notion of a distance in both mathematics and the everyday
life is usually of a variational nature. We natural perceive the
shortest distance between two points delivered by the straight lines
and only then can define it for curves through an approximation.  This
variational nature echoes also in the following definition.
\begin{definition}
  \label{de:distance}
  The \((\sigma,\bs)\)-\emph{distance} between two points is the extremum of
  \(\bs\)-dia\-me\-ters for all \(\sigma\)-cycles passing through both points.
\end{definition}

During geometry classes we oftenly make measurements with a compass,
which is based on the idea that \emph{a cycle is locus of points equidistant
from its centre}. We can expand it for all cycles in the following
definition:
\begin{definition}
  \label{de:length}
  The \(\bs\)-\emph{length} from a \(\rs\)-centre or from a \(\rs\)-focus
  of a directed interval \(\lvec{AB}\) is the \(\bs\)-radius of the
  \(\sigma\)-cycle with its \(\rs\)-\emph{centre} or
  \(\rs\)-\emph{focus} correspondingly at the point \(A\) which passes
  through \(B\). These lengths are denoted by \(l_c(\lvec{AB})\) and
  \(l_f(\lvec{AB})\) correspondingly.
\end{definition}
\begin{remark}
  \begin{enumerate}
  \item Note that the distance is a symmetric functions of two points
    by its definition and this is not necessarily true for
    lengths. For modal logic of non-symmetric distances
    see, for example, \cite{KuruczWolterZakharyaschev05}. However the
    first axiom (\(l(x,y)=0\) iff  \(x=y\)) should be modified as follows:
    \begin{displaymath}
     ( l(x,y)=0 \text{ and } l(x,y)=0 ) \text{ iff } x=y.
    \end{displaymath}
  \item A cycle is uniquely defined by elliptic or hyperbolic centre
    and a point which it passes. However the parabolic centre is not
    so useful. Correspondingly \((\sigma,0)\)-length from parabolic
    centre is not properly defined. 
  \end{enumerate}
\end{remark}
\begin{lemma}
  \label{le:distance-first} 
  \begin{enumerate}
  \item The cycle of the form~\eqref{eq:zero-cycle-matrix} has zero
  radius. 
  \item The distance between two points \(y=e_0u+e_1v\)
       and  \(y'=e_0u'+e_1v'\) in the elliptic or hyperbolic spaces is
    \begin{equation}
      \label{eq:distance-first-ell-hyp}
      d^2(y, y') = \frac{ \bs ((u-u')^2-\sigma(v- v')^2) +4(1-\sigma\bs) v v'}
      {(u- u')^2 \bs-(v-v')^2} ((u-u')^2 -\sigma(v- v')^2),
    \end{equation}
    and in parabolic case it is (see Fig.~\ref{fig:distances}(b)
    and~\cite{Yaglom79}*{p.~38, (5)}) 
    \begin{equation}
      \label{eq:distance-first-par}
      d^2(y, y') = (u-u')^2.
    \end{equation}
  \end{enumerate}
\end{lemma}
\begin{proof}
  Let \(\cycle{\sigma}{s}(l)\) be the family of cycles passing through
  both points \((u, v)\) and \((u', v')\) (under the assumption
  \(v\neq v'\)) and parametrised by its coefficient \(l\) in the
  defining equation~\eqref{eq:cycle-def}.  By a calculation done in
  \GiNaC~\cite{Kisil05b}*{\S~\ref{G-sec:dist-betw-points}} we found
  that the only critical point of \(\det(\cycle{\sigma}{s}(l))\) is:
  \begin{equation}
    \label{eq:critical-point}
    l_0 = 
    \frac{1}{2}\left((u'+u) +
    (\bs\sigma-1)\frac{(u'-u)(v^2-v'^2)}
    {(u'- u)^2 \bs-(v-v')^2} \right),
  \end{equation}
  [Note that in the case \(\sigma\bs=1\), i.e. both points and cycles
  spaces are simultaneously either elliptic or hyperbolic, this
  expression reduces to the expected midpoint
  \(l_0=\frac{1}{2}(u+u')\).] Since in the elliptic or hyperbolic case
  the parameter \(l\) can take any real value, the extremum of
  \(\det(\cycle{\sigma}{s}(l))\) is reached in \(l_0\) and is equal
  to~\eqref{eq:distance-first-ell-hyp} (calculated by
  \GiNaC~\cite{Kisil05b}*{\S~\ref{G-sec:dist-betw-points}}). A separate
  calculation for the case \(v= v'\) gives the same answer.

  In the parabolic case the possible values of \(l\) are either in
  \((-\infty, \frac{1}{2}(u+u'))\), or \((\frac{1}{2}(u+u'),\infty)\),
  or the only value is \(l=\frac{1}{2}(u+u')\) since for that value a
  parabola should flip between upward and downward directions of its
  branches. In any of those cases the extremum value corresponds to the
  boundary point \(l=\frac{1}{2}(u+u')\) and is equal
  to~\eqref{eq:distance-first-par}.
\end{proof}
\begin{corollary}
  \label{co:inner-product-dist}
  If cycles \(\cycle{s}{\bs}\) and \(\cycle[\tilde]{s}{\bs}\) are
  normalised by conditions \(k=1\) and \(\tilde{k}=1\) then
  \begin{displaymath}
    \scalar{\cycle{s}{\bs}}{\cycle[\tilde]{s}{\bs}} =
    \modulus{c-\tilde{c}}_{\bs}^2-r_{\bs}^2-\tilde{r}_{\bs}^2, 
  \end{displaymath}
  where  \(\modulus{c-c}_{\bs}^2=(l-\tilde{l})^2-\bs(n-\tilde{n})^2\) is
  the square of \(\bs\)-distance between cycles' centres,
  \(r_{\bs}\) and \(\tilde{r}_{\bs}\) are \(bs\)-radii of the
  respective cycles.
\end{corollary}
To get feeling of the identity~\eqref{eq:distance-first-ell-hyp} we may observe, that:
\begin{eqnarray*}
  d^2(y, y') &=& (u-u')^2 + (v- v')^2,\qquad \textrm{ for elliptic
    values } \sigma=\bs=-1;\\
  d^2(y, y') &=& (u-u')^2 - (v- v')^2,\qquad \textrm{ for
    hyperbolic values } \sigma=\bs=1;
\end{eqnarray*}
i.e. these are familiar expressions for the elliptic and hyperbolic
spaces. However four other cases (\(\sigma\bs=-1\) or \(0\))
gives quite different results. For example, \(d^2(y, y')\not\rightarrow
0\) if \(y\) tense to \(y'\) in the usual sense.
\begin{remark}
  \begin{enumerate}
  \item In the three cases \(\sigma=\bs=-1\), \(0\) or \(1\), which
    were typically studied before, the above distances are conveniently defined
    through the Clifford algebra multiplications
    \begin{displaymath}
      d_{e,p,h}^2(ue_0+ve_1)=-(ue_0+ve_1)^2.    
    \end{displaymath}
  \item Unless \(\sigma=\bs\) the parabolic
    distance~\eqref{eq:distance-first-par} is not received
    from~\eqref{eq:distance-first-ell-hyp} by the substitution
    \(\sigma =0\).
  \end{enumerate}
\end{remark}
Now we turn to calculations of the lengths.
\begin{lemma}
  \label{le:distance-second}
  \begin{enumerate}
  \item \label{it:length-centre}
    The \(\bs\)-length from the \(\rs\)-centre between two
    points \(y=e_0u+e_1v\) and \(y'=e_0u'+e_1v'\) is
    \begin{equation}
      \label{eq:k-center-point}
      l_{c_{\bs}}^2(y, y')   =  
 (u-u')^2-\sigma v'^2+2\rs v v' -\bs v^2.
    \end{equation}
  \item \label{it:length-focus}
    The \(\bs\)-length from the \(\rs\)-focus between two
    points \(y=e_0u+e_1v\) and \(y'=e_0u'+e_1v'\) is
    \begin{equation}
      \label{eq:k-focus-point}
      l_{f_{\bs}}^2(y, y')  = (\rs-\bs) p^2-2vp,
    \end{equation}
    where:
    \begin{eqnarray}
      \label{eq:focal-length}
      p & = & \rs\left(-(v'-v)\pm\sqrt{\rs(u'-u)^2+(v'-v)^2-\sigma\rs
          v'^2}\right), \quad\text{if } \rs\neq0;\\
      \label{eq:parab-focal-length}
      p & = & \frac{(u'-u)^2-\sigma v'^2}{2(v'-v)}, \quad\text{if } \rs=0.
    \end{eqnarray}
  \end{enumerate}
\end{lemma}
\begin{proof}
  Identity~\eqref{eq:k-center-point} is verified in
  \GiNaC~\cite{Kisil05b}*{\S~\ref{G-sec:length-direct-interv}}. For the
  second part we observe that the parabola with the focus \((u,v)\)
  passing through \((u',v')\) has the following parameters:
  \begin{displaymath}
   k=1,\quad l= u,\quad n=p,\quad  m = 2\rs pv'-u'^2+2uu'+\sigma v'^2.
  \end{displaymath}
  Then the formula~\eqref{eq:k-focus-point} is verified by the \GiNaC\
  calculation~\cite{Kisil05b}*{\S~\ref{G-sec:length-interv-from}}.
\end{proof}
\begin{remark}
  \label{re:focal-lenght}
  \begin{enumerate}
  \item     The value of  \(p\) in~\eqref{eq:focal-length} is the
    focal length of either of the two cycles, which are in the 
    parabolic case upward or downward parabolas (corresponding to the
    plus or minus signs) with focus at \((u, v)\) and passing through
    \((u', v')\).
  \item \label{it:more-lengths}
    In the case \(\sigma\bs=1\) the
    length~\eqref{eq:k-center-point} became the standard elliptic
    or hyperbolic distance \((u-u')^2-\sigma (v-v')^2\) obtained
    in~\eqref{eq:distance-first-ell-hyp}. Since these expressions
    appeared both as distances and lengths they are widely used. 

    On the other hand in the parabolic
    space we get three additional lengths besides of
    distance~\eqref{eq:distance-first-par}:   
    \begin{displaymath}
      l_{c_{\bs}}^2(y, y')    =  (u-u')^2+2 v v'-\bs v^2
    \end{displaymath}
    parametrised by \(\bs\) (cf.~Remark~\ref{rem:prejudice}.\ref{it:par-is-zero}). 
  \item The parabolic distance~\eqref{eq:distance-first-par} can
    be expressed as 
    \begin{displaymath}
      d^2(y, y') = p^2+2(v-v')p
    \end{displaymath}
    in terms of the focal length~\eqref{eq:focal-length}, which is an
    expression similar to~\eqref{eq:k-focus-point}.
  \end{enumerate}
\end{remark}

\subsection[Conformal Properties of Moebius Maps]{Conformal Properties of M\"obius Maps}
\label{sec:conf-prop-moeb}

All \emph{lengths} \(l(\lvec{AB})\) in \(\Space{R}{\sigma}\) from
Definition~\ref{de:length} are such that for a fixed point \(A\) all
level curves of \(l(\lvec{AB})=c\) are corresponding cycles: circles,
parabolas or hyperbolas, which are covariant objects in the
appropriate geometries. Thus we can expect some covariant properties
of distances and lengths.
\begin{definition}
  \label{de:conformal}
  We say that a distance or a length \(d\) is \(\SL\)-\emph{conformal} if for
  fixed \(y\), \(y'\in\Space{R}{\sigma}\) the limit:
  \begin{equation}
    \label{eq:conformal-lim}
    \lim_{t\rightarrow 0} \frac{d(g\cdot y, g\cdot(y+ty'))}{d(y,
      y+ty')}, \qquad
    \text{ where } g\in\SL, 
  \end{equation}
  exists and its value depends only from \(y\) and \(g\) and is
  independent from \(y'\). 
\end{definition}
The following proposition shows that \(\SL\)-conformality is not rare. 
\begin{proposition}
  \label{pr:conformity}
  \begin{enumerate}
  \item \label{it:conformity-dist}
    The distance~\eqref{eq:distance-first-ell-hyp} is conformal if
    and only if the type of point and cycle spaces are the same, i.e.
    \(\sigma\bs=1\). The parabolic
    distance~\eqref{eq:distance-first-par} is conformal only in the
    parabolic point space.
  \item \label{it:conformity-length-centre}
    The lengths from centres~\eqref{eq:k-center-point} are
    conformal for any combination of values of \(\sigma\),
    \(\bs\)  and \(\rs\). 
  \item \label{it:conformity-length-foci} The lengths from
    foci~\eqref{eq:k-focus-point} are conformal for \(\rs\neq 0\) and
    any combination of values of \(\sigma\) and \(\bs\).
  \end{enumerate}\end{proposition}
\begin{proof}
  This is another straightforward calculation in
  \GiNaC~\cite{Kisil05b}*{\S~\ref{G-sec:check-conformity}}. 
\end{proof}
The conformal property of the
distance~\eqref{eq:distance-first-ell-hyp}--\eqref{eq:distance-first-par}
from Prop.~\ref{pr:conformity}.\ref{it:conformity-dist} is well-known, of course,
see~\cites{Cnops02a,Yaglom79}. However the same property of
non-symmetric lengths from Prop.~\ref{pr:conformity}.\ref{it:conformity-length-centre} and
\ref{pr:conformity}.\ref{it:conformity-length-foci} could be hardly expected. The smaller
group \(\SL\) (in comparison to all linear-fractional transforms of
\(\Space{R}{2}\)) generates bigger number of conformal metrics, cf.
Rem.~\ref{re:real-line-inv-0}.

The exception of the case \(\rs=0\) from the conformality
in~\ref{pr:conformity}.\ref{it:conformity-length-foci} looks disappointing on the first
glance, especially in the light of the parabolic Cayley transform
considered later in \S~\ref{sec:parab-cayl-transf}. However a detailed
study of algebraic structure invariant under parabolic
rotations~\citelist{\cite{Kisil07a} \cite{Kisil09c}} removes obscurity from this case. Indeed our
Definition~\ref{de:conformal} of conformality heavily depends on the
underlying linear structure in \(\Space{R}{a}\): we measure a distance
between points \(y\) and \(y+ty'\) and intuitively expect that it is
always small for small \(t\). As explained
in~\cite{Kisil07a}*{\S~\ref{W-sec:invar-line-algebra}} the standard
linear structure is incompatible with the parabolic rotations and thus
should be replaced by a more relevant one. More precisely, instead of
limits \(y'\rightarrow y\) along the straight lines towards \(y\) we
need to consider limits along vertical lines, see
Fig.~\ref{fig:concentric-equidist}
and~\cite{Kisil07a}*{Fig.~\ref{W-fig:p-rotations} and
  Rem.~\ref{W-re:conformality}}.  
\begin{proposition}
  \label{pr:parab-conf}
  Let the focal length is given by the
  identity~\eqref{eq:k-focus-point}  with \(\sigma=\rs=0\), e.g.:
  \begin{displaymath}
    l_{f_{\bs}}^2(y, y')  = -\bs p^2-2vp, \qquad \text{where}
    \quad
    p = \frac{(u'-u)^2
    }{2(v'-v)}
  \end{displaymath}
  Then it is conformal in the sense that for any constant \(y=ue_0+ve_1\)
  and \(y'=u'e_0+v'e_1\) with a fixed \(u'\) we have:
  \begin{equation}
    \label{eq:parab-conf-factor}
    \lim_{v'\rightarrow \infty} \frac{l_{f_{\bs}}(g\cdot y,
      g\cdot y')}{l_{f_{\bs}}(y, y')}
    =\frac{1}{(cu+d)^2},\qquad \text{where}\quad
    g=
    \begin{pmatrix}
      a&be_0\\-ce_0&d
    \end{pmatrix}.
  \end{equation}
\end{proposition}
We also revise the parabolic case of conformality in
\S~\ref{sec:infin-conf} with a related definition based on infinitesimal
cycles.
\begin{remark}
  The expressions of
  lengths~\eqref{eq:k-center-point}--\eqref{eq:k-focus-point} are
  generally non-symmetric and this is a price one should pay for its
  non-triviality. All symmetric distances lead to nine two-dimensional
  Cayley-Klein geometries, see~\citelist{\cite{Yaglom79}*{App.~B}
  \cite{HerranzSantander02a} \cite{HerranzSantander02b}}. In the parabolic
  case a symmetric distance of a vector \((u,v)\) is always a function
  of \(u\) alone, cf. Rem.~\ref{re:par-is-not-limit}. For such a
  distance a parabolic unit circle consists from two vertical lines
  (see dotted vertical lines in the second rows on
  Figs.~\ref{fig:orthogonality1} and~\ref{fig:orthogonality2}), which
  is not aesthetically attractive. On the other hand the parabolic
  ``unit cycles'' defined by lengths~\eqref{eq:k-center-point}
  and~\eqref{eq:k-focus-point} are parabolas, which makes the parabolic
  Cayley transform (see Section~\ref{sec:parab-cayl-transf}) very
  natural.
\end{remark}

We can also consider a distance between points in the upper half-plane
which is preserved by M\"obius transformations, see~\cite{Kisil08a}. 
\begin{lemma}
  Let the line element be \(\modulus{dy}^2=du^2-\sigma dv^2\) and
  the ``length of a curve'' is given by the corresponding line
  integral, cf.~\cite{Beardon05a}*{\S~15.2}: 
  \begin{equation}
    \label{eq:curve-length}
    \int_\Gamma \frac{\modulus{dy}}{v}.
  \end{equation}
  Then the length of the curve is preserved under the M\"obius
  transformations. 
\end{lemma}
\begin{proof}
  The proof is based on the following three observations:
  \begin{enumerate}
  \item The line element \(\modulus{dy}^2=du^2-\sigma dv^2\) at the point \(e_1\)
    is invariant under action of the respective fix-group of this
    point (see Lem.~\ref{le:fix-subgroups}).
  \item The fraction \(\frac{\modulus{dy}}{v}\) is invariant under
    action of the \(ax+b\)-group.
  \item M\"obius action of \(\SL\) in each EPH case is generated by
    \(ax+b\) group and the corresponding fix-subgroup, see
    Lem.~\ref{le:fix-ax+b-gen-sl2}. 
  \end{enumerate}
\end{proof}

It is known~\cite{Beardon05a}*{\S~15.2} in the elliptic case that the
curve between two points with the shortest
length~\eqref{eq:curve-length} is an arc of the circle orthogonal to
the real line. M\"obius transformations map such arcs to arcs with the
same property, therefore the length of such arc calculated
in~\eqref{eq:curve-length} is invariant under the M\"obius
transformations. 

Analogously in the hyperbolic case the longest curve between two
points is an arc of hyperbola orthogonal to the real line. However in
the parabolic case there is no curve delivering the shortest
length~\eqref{eq:curve-length}, the infimum is \(u-u'\),
see~\eqref{eq:distance-first-par} and
Fig.~\ref{fig:distances}. However we can still define an invariant
distance in the parabolic case in the following way:
\begin{lemma}[\cite{Kisil08a}]
  Let two points \(w_1\) and \(w_2\) in the upper half-plane are
  linked by an arc of a parabola with zero \(\bs\)-radius. Then the
  length~\eqref{eq:curve-length} along the arc is invariant under
  M\"obius transformations.
\end{lemma}

\subsection{Perpendicularity and Orthogonality}
\label{sec:perp-orth}
In a Euclidean space the shortest distance from a point to a line is
provided by the corresponding perpendicular. Since we have already
defined various distances and lengths we may use them for a definition
of corresponding notions of perpendicularity.

\begin{definition}
  \label{de:perpendicular}
  Let \(l\) be a length or distance.  We say that a vector \(\lvec{AB}\) is
  \emph{\(l\)-perpendicular} to a vector \(\lvec{CD}\) if function
  \(l(\lvec{AB}+\varepsilon \lvec{CD})\) of a variable \(\varepsilon\) has a
  local extremum at \(\varepsilon=0\). This is denoted by
  \(\lvec{AB}\leftthreetimes \lvec{CD}\).
\end{definition}
\begin{remark}
  \begin{enumerate}
  \item Obviously the \(l\)-perpendicularity is not a symmetric notion
    (i.e. \(\lvec{AB}\leftthreetimes \lvec{CD}\) does not imply
    \(\lvec{CD}\leftthreetimes \lvec{AB}\)) similarly to
    f-orthogonality, see subsection~\ref{sec:focal-orthogonality}.
  \item \(l\)-perpendicularity is obviously linear in  \(\lvec{CD}\),
    i.e.  \(\lvec{AB}\leftthreetimes \lvec{CD}\) implies
    \(\lvec{AB}\leftthreetimes r\lvec{CD}\) for any real non-zero
    \(r\). However \(l\)-perpendicularity is not generally linear in
    \(\lvec{AB}\), i.e. \(\lvec{AB}\leftthreetimes \lvec{CD}\) does
    not necessarily imply \(r\lvec{AB}\leftthreetimes \lvec{CD}\).
  \end{enumerate}
\end{remark}
There is the following obvious connection between perpendicularity and
orthogonality.
\begin{lemma}
  \label{le:parp-ortho}
  Let \(\lvec{AB}\) be \(l_{c_{\bs}}\)-perpendicular
  (\(l_{f_{\bs}}\)-perpendicular) to a vector \(\lvec{CD}\). Then the
  flat cycle (straight line) \(AB\), is (s-)orthogonal to the cycle
  \(\cycle{s}{\sigma}\) with centre (focus) at \(A\) passing through
  \(B\). The vector \(\lvec{CD}\) is tangent to \(\cycle{s}{\sigma}\)
  at \(B\).
\end{lemma}
\begin{proof}
  This follows from the relation of centre of (s-)ghost cycle to
  centre (focus) of (s-)orthogonal cycle stated in
  Props.~\ref{pr:ghost-cycle} and~\ref{pr:focal-ortho}
  correspondingly. 
\end{proof}

Consequently the perpendicularity of vectors \(\lvec{AB}\) and
\(\lvec{CD}\) is reduced to the orthogonality of the corresponding
flat cycles only in the cases, when orthogonality itself is reduced to
the local notion at the point of cycles intersections (see
Rem.~\ref{re:local-ortho}).  

Obviously, \(l\)-perpendicularity turns to be the usual orthogonality
in the elliptic case, cf. Lem.~\ref{lem:orthogonal1}.\ref{it:ell-perpendic} below. For two
other cases the description is given as follows:
\begin{lemma}
  \label{le:perp-explicit}  
  Let \(A=(u,v)\) and \(B=(u',v')\). Then
  \begin{enumerate}
  \item \(d\)-perpendicular (in the sense
    of~\eqref{eq:distance-first-ell-hyp}) to \(\lvec{AB}\) in the
    \emph{elliptic} or \emph{hyperbolic} cases is a multiple of the vector
    \begin{displaymath}
      (\sigma (v-v')^3-(u-u')^2 (v+v'(1-2 \sigma \bs)), \bs(u-u')^3-(u-u')(v- v')(-2 v' +(v+v')
      \bs \sigma)),    
    \end{displaymath}
    which for \(\sigma\bs=1\) reduces to the expected
    value \((v-v', \sigma(u-u'))\).
  \item \label{it:d-perp-parab} \(d\)-perpendicular (in the sense
    of~\eqref{eq:distance-first-par}) to \(\lvec{AB}\) in the
    \emph{parabolic} case is \((0, t)\), \(t\in\Space{R}{}\) which
    coincides with the \emph{Galilean orthogonality} defined
    in~\textup{\cite{Yaglom79}*{\S~3}}.
  \item \(l_{c_{\bs}}\)-perpendicular (in the sense
    of~\eqref{eq:k-center-point}) to \(\lvec{AB}\) is a multiple of
    \((\sigma v'-\rs v, u-u')\).
  \item \label{it:focal-perpendicularity}
    \(l_{f_{\bs}}\)-perpendicular (in the sense
    of~\eqref{eq:k-focus-point}) to \(\lvec{AB}\) is a multiple of
    \((\sigma v'+p, u-u')\), where \(p\) is defined
    either by~\eqref{eq:focal-length} or
    by~\eqref{eq:parab-focal-length} for corresponding values of
    \(\rs\). 
  \end{enumerate}
\end{lemma}
\begin{proof}
  The perpendiculars are calculated by
  \GiNaC~\cite{Kisil05b}*{\S~\ref{G-sec:calc-perp}}.
\end{proof}
It is worth to have an idea about different types of perpendicularity
in the terms of the standard Euclidean geometry. Here are some
examples.
\begin{lemma} 
  \label{lem:orthogonal1}
  Let \(\lvec{AB}=u e_0+v e_1\) and
  \(\lvec{CD}=u'e_0+v'e_1\), then:
  \begin{enumerate} 
    \renewcommand{\theenumi}{(\ephname{enumi})}
  \item \label{it:ell-perpendic}
    In the elliptic case the \(d\)-perpendicularity for \(\bs=-1\)
    means that \(\lvec{AB}\) and \(\lvec{CD}\) form a right angle, or
    analytically \(u u'+v v'=0\).
  \item \label{item:par-orthogon} In the parabolic case
    the \(l_{f_{\bs}}\)-perpendicularity for \(\bs=1\) means that
    \(\lvec{AB}\) bisect the angle between \(\lvec{CD}\) and the
    vertical direction or analytically:
    \begin{equation}
      \label{eq:par-s-perp}
      u'u-v'p=u'u-v'(\sqrt{u^2+v^2}-v)=0,
    \end{equation}
    where \(p\) is the focal length~\eqref{eq:focal-length}
  \item In the hyperbolic case the \(d\)-perpendicularity for
    \(\bs=-1\) means that the angles between \(\lvec{AB}\) and
    \(\lvec{CD}\) are bisected by lines parallel to \(u=\pm v\), or
    analytically \(u' u-v' v=0\).
  \end{enumerate}
\end{lemma}
\begin{remark}
  \label{re:par-is-not-limit}
  If one attempts to devise a parabolic length as a limit or an
  intermediate case for the elliptic \(l_e=u^2+v^2\) and hyperbolic
  \(l_p=u^2-v^2\) lengths then the only possible guess is
  \(l'_p=u^2\)~\eqref{eq:distance-first-par}, which is too trivial for
  an interesting geometry.
  
  Similarly the only orthogonality conditions linking the
  elliptic \(u_1 u_2+v_1 v_2=0\) and the hyperbolic \(u_1 u_2-v_1
  v_2=0\) cases seems to be \(u_1 u_2=0\) (see~\cite{Yaglom79}*{\S~3}
  and Lem.~\ref{le:perp-explicit}.\ref{it:d-perp-parab}), which is again too trivial.  This
  support our Remark~\ref{rem:prejudice}.\ref{it:par-is-limit}.
\end{remark}

\section{Invariants of Infinitesimal Scale}
\label{sec:invar-infin-scale}

Although parabolic zero-radius cycles defined
in~\ref{de:zero-radius-cycle} do not satisfy our expectations for
``zero-radius'' but they are often technically suitable for the same
purposes as elliptic and hyperbolic ones. Yet we may want to find
something which fits better for our intuition on ``zero sized'' object.
Here we present an approach based on non-Archimedean (non-standard)
analysis~\cites{Devis77,Uspenskii88}.

\subsection{Infinitesimal Radius Cycles}
\label{sec:zero-length-cycles}
\begin{figure}[htbp]
  \centering
  \quad(a)\includegraphics[scale=1.2]{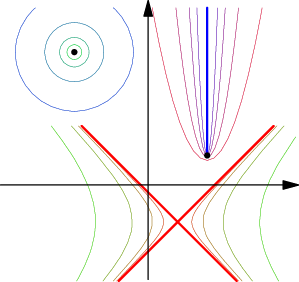}\hfill
  {(b)\includegraphics[scale=1.2]{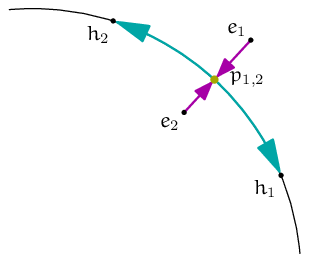}\quad}
  \caption[Zero-radius cycles and ``phase'' transition]{(a)
    Zero-radius cycles in elliptic (black point) and hyperbolic (the
    red light cone). Infinitesimal radius parabolic
    cycle is the blue vertical ray starting at the focus.\\
    (b) Elliptic-parabolic-hyperbolic phase transition between fixed
    points of the subgroup \(K\).}
  \label{fig:infinites-radius}
\end{figure}

Let \(\varepsilon\) be a positive infinitesimal number, i.e. \(0 <
n\varepsilon <1\) for any \(n\in\Space{N}{}\)~\cites{Devis77,Uspenskii88}.
\begin{definition}
  A cycle \(\cycle{s}{\bs}\) such that \(\det\cycle{s}{\bs}\) is an
  infinitesimal number is called
  \emph{infinitesimal radius cycle}. 
\end{definition}
\begin{lemma}
  Let \(\bs\) and \(\rs\) be two metric signs and let a point
  \((u_0,v_0)\in\Space{R}{p}\) with \(v_0>0\).  Consider a cycle
  \(\cycle{s}{\bs}\) defined by
  \begin{equation}
    \label{eq:inf-cycle}
    \cycle{s}{\bs}=(1, u_0, n, u_0^2+2 n v_0
    -\rs n^2),
  \end{equation}
  where
  \begin{equation}
    \label{eq:inf-cycle-n-value}
    n=\left\{
      \begin{array}{ll}
        \displaystyle    \frac{v_0-\sqrt{v_0^2-(\rs-\bs)\varepsilon^2}}{\rs-\bs},&\text{ if } \rs\neq \bs;\\
       \displaystyle \frac{\varepsilon^2}{2v_0},&\text{ if } \rs=\bs.
      \end{array}\right.
  \end{equation}
  Then 
  \begin{enumerate}
  \item The point \((u_0, v_0)\) is \(\rs\)-focus of the cycle.
  \item The square of \(\bs\)-radius is exactly \(-\varepsilon^2\),
    i.e.~\eqref{eq:inf-cycle} defines an infinitesimal radius cycle. 
  \item The focal length of the cycle is an infinitesimal
    number of order \(\varepsilon^2\).
  \end{enumerate}
\end{lemma}
\begin{proof} 
  The cycle~\eqref{eq:inf-cycle} has the squared \(\bs\)-radius equal
  to \(-\varepsilon^2\) if \(n\) is a root to the equation:
  \begin{displaymath}
    (\rs-\bs)n^2-2v_0n+\varepsilon^2=0.
  \end{displaymath}
  Moreover only the root from~\eqref{eq:inf-cycle-n-value} of the
  quadratic case gives an infinitesimal focal length.  This also is
  supported by calculations done in \GiNaC,
  see~\cite{Kisil05b}*{\S~\ref{G-sec:basic-prop-infin}}.
\end{proof}
The graph of cycle~\eqref{eq:inf-cycle} in the parabolic space drawn
at the scale of real numbers looks like a vertical ray started at its
focus, see Fig.~\ref{fig:infinites-radius}(a), due to
the following Lemma.
\begin{lemma}{\cite{Kisil05b}*{\S~\ref{G-sec:basic-prop-infin}}}
  Infinitesimal cycle~\eqref{eq:inf-cycle} consists of points, which are infinitesimally
  close (in the sense of length from focus~\eqref{eq:k-focus-point})
  to  its focus \(F=(u_0, v_0)\):
  \begin{equation}
    \label{eq:inf-cycle-points}
    (u_0+\varepsilon u, v_0+v_0u^2+((\bs-\rs) u^2-\rs)
    \frac{\varepsilon^2}{4v_0}+O(\varepsilon^3)).
  \end{equation}
\end{lemma}
Note that points below of \(F\) (in the ordinary scale) are not
infinitesimally close to \(F\) in the sense of
length~\eqref{eq:k-focus-point}, but are in the sense of
distance~\eqref{eq:distance-first-par}.
Figure~\ref{fig:infinites-radius}(a) shows elliptic, hyperbolic
\emph{concentric} and parabolic \emph{confocal} cycles of decreasing
radii which shrink to the corresponding infinitesimal radius cycles.

It is easy to see that infinitesimal radius cycles has properties
similar to zero-radius ones, cf. Lemma~\ref{le:moeb-conj-z-cycle}. 
\begin{lemma}
  \label{le:infinites-cycle-inv}
  The image of \(\SL\)-action on an infinitesimal radius
  cycle~\eqref{eq:inf-cycle} by
  conjugation~\eqref{eq:cycle-transform-short} is an infinitesimal
  radius cycle of the same order. 

  Image of an infinitesimal cycle under cycle conjugation is an
  infinitesimal cycle of the same or lesser order. 
  \end{lemma}
\begin{proof}
  These are calculations done in \GiNaC,
  see~\cite{Kisil05b}*{\S~\ref{G-sec:mobi-transf-infin}}. 
\end{proof}

The consideration of infinitesimal numbers in the elliptic and
hyperbolic case should not bring any advantages since the (leading)
quadratic terms in these cases are non-zero. However 
non-Archimedean numbers in the parabolic case provide a more
intuitive and efficient presentation. For example zero-radius cycles
are not helpful for the parabolic Cayley transform (see
subsection~\ref{sec:parab-cayl-transf}) but infinitesimal cycles are
their successful replacements. 

The second part of the following result is a useful substitution for
Lem.~\ref{le:zero-radius-otho}. 
\begin{lemma}
  \label{le:infinitesimal-ortho} 
  Let \(\cycle{s}{\bs}\) be the
  infinitesimal cycle~\eqref{eq:inf-cycle} and \(\cycle[\breve]{s}{\bs}=(k,l,n,m)\) be a
  generic cycle. Then
  \begin{enumerate}
  \item The orthogonality condition~\eqref{eq:orthogonality-first}
    \(\cycle{s}{\bs}\perp \cycle[\breve]{s}{\bs}\) and the
    f-orthogonality~\eqref{eq:f-orthog-def} \(\cycle[\breve]{s}{\bs}
    \sperp \cycle{s}{\bs}\) both are given by:
    \begin{displaymath}
      ku_0^2-2lu_0+m=O(\varepsilon).
    \end{displaymath}
    In other words the cycle \(\cycle[\breve]{s}{\bs}\) has root
    \(u_0\) in the parabolic space. 
  \item \label{it:infinitesimal-f-ortho}
    The f-orthogonality~\eqref{eq:f-orthog-def} \(\cycle{s}{\bs}
    \sperp \cycle[\breve]{s}{\bs}\) is given by:
    \begin{equation}
      \label{eq:f-orthog-infinites}
      ku_0^2-2lu_0-2nv_0+m=O(\varepsilon).
    \end{equation}
    In other words the cycle \(\cycle[\breve]{s}{\bs}\) passes focus
    \((u_0,v_0)\) of the infinitesimal cycle in the parabolic space. 
  \end{enumerate}
\end{lemma}
\begin{proof}
  These are \GiNaC\ calculations~\cite{Kisil05b}*{\S~\ref{G-sec:orth-with-infin}}.
\end{proof}
It is interesting to note that the exotic f-orthogonality became
warranted replacement of the usual one for the infinitesimal cycles. 

\subsection{Infinitesimal Conformality}
\label{sec:infin-conf}
An intuitive idea of conformal maps, which is oftenly provided in the
complex analysis textbooks for illustration purposes, is ``they send
small circles into small circles with respective centres''. Using
infinitesimal cycles one can turn it into a precise definition.
\begin{definition}
  \label{de:infinites-conform}
  A map of a region of \(\Space{R}{a}\) to another region is
  \(l\)-infinitesimally conformal for a length \(l\) (in the sense of
  Defn.~\ref{de:length}) if for any \(l\)-infinitesimal cycle: 
  \begin{enumerate}
  \item Its image is an \(l\)-infinitesimal cycle of the same order;
  \item The image of its centre/focus is displaced from the centre/focus of its
    image by an infinitesimal number of a greater order than its radius.
  \end{enumerate}
\end{definition}
\begin{remark}
  Note that in comparison with Defn.~\ref{de:conformal} we now work
  ``in the opposite direction'': former we had the fixed group of
  motions and looked for corresponding conformal lengths/distances,
  now we take a distance/length (encoded in the infinitesimally equidistant
  cycle) and check which motions respect it.
\end{remark}

Natural conformalities for lengths from centre in the elliptic and
parabolic cases are already well studied.  Thus we are mostly
interested here in conformality in the parabolic case, where lengths
from focus are better suited. The image of an infinitesimal
cycle~\eqref{eq:inf-cycle} under \(\SL\)-action is a cycle, moreover
its is again an infinitesimal cycle of the same order by
Lemma~\ref{le:infinites-cycle-inv}. This provides the first condition
of Defn.~\ref{de:infinites-conform}. The second part is delivered by
the following statement:
\begin{proposition}
  \label{pr:infinit-conf}
  Let \(\cycle[\breve]{s}{\bs}\) be the image under \(g\in\SL\) of an infinitesimal
  cycle \(\cycle{s}{\bs}\) from~\eqref{eq:inf-cycle}.
  Then \(\rs\)-focus of \(\cycle[\breve]{s}{\bs}\) is  displaced from
  \(g(u_0,v_0)\) by infinitesimals of order \(\varepsilon^2\) (while
  both cycles have \(\bs\)-radius of order \(\varepsilon\)).

  Consequently \(\SL\)-action is infinitesimally conformal in the
  sense of Defn.~\ref{de:infinites-conform} with
  respect to the length from focus (Defn.~\ref{de:length}) for all
  combinations of \(\sigma\), \(\bs\) and \(\rs\).
\end{proposition}
\begin{proof}
    These are \GiNaC\ calculations~\cite{Kisil05b}*{\S~\ref{G-sec:mobi-transf-infin}}.
\end{proof}
Infinitesimal conformality seems intuitively to be close to
Defn.~\ref{de:conformal}.  Thus it is desirable to give a reason for
the absence of exclusion clauses in Prop.~\ref{pr:infinit-conf} in
comparison to Prop.~\ref{pr:conformity}.\ref{it:conformity-length-foci}. As shows
calculations~\cite{Kisil05b}*{\S~\ref{G-sec:check-conformity}} the
limit~\eqref{eq:conformal-lim} at point \(y_0=u_0e_0+v_0e_1\) do exist but
depends from the direction \(y=ue_0+ve_1\):
\begin{equation} 
  \label{eq:focal-length-factor}
  \lim_{t\rightarrow 0} \frac{d(g\cdot y_0, g\cdot(y_0+ty))}{d(y_0,
    y_0+ty)}=\frac{1}{(d+cu_0)^2+\sigma c^2 v_0^2 -2 K c v_0(d+c u_0)},
\end{equation}
where \(\displaystyle K=\frac{u}{v}\) and \(g=
\begin{pmatrix}
  a&b\\c&d
\end{pmatrix}\). However if we consider points
~\eqref{eq:inf-cycle-points} of the infinitesimal cycle then
\(\displaystyle K=\frac{\varepsilon u}{v_0 u^2}=
\frac{\varepsilon}{v_0 u}\). Thus the value of the
limit~\eqref{eq:focal-length-factor} at the infinitesimal scale is
independent from \(y=ue_0+ve_1\). It also coincides (up to an
infinitesimal number) with the value in~\eqref{eq:parab-conf-factor}. 

\begin{remark}
  \label{re:non-standard}
  There is another connection between parabolic function theory and
  non-standard analysis. As was mentioned in
  \S~\ref{sec:ellipt-parab-hyperb}, the Clifford algebra \(\Cliff{p}\)
  corresponds to the set of \emph{dual numbers} \(u+\varepsilon v\)
  with \(\varepsilon^2=0\)~\cite{Yaglom79}*{Supl.~C}. On the other hand
  we may consider the set of numbers \(u+\varepsilon v\) within the
  non-standard analysis, with \(\varepsilon\) being an infinitesimal.
  In this case \(\varepsilon^2\) is a higher order infinitesimal than
  \(\varepsilon\) and effectively can be treated as \(0\) at
  infinitesimal scale of \(\varepsilon\), i.e. we again get the dual
  numbers condition \(\varepsilon^2=0\). This explains why many results
  of differential calculus can be naturally deduced within dual numbers
  framework~\cite{CatoniCannataNichelatti04}.
\end{remark}
Infinitesimal cycles are also a convenient tool for calculations of
invariant measures, Jacobians, etc.

\section{Global Properties}
\label{sec:global-properties}
So far we were interested in individual properties of cycles and local
properties of the point space. Now we describe some global properties
which are related to the set of cycles as the whole.

\subsection[Compactification of R2]{Compactification of $\Space{R}{\sigma}$}
\label{sec:comp-spra}
Giving Definition~\ref{de:cycle-2-matrix} of maps \(Q\) and \(M\) we
did not consider properly their domains and ranges. For example, the
image of \((0,0,0,1)\in\Space{P}{3}\), which is transformed by \(Q\)
to the equation \(1=0\), is not a valid conic section in
\(\Space{R}{\sigma}\). We also did not investigate yet accurately
singular points of the M\"obius map~\eqref{eq:moebius-def}. It
turns out that both questions are connected.

One of the standard approaches~\cite{Olver99}*{\S~1} to deal with
singularities of the M\"obius map is to consider projective coordinates
on the plane. Since we have already a projective space of cycles, we
may use it as a model for compactification which is even more
appropriate. The identification of points with zero-radius cycles
plays an important r\^ole here.
\begin{definition}
  The only irregular point \((0,0,0,1)\in\Space{P}{3}\) of the map
  \(Q\) is called \emph{zero-radius cycle at infinity} and denoted by
  \(\zcycle{}{\infty}\). 
\end{definition}
The following results are easily obtained by direct calculations even
without a computer:
\begin{lemma}
  \begin{enumerate}
  \item \(\zcycle{}{\infty}\) is the image of the zero-radius cycle
    \(\zcycle{}{(0,0)}=(1,0,0,0)\) at the origin under reflection
    (inversion) into the unit cycle \((1, 0,0,-1)\), see blue cycles
    in Fig.~\ref{fig:inversions}(b)-(d).
  \item The following statements are equivalent
    \begin{enumerate}
    \item A point \((u,v)\in\Space{R}{\sigma}\) belongs to the
      zero-radius cycle \(\zcycle{}{(0,0)}\) centred at the origin;
     \item The zero-radius cycle \(\zcycle{}{(u,v)}\) is \(\sigma\)-orthogonal to 
      zero-radius cycle \(\zcycle{}{(0,0)}\);
    \item The inversion \(z\mapsto \frac{1}{z}\) in the unit cycle is
      singular in the point \((u,v)\);
    \item The image of \(\zcycle{}{(u,v)}\) under inversion in the unit
      cycle is orthogonal to  \(\zcycle{}{\infty}\). 
    \end{enumerate}
    If any from the above is true we also say that image of \((u,v)\) under
    inversion in the unit cycle belongs to zero-radius cycle at infinity.
  \end{enumerate}
\end{lemma}

\begin{figure}[htbp]
  \centering
  \includegraphics[scale=1.3]
 {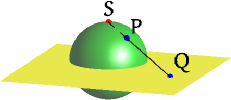} \hfill
  \includegraphics[scale=1.3]
{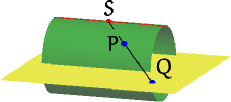} \hfill
  \includegraphics[scale=1.3]
{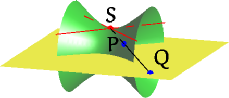}\\
(a)\hfill(b)\hfill(c)\hfill
  \caption[Compactification and stereographic projections]{Compactification of \(\Space{R}{\sigma}\) and stereographic projections.}
  \label{fig:compactifications}
\end{figure}

In the elliptic case the compactification is done by addition to
\(\Space{R}{e}\) a point \(\infty\) at infinity, which is the elliptic
zero-radius cycle. However in the parabolic and hyperbolic cases the
singularity of the M\"obius transform is not localised in a single
point---the denominator is a zero divisor for the whole zero-radius
cycle.  Thus in each EPH case the correct compactification is made by
the union \(\Space{R}{\sigma}\cup\zcycle{}{\infty}\).

 It is common to identify the compactification \(\DSpace{R}{e}\) of the
 space \(\Space{R}{e}\) with a Riemann sphere.
 This model can be visualised by the stereographic
 projection, see~\cite{BergerII}*{\S~18.1.4} and
 Fig.~\ref{fig:compactifications}(a). A similar model can be 
 provided for the parabolic and hyperbolic spaces as
 well, see~\cite{HerranzSantander02b} and
 Fig.~\ref{fig:compactifications}(b),(c). Indeed the space 
 \(\Space{R}{\sigma}\) can be identified with a corresponding surface
 of the constant curvature: the sphere (\(\sigma=-1\)), the cylinder
 (\(\sigma=0\)), or the one-sheet hyperboloid (\(\sigma=1\)). The map
 of a surface to \(\Space{R}{\sigma}\) is given by the polar
 projection, see~\cite{HerranzSantander02b}*{Fig.~1} and
 Fig.~\ref{fig:compactifications}(a)-(c). These
 surfaces provide ``compact'' model of the corresponding
 \(\Space{R}{\sigma}\) in the sense that M\"obius transformations
 which are lifted from \(\Space{R}{\sigma}\) by the projection are
 not singular on these surfaces. 

However the hyperbolic case has its own caveats which may be easily
oversight as in the paper cited above, for example. A compactification of
the hyperbolic space \(\Space{R}{h}\) by a light cone (which the
hyperbolic zero-radius cycle) at infinity will indeed
produce a closed M\"obius invariant object. However it will not be
satisfactory for some other reasons explained in the next subsection.

\subsection{(Non)-Invariance of The Upper Half-Plane}
\label{sec:invar-upper-half}
The important difference between the hyperbolic case and the two
others is that
\begin{lemma}
  In the elliptic and parabolic cases the upper halfplane in
  \(\Space{R}{\sigma}\) is preserved by M\"obius transformations from
  \(\SL\). However in the hyperbolic case any point \((u,v)\) with
  \(v>0\) can be mapped to an arbitrary point \((u',v')\) with
  \(v'\neq 0\).
\end{lemma}
This is illustrated by Fig.~\ref{fig:k-orbit-sect}: any cone from the
family~\eqref{eq:cones-fam} is intersecting the both planes \(EE'\) and
\(PP'\) over a connected curve, however intersection with the 
plane \(HH'\) has two branches.

The lack of invariance in the hyperbolic case has many important
consequences in seemingly different areas, for example:
\begin{figure}[htbp]
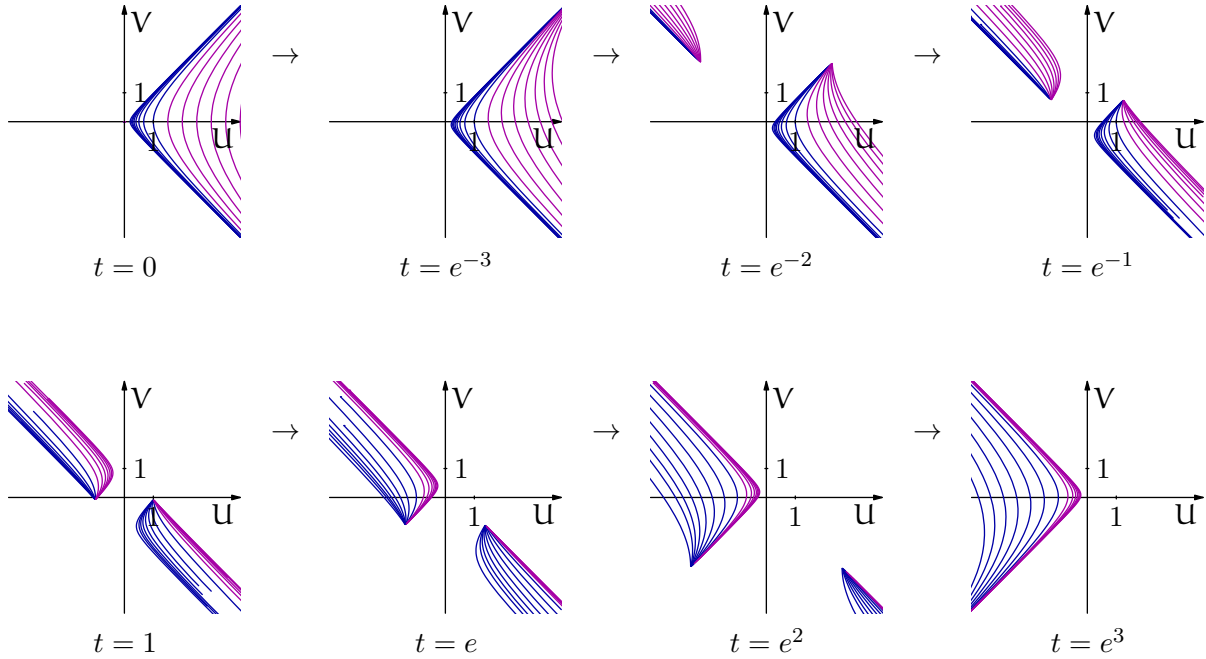

  \parbox[t]{.2\textwidth}{
    \begin{center}
      \includegraphics[scale=1.1]{parabolic0.20}\\
      \(t=0\)
    \end{center}
  } \hfill\raisebox{-1.2cm}{\(\to\)} \hfill
  \parbox[t]{.2\textwidth}{
    \begin{center}
      \includegraphics[scale=1.1]{parabolic0.21}\\
      \(t=e^{-3}\)
    \end{center}
  }   \hfill \raisebox{-1.2cm}{\(\to\)} \hfill
  \parbox[t]{.2\textwidth}{
    \begin{center}
      \includegraphics[scale=1.1]{parabolic0.22}\\
      \(t=e^{-2}\)
    \end{center}
  } \hfill \raisebox{-1.2cm}{\(\to\)} \hfill
    \parbox[t]{.2\textwidth}{
    \begin{center}
      \includegraphics[scale=1.1]{parabolic0.23}\\
      \(t=e^{-1}\)
    \end{center}
  }\\[5mm]
  \parbox[t]{.2\textwidth}{
    \begin{center}
      \includegraphics[scale=1.1]{parabolic0.24}\\
      \(t=1\)
    \end{center}
  }  \hfill \raisebox{-1.2cm}{\(\to\)} \hfill
  \parbox[t]{.2\textwidth}{
    \begin{center}
      \includegraphics[scale=1.1]{parabolic0.25}\\
      \(t=e\)
    \end{center}
  } \hfill \raisebox{-1.2cm}{\(\to\)} \hfill
  \parbox[t]{.2\textwidth}{
    \begin{center}
      \includegraphics[scale=1.1]{parabolic0.26}\\
      \(t=e^2\)
    \end{center}
  } \hfill \raisebox{-1.2cm}{\(\to\)} \hfill
  \parbox[t]{.2\textwidth}{
    \begin{center}
      \includegraphics[scale=1.1]{parabolic0.27}\\
      \(t=e^3\)
    \end{center}
  }
  \caption[Continuous transformation from future
    to the past]{Eight frames from a continuous transformation from future
    to the past parts   of the light cone.}
  \label{fig:future-to-past}
\end{figure}
\begin{description}
\item[\textbf{Geometry}] \(\Space{R}{h}\) is not split by the real
  axis into two disjoint pieces: there is a continuous path (through
  the light cone at infinity) from the upper half-plane to the lower
  which does not cross the real axis (see \(\sin\)-like curve joined
  two sheets of the hyperbola in
  Fig.~\ref{fig:hyp-upper-half-plane}(a)).
\item[\textbf{Physics}] There is no M\"obius invariant way to separate
  ``past'' and ``future'' parts of the light cone~\cite{Segal76}, i.e.
  there is a continuous family of M\"obius transformations reversing
  the arrow of time. For example, the family of matrices \(
  \begin{pmatrix}
    1&-te_1\\te_1&1
  \end{pmatrix}\), \(t\in [0,\infty)\) provides such a transformation.
  Fig.~\ref{fig:future-to-past} illustrates this by corresponding
  images for eight subsequent values of \(t\).
\item[\textbf{Analysis}] There is no a possibility to split
  \(\FSpace{L}{2}(\Space{R}{})\) space of function into a direct sum
  of the Hardy type space of functions having an analytic extension into
  the upper half-plane and its non-trivial complement, i.e. any function
  from \(\FSpace{L}{2}(\Space{R}{})\) has an ``analytic extension''
  into the upper half-plane in the sense of hyperbolic function
  theory,  see~\cite{Kisil97c}.
\end{description}
\begin{figure}[htbp]
  \centering
     (a)\includegraphics[scale=1.1]{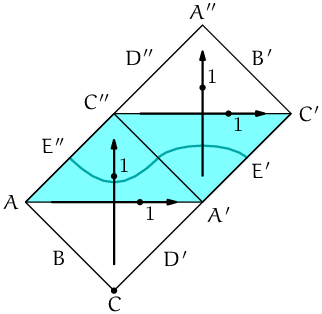} \hspace{.1\textwidth}
    (b)\includegraphics[scale=1.1
]{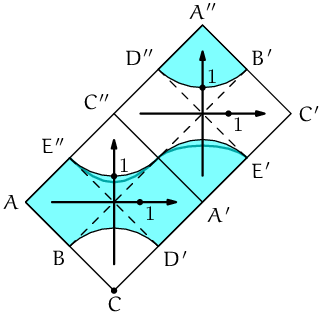}
  \caption[Hyperbolic objects in the double cover]{Hyperbolic objects
    in the double cover of \(\Space{R}{h}\):\\ 
  (a) the ``upper'' half-plane;\qquad (b) the unit circle.}
  \label{fig:hyp-upper-half-plane}
\end{figure} 
All the above problems can be resolved in the following
way~\cite{Kisil97c}*{\S~A.3}.  
We take two copies \(\Space[+]{R}{h}\)
and \(\Space[-]{R}{h}\) of \( \Space{R}{h} \), depicted by the squares
\(ACA'C''\) and \(A'C'A''C''\) in Fig.~\ref{fig:hyp-upper-half-plane}
correspondingly. The boundaries of these squares are light cones at
infinity and we glue \(\Space[+]{R}{h}\) and \(\Space[-]{R}{h}\) in
such a way that the construction is invariant under the natural action
of the M\"obius transformation.  That is achieved if the letters
\(A\), \(B\), \(C\), \(D\), \(E\) in
Fig.~\ref{fig:hyp-upper-half-plane} are identified regardless of the
number of primes attached to them. 
The corresponding model through a stereographic projection is presented
on Fig.~\ref{fig:compact-2}, compare with Fig.~\ref{fig:compactifications}(c).
\begin{figure}[htbp]
  \centering
  \includegraphics[scale=1.2]
 {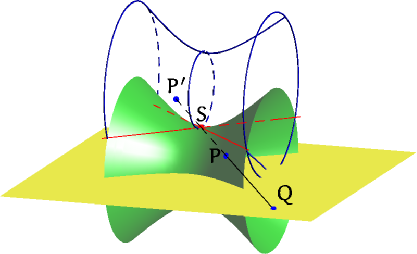}
  \caption[Double cover of the hyperbolic space]{Double cover of the
    hyperbolic space, cf. Fig.~\ref{fig:compactifications}(c). The
  second hyperboloid is shown as a  blue skeleton. It is attached to the
  first one along the light cone at infinity, which is represented by
  two red lines.}  
  \label{fig:compact-2}
\end{figure}

This aggregate denoted by \(\TSpace{R}{h}\) is a two-fold cover of
\(\Space{R}{h}\). The hyperbolic ``upper'' half-plane in
\(\TSpace{R}{h}\) consists of the upper halfplane in
\(\Space[+]{R}{h}\) and the lower one in \(\Space[-]{R}{h}\).  A
similar conformally invariant two-fold cover of the Minkowski
space-time was constructed in~\cite{Segal76}*{\S~III.4} in connection
with the red shift problem in extragalactic astronomy. 

\begin{remark}
  \label{rem:two-cver}
  \begin{enumerate}
  \item \label{it:hyp-object} The hyperbolic orbit of the \(K\)
    subgroup in the \(\TSpace{R}{h}\) consists of two branches of the
    hyperbola passing through \((0,v)\) in \(\Space[+]{R}{h}\) and
    \((0,-v^{-1})\) in \(\Space[-]{R}{h}\), see
    Fig.~\ref{fig:hyp-upper-half-plane}.  If we watch the rotation of
    a straight line generating a cone~\eqref{eq:cones-fam} then its
    intersection with the plane \(HH'\) on
    Fig.~\ref{fig:k-orbit-sect}(d) will draw the both branches. 
    As mentioned in Rem.~\ref{rem:orbits2}.\ref{it:hyp-orbit} they have the same
    focal length.
  \item The ``upper'' halfplane is bounded by two disjoint ``real''
    axes denoted by \(AA'\) and \(C'C''\) in
    Fig.~\ref{fig:hyp-upper-half-plane}.
  \end{enumerate}
\end{remark}

For the hyperbolic Cayley transform in the next subsection we 
need the conformal version of the hyperbolic unit disk.
 We define it in \( \TSpace{R}{h} \) as follows:
\begin{eqnarray*}
  \TSpace{D}{}&=&\{(ue_0+ve_1) \such u^2-v^2>-1,\ u\in \Space[+]{R}{h} \}\\
  &&{}\cup \{(ue_0+ve_1) \such u^2-v^2<-1,\ u\in \Space[-]{R}{h} \}.
\end{eqnarray*}
It can be shown that \( \TSpace{D}{}\) is conformally invariant and 
has a boundary \(\TSpace{T}{}\)---two copies of the unit circles in 
\(\Space[+]{R}{h}\) and \(\Space[-]{R}{h}\). We call \(\TSpace{T}{}\) 
the \emph{(conformal) unit circle} in \(\Space{R}{h}\). 
Fig.~\ref{fig:hyp-upper-half-plane}(b) illustrates
the geometry of the conformal unit disk in \(\TSpace{R}{h}\) in
comparison with the ``upper'' half-plane.

\section{The Cayley Transform and the Unit Cycle}
\label{sec:unit-circles}

The upper half-plane is the universal starting point for an analytic
function theory of any EPH type. However universal models are rarely
best suited to particular circumstances. For many reasons it is more
convenient to consider analytic functions in the unit disk rather than
in the upper half-plane, although both theories are completely
isomorphic, of course. This isomorphism is delivered by the
\emph{Cayley transform}. Its drawback is that there is no a
``universal unit disk'', in each EPH case we obtain something
specific from the same upper half-plane.

\subsection{Elliptic and Hyperbolic Cayley Transforms}
\label{sec:ellipt-hyperb-caley}

In the elliptic and hyperbolic cases~\cite{Kisil97c} the Cayley transform is given by the
matrix \(C=\begin{pmatrix} 1 &  -e_1 \\ \sigma e_1 &
  1 \end{pmatrix}\), where
\(\sigma=e_1^2\)~\eqref{eq:alg-mult} and \(\det C =2\). It can be applied as the
M\"obius transformation
\begin{equation}
  \label{eq:cayley-points}
  \begin{pmatrix}
  1 & -e_1 \\ \sigma e_1 & 1
\end{pmatrix}:\  w=(ue_0+ve_1) \mapsto Cw=
\frac{(ue_0+ve_1)- e_1}{ \sigma e_1 (ue_0+ve_1)+1}
\end{equation}
to a point \((ue_0+ve_1)\in\Space{R}{\sigma}\).  Alternatively it acts by conjugation
\(g_C=\frac{1}{2}CgC^{-1}\) on an element \(g\in\SL\):
\begin{equation}
  \label{eq:cayley-matr}
  g_C= \frac{1}{2}
  \begin{pmatrix}
    1 & - e_1 \\ \sigma e_1 & 1
  \end{pmatrix}
  \begin{pmatrix}
    a & be_0 \\ -c e_0 & d
  \end{pmatrix}
  \begin{pmatrix}
    1 &  e_1 \\  -\sigma e_1 & 1
  \end{pmatrix} 
\end{equation} The connection between the two
forms~\eqref{eq:cayley-points} and~\eqref{eq:cayley-matr} of the Cayley
transform is given by \(g_C Cw= C(gw)\), i.e. \(C\) intertwines the
actions of \(g\) and \(g_C\).

The Cayley transform \((u'e_0+v'e_1)= C(ue_0+ve_1)\) in the elliptic
case is very important~\citelist{\cite{Lang85}*{\S~IX.3}
\cite{MTaylor86}*{Ch.~8, (1.12)}} both for complex analysis and representation theory
of \(\SL\).  The transformation \(g\mapsto
g_C\)~\eqref{eq:cayley-matr} is an isomorphism of the groups \(\SL\)
and \(\mathrm{SU}(1,1)\) namely in \(\Cliff{e}\) we have
\begin{equation}
  \label{eq:cayley-elliptic}
  g_C= 
  \frac{1}{2}\begin{pmatrix}
    f & h \\ -h & f
\end{pmatrix}, \textrm{ with }
f= (a+d)-(b-c)e_1e_0 \textrm{ and } h= (a-d)e_1+(b+c)e_0.
\end{equation}

Under the map \(\Space{R}{e}\to \Space{C}{}\)~\eqref{eq:complexification}
this matrix becomes \(
\begin{pmatrix}
  \alpha & \beta \\ \bar{\beta} & \bar{\alpha}
\end{pmatrix}
\), i.e. the standard form of elements of
\(\mathrm{SU}(1,1)\)~\citelist{\cite{Lang85}*{\S~IX.1}
\cite{MTaylor86}*{Ch.~8, (1.11)}}.  

The images of elliptic actions of subgroups \(A\), \(N\), \(K\) are
given in Fig.~\ref{fig:unit-disks}(\(E\)).  The types of orbits can be
easily distinguished by the number of \emph{fixed points on the
  boundary}: two, one and zero correspondingly. Although a closer
inspection demonstrate that there are always two fixed points, either:
\begin{itemize}
\renewcommand{\theenumi}{(\ephname{enumi})}
\item one strictly inside and one strictly outside of the unit circle;
  or
\item one double fixed point on the unit circle; or
\item two different fixed points exactly on the circle. 
\end{itemize}
Consideration of Figure~\ref{fig:infinites-radius}(b) shows that the parabolic
subgroup \(N\) is like a phase transition between the elliptic
subgroup \(K\) and hyperbolic \(A\), cf.~\eqref{eq:eph-class}.

In some sense the elliptic Cayley transform swaps complexities: by
contract to the upper half-plane the \(K\)-action is now simple but
\(A\) and \(N\) are not. The simplicity of \(K\) orbits is explained
by diagonalisation of matrices:
\begin{equation}
  \label{eq:k-diagonalisation}
   \frac{1}{2}
   \begin{pmatrix}
    1 & - e_1 \\ - e_1 & 1
  \end{pmatrix}
  \begin{pmatrix}
    \cos \phi  & -e_0\sin \phi  \\ - e_0\sin \phi & \cos \phi
  \end{pmatrix}
  \begin{pmatrix}
    1 & e_1 \\  e_1 & 1
  \end{pmatrix} 
  =
  \begin{pmatrix}
    \rme^{\rmi \phi} & 0\\ 0 & \rme^{\rmi \phi}
  \end{pmatrix},
\end{equation}
where \(\rmi =e_1 e_0\) behaves as the complex imaginary unit, i.e. \(\rmi^2=-1\).

A hyperbolic version of the Cayley transform was used
in~\cite{Kisil97c}. The above formula~\eqref{eq:cayley-matr} in
\(\Space{R}{h}\) becomes as follows: 
\begin{equation}
  \label{eq:cayley-hyp}
   g_C= 
  \frac{1}{2}\begin{pmatrix}
    f & h \\ h & f
\end{pmatrix}, \textrm{ with }
f =a+d-(b+c)e_1e_0  \textrm{ and }
h = (d-a) e_1+ (b-c)e_0,
\end{equation}
with some subtle differences in comparison
with~\eqref{eq:cayley-elliptic}.  The corresponding \(A\), \(N\) and
\(K\) orbits are given on Fig.~\ref{fig:unit-disks}(\(H\)). However
there is an important distinction between the elliptic and hyperbolic
cases similar to one discussed in
subsection~\ref{sec:invar-upper-half}.

\begin{figure}[htbp]
  \centering
\raisebox{1.5cm}{(\(E\))}  \quad
\includegraphics[scale=.75
]{parabolic0.40}\hfill
  \includegraphics[scale=.75
]{parabolic0.41}\hfill
  \includegraphics[scale=.75
]{parabolic0.42}\quad\\[2mm]
\raisebox{1cm}{(\(P_e\))}
\includegraphics[scale=.75
]{parabolic0.50}\hfill
  \includegraphics[scale=.75
]{parabolic0.51}\hfill
  \includegraphics[scale=.75
]{parabolic0.52}\\[2mm]
\raisebox{1.7cm}{(\(P_p\))} 
 \includegraphics[scale=.75
]{parabolic0.56}\hfill
  \includegraphics[scale=.75
]{parabolic0.57}\hfill
  \includegraphics[scale=.75
]{parabolic0.58}\\[2mm]
\raisebox{1cm}{(\(P_e\))}
\includegraphics[scale=.75
]{parabolic0.53}\hfill
  \includegraphics[scale=.75
]{parabolic0.54}\hfill
  \includegraphics[scale=.75
]{parabolic0.55}\\[2mm]
\raisebox{1.7cm}{(\(H\))}
  \includegraphics[scale=.75
]{parabolic0.60}\hfill
  \includegraphics[scale=.75
]{parabolic0.61}\hfill
  \includegraphics[scale=.75
]{parabolic0.62}\hfill
\caption[The elliptic, parabolic and hyperbolic unit disks]{The
  unit disks and orbits of subgroups \(A\),
  \(N\) and \(K\):\\
  (\(E\)): The elliptic unit disk;\\
  (\(P_e\)), (\(P_p\)), (\(P_h\)): The elliptic, parabolic and
  hyeprbolic flavour of the parabolic unit disk (the pure
  \emph{parabolic type} (\(P_p\)) transform is very similar with Figs.~\ref{fig:a-n-action}
  and~\ref{fig:k-subgroup}(\(K_p\))). \\
  \ \qquad(\(H\)): The hyperbolic unit disk. 
}
  \label{fig:unit-disks}
\end{figure}

\begin{lemma}
  \begin{enumerate}
  \item In the \emph{elliptic case} the ``real axis'' \(U\) is transformed
    to the unit circle and the upper half-plane---to the unit disk: 
    \begin{eqnarray}
      \label{eq:unit-circle-ell}
      \{(u,v) \such v = 0\} &\to& \{ (u',v') \such l^2_{c_e}(u'e_0+v'e_1)= u'^2+v'^2=1\} \\
      \label{eq:unit-disk-ell}
      \{(u,v) \such v > 0\} &\to& \{(u',v') \such
      l^2_{c_e}(u'e_0+v'e_1)= u'^2+v'^2<1\},
    \end{eqnarray}
    where the length from centre \(l^2_{c_e}\) is given
    by~\eqref{eq:k-center-point} for \(\sigma=\bs=-1\).

    On both sets \(\SL\) acts transitively and the unit circle is generated,
    for example, by the point \((0, 1)\) and the unit disk is generated by
    \((0,0)\). 
  \item In the \emph{hyperbolic case} the ``real axis'' \(U\) is transformed
    to the hyperbolic unit circle: 
    \begin{equation}
      \label{eq:unit-circle-hyp}
      \{(u,v) \such v = 0\} \to \{ (u',v') \such l^2_{c_h}(u',v')= u'^2-v'^2=-1\}, 
    \end{equation}
    where the length from centre \(l^2_{c_h}\) is given
    by~\eqref{eq:k-center-point} for \(\sigma=\bs=1\).

    On the  hyperbolic unit circle \(\SL\) acts transitively and it is
    generated, for example, by point \((0,1)\). 

    \(\SL\) acts also \emph{transitively on the whole complement}
    \begin{displaymath}
      \{(u',v') \such l^2_{c_h}(u'e_0+v'e_1)\neq -1\}
    \end{displaymath}
    to the unit circle, i.e. on its ``inner'' and ``outer'' parts
    together.
  \end{enumerate}
\end{lemma} The last feature of the hyperbolic Cayley transform can be
treated in a way described in the end of
subsection~\ref{sec:invar-upper-half}, see also
Fig.~\ref{fig:hyp-upper-half-plane}(b).  With such an arrangement
the hyperbolic Cayley transform maps the ``upper'' half-plane from
Fig.~\ref{fig:hyp-upper-half-plane}(a) onto the ``inner'' part of
the unit disk from Fig.~\ref{fig:hyp-upper-half-plane}(b) .

\begin{figure}[htbp]
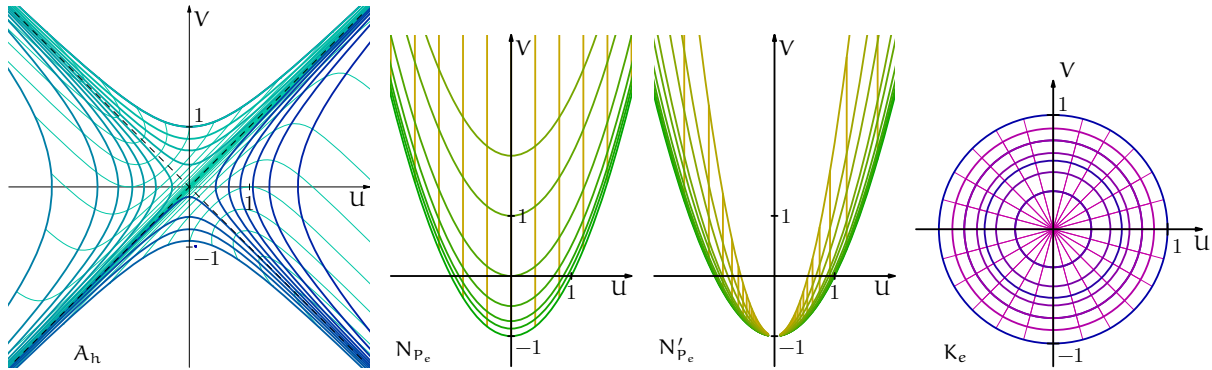

  \centering
\includegraphics[scale=.8
]{parabolic0.63}\hfill
\includegraphics[scale=.8
]{parabolic0.73}\hfill
\includegraphics[scale=.8
]{parabolic0.72}\hfill
\includegraphics[scale=.8
]{parabolic0.42}
  \caption[Cayley transform of the fix subgroups]{Concentric/confocal
    orbits of one parametric subgroups, cf. Fig.~\ref{fig:fix-sbroups}}
  \label{fig:concentric-equidist}
\end{figure}

One may wish that the hyperbolic Cayley transform diagonalises the
action of subgroup \(A\), or some conjugated, in a fashion similar
to the elliptic case~\eqref{eq:k-diagonalisation} for \(K\).
Geometrically it will correspond to hyperbolic rotations of hyperbolic
unit disk around the origin. Since the origin is the image of the
point \(e_1\) in the upper half-plane under the Cayley transform, we
will use the fix subgroup \(A'_h\)~\eqref{eq:hyp-fix-subgroup}
conjugated to \(A\) by \(\begin{pmatrix} 1 & e_0 \\ e_0 & 1
\end{pmatrix}\in \SL\).
Under the Cayley map~\eqref{eq:cayley-hyp} the subgroup \(A'_h\)
became, cf.~\cite{Kisil97c}*{(3.6--3.7)}: 
\begin{eqnarray*}
   \frac{1}{2}
   \begin{pmatrix}
     1&e_1\\
     -e_1&1
   \end{pmatrix}
  \begin{pmatrix}
    \cosh t & -e_0 \sinh t \\ 
    e_0\sinh t & \cosh t
  \end{pmatrix}
  \begin{pmatrix}
    1&-e_1\\
    e_1&1
  \end{pmatrix}
  =
  \begin{pmatrix}
    \exp(e_1e_0 t) &0\\
    0 & \exp(e_1e_0 t)
  \end{pmatrix},
\end{eqnarray*}
where \(\exp(e_1e_0t) = cosh(t)+e_1 e_0 sinh(t)\).  This obviously
corresponds to hyperbolic rotations of \(\Space{R}{h}\).  Orbits of
the fix subgroups \(A'_h\), \(N'_p\) and \(K'_e\) from
Lem.~\ref{le:fix-subgroups} under the Cayley transform are shown on
Fig.~\ref{fig:concentric-equidist}, which should be compared with
Fig.~\ref{fig:fix-sbroups}. However the parabolic Cayley transform
requires a separate discussion.

\subsection{Parabolic Cayley Transforms }
\label{sec:parab-cayl-transf}

This case benefits from a bigger variety of
choices.  The first natural attempt to define a Cayley
transform can be taken from the same formula~\eqref{eq:cayley-points}
with the parabolic value \(\sigma=0\). The corresponding
transformation defined by the matrix \(
\begin{pmatrix}
  1 & -e_1 \\ 0 & 1
\end{pmatrix}\) and defines the shift one unit down.

However within the framework of this paper a more general version of
parabolic Cayley transform is possible. It is given by
the matrix 
\begin{equation}
  \label{eq:parab-cayley-matr}
   C_{\bs}=\begin{pmatrix}
     1 & -e_1 \\ \bs e_1 & 1
   \end{pmatrix}, \qquad \text{ where } \bs=-1, 0, 1 \text{ and } \det
   C_{\bs}= 1 \text{ for all } \bs.
 \end{equation}
Here \(\bs=-1\) corresponds to the parabolic
Cayley transform \(P_e\) with the elliptic flavour, \(\bs=1\) ---
to the parabolic Cayley transform \(P_h\) with the hyperbolic flavour,
cf.~\cite{Kisil04b}*{\S~2.6}.
Finally the parabolic-parabolic transform is given by
an upper-triangular matrix from the end of the previous paragraph.
 
Fig.~\ref{fig:unit-disks} presents these transforms in rows (\(P_e\)),
(\(P_p\)) and (\(P_h\)) correspondingly. The row (\(P_p\)) almost coincides
with Figs.~\ref{fig:a-n-action}(\(A_a\)),
\ref{fig:a-n-action}(\(N_a\)) and \ref{fig:k-subgroup}(\(K_p\)).
Consideration of Fig.~\ref{fig:unit-disks} by columns from top to
bottom gives an impressive mixture of many common properties (e.g. the
number of fixed point on the boundary for each subgroup) with several
gradual mutations.

The description of the parabolic ``unit disk'' admits several
different interpretations in terms lengths from Defn.~\ref{de:length}. 
\begin{lemma}
  \label{le:parabolic-disk}
  Parabolic Cayley transform \(P_{\bs}\) as defined by the matrix
  \(C_{\bs}\)~\eqref{eq:parab-cayley-matr} acts on the \(V\)-axis
  always as a shift one unit down.

  Its image can be described in term of various lengths as follows:
  \begin{enumerate}
  \item \label{it:parab-unit-disk-cen} \(P_{\bs}\) for \(\bs\neq0\)
      transforms the ``real axis'' \(U\) to the p-cycle with the
      p-length squared \(-\bs\) from its e-centre
      \((0,-\frac{\bs}{2})\), cf.~\eqref{eq:unit-circle-ell}:
    \begin{equation}
      \label{eq:unit-circle-par-cen}
      \{(u,v) \such v = 0\} \to \{ (u',v') \such  l_{c_{e}}^2((0,-\frac{\bs}{2}), (u',v'))\cdot(-\bs)=1\}, 
    \end{equation}
    where \(l_{c_{e}}^2((0,-\frac{\bs}{2}), (u',v'))= u'^2+\bs
    v'\), see~\eqref{eq:k-center-point}. 

    The image of upper halfplane is:
    \begin{equation}
      \label{eq:unit-disk-par-cen}
      \{(u,v) \such v > 0\} \to 
         \{(u',v') \such l_{c_{e}}^2((0, -\frac{\bs}{2}), (u',v'))\cdot(-\bs)< 1\}.
    \end{equation}

  \item \label{it:parab-unit-disk-foc} \(P_{\bs}\) with \(\bs\neq0\)
    transforms the ``real axis'' \(U\) to the p-cycle with p-length squared
    \(-\bs\)~\eqref{eq:focal-length} from its h-focus
    \((0,-1-\frac{\bs}{4})\),  and the upper half-plane---to the ``interior'' part of it,
    cf.~\eqref{eq:unit-circle-ell}:
    \begin{eqnarray}
      \label{eq:unit-circle-par-foc}
      \{(u,v) \such v = 0\} &\to& \{ (u',v') \such
      l_{f_{h}}^2((0,{\textstyle-1-\frac{\bs}{4}}), 
      (u',v'))\cdot( -\bs)=1\}, \\
      \label{eq:unit-disk-par-foc}
      \{(u,v) \such v > 0\} &\to& \{(u',v') \such
      l_{f_{h}}^2((0,{\textstyle-1-\frac{\bs}{4}}), (u',v'))\cdot( -\bs)<1\}. 
    \end{eqnarray}
  \item \label{it:parab-unit-disk-p-foc} \(P_{\bs}\)
    transforms the ``real axis'' \(U\) to the cycle with p-length
    \(-\bs\) from its p-focus
    \((0,-1)\),  and the upper half-plane---to the ``interior'' part of it,
    cf.~\eqref{eq:unit-circle-ell}:
    \begin{eqnarray}
      \label{eq:unit-circle-par-p-foc}
      \{(u,v) \such v = 0\} &\to& \{ (u',v') \such
      l_{f_{p}}^2((0,-1), (u',v')) \cdot( -\bs)=1\}, \\
      \label{eq:unit-disk-par-p-foc}
      \{(u,v) \such v > 0\} &\to& \{(u',v') \such
      l_{f_{p}}^2((0,-1), (u',v'))\cdot( -\bs)<1\},
    \end{eqnarray}
    where \(l_{f_{p}}^2((0,-1), (u',v'))=\frac{u'^2}{v'+1}\)~\eqref{eq:focal-length}.
  \end{enumerate}
\end{lemma}
\begin{remark}
  Note that the both elliptic~\eqref{eq:unit-circle-ell} and
  hyperbolic~\eqref{eq:unit-circle-hyp} unit circle descriptions can
  be written uniformly with parabolic
  descriptions~\eqref{eq:unit-circle-par-cen}--\eqref{eq:unit-circle-par-p-foc} as
  \begin{displaymath}
     \{ (u',v') \such l^2_{c_{\bs}}(u'e_0+v'e_1)\cdot(-\bs)=1\} .
  \end{displaymath}
\end{remark}

The above descriptions~\ref{le:parabolic-disk}.\ref{it:parab-unit-disk-cen}
and~\ref{le:parabolic-disk}.\ref{it:parab-unit-disk-p-foc} are attractive for reasons given in
the following two lemmas. Firstly, the \(K\)-orbits in the elliptic
case (Fig.~\ref{fig:concentric-equidist}(\(K_e\))) and the
\(A\)-orbits in the hyperbolic case
(Fig.~\ref{fig:concentric-equidist}(\(A_h\))) of Cayley transform are
\emph{concentric}.
\begin{lemma}
  \label{le:n-orbits-concentric}
  \(N\)-orbits in the parabolic cases
  (Fig.~\ref{fig:unit-disks}(\(N_{P_{e}}\), \(N_{P_p}\),
  \(N_{P_{h}}\))) are \emph{concentric parabolas} (or straight lines)
  in the sense of Defn.~\ref{de:center-first} with e-centres at
  \((0,\frac{1}{2})\), \((0,\infty)\), \((0,-\frac{1}{2})\)
  correspondingly.
\end{lemma}

Secondly, Calley images of the fix subgroups' orbits in elliptic and
hyperbolic spaces in Fig.~\ref{fig:concentric-equidist}(\(A_h\))
and~(\(K_e\)) are equidistant from the origin in the corresponding
metrics.
\begin{lemma}
  \label{le:np-orbits-p-confocal}
  The Cayley transform of orbits of the parabolic fix subgroup in
  Fig.~\ref{fig:concentric-equidist}(\(N'_{P_e}\)) are parabolas
  consisting of points on the same
  \(l_{f_p}\)-length~\eqref{eq:k-focus-point} from the point \((0,-1)\),
  cf.~\ref{le:parabolic-disk}.\ref{it:parab-unit-disk-p-foc}. 
\end{lemma}
Note that parabolic rotations of the parabolic unit disk are
incompatible with the algebraic structure provided by the algebra of
dual numbers. However we can introduce~\cites{Kisil07a,Kisil09c} a linear
algebra structure and vector multiplication which will rotationally
invariant under action of subgroups \(N\) and \(N'\).

\begin{remark}
  \label{re:par-more-cayley}
  We see that the varieties of possible Cayley transforms in the
  parabolic case is bigger than in the two other cases. It is
  interesting that this parabolic richness is a consequence of the
  parabolic degeneracy of the generator \(e_1^2=0\). Indeed for both
  the elliptic and the hyperbolic signs in \(e_1^2=\pm 1\) only one
  matrix~\eqref{eq:cayley-points} out of two possible \(
  \begin{pmatrix}
    1 & e_1 \\ \pm \sigma e_1 & 1
  \end{pmatrix}\)\vspace{3pt} has a non-zero determinant. And only for
  the degenerate parabolic value \(e_1^2=0\) both these matrices are
  non-singular!
\end{remark}

\subsection{Cayley Transforms of Cycles}
\label{sec:cayl-transf-cycl}

\begin{figure}[htbp]
\centering
  \includegraphics[scale=1.2]{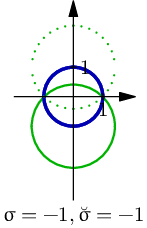}\hfill
  \includegraphics[scale=1.2]{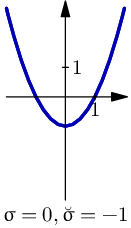}\hfill
  \includegraphics[scale=1.2]{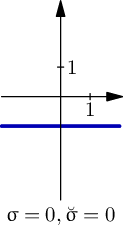}\hfill
  \includegraphics[scale=1.2]{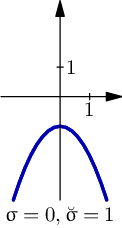}\hfill
  \includegraphics[scale=1.2]{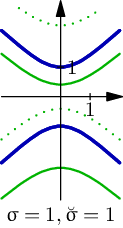}
  \caption[Cayley transforms in elliptic, parabolic and hyperbolic
  spaces]{Cayley transforms in elliptic (\(\sigma=-1\)), parabolic
    (\(\sigma=0\)) and hyperbolic (\(\sigma=1\)) spaces. On the each
    picture the reflection of the real line in the green cycles (drawn
    continuously or dotted) is the is the blue ``unit
    cycle''. Reflections in the solidly drawn cycles send the upper
    half-plane to the unit disk, reflection in the dashed cycle---to
    its complement.
    Three Cayley transforms in the parabolic space (\(\sigma=0\)) are
    themselves elliptic (\(\bs=-1\)), parabolic (\(\bs=0\)) and
    hyperbolic (\(\bs=1\)), giving a gradual transition between proper
    elliptic and hyperbolic cases.}
  \label{fig:unit-disk-all-cases}
\end{figure}

The next natural step within the SFSCc is to expand the Cayley
transform to the space of cycles. This is performed as follows:
\begin{lemma} Let \(\cycle{s}{a}\) be a cycle in \(\Space{R}{\sigma}\).
  \begin{itemize}
  \item[(e, h)] In the elliptic or hyperbolic cases the Cayley
    transform maps a cycle \(\cycle{s}{\bs}\) to the composition of its inversion 
    with the reflection
    \(\cycle[\hat]{s}{\bs} \cycle{s}{\bs}\cycle[\hat]{s}{\bs}\) in
    the cycle \(\cycle{s}{\bs}\), where \(\cycle[\hat]{s}{\bs}=
    \begin{pmatrix}
      \pm \se{1} && 1\\
      1 && \mp \se{1}
    \end{pmatrix}\)  with \(\bs=\pm1\) (see the first and last drawings on
    Fig.~\ref{fig:unit-disk-all-cases}). 
  \item[(p)] 
    \label{it:cayley-cycle-parab}
    In the parabolic case the Cayley transform maps a cycle
    \((k,l,n,m)\) to the cycle \((k-2\bs n, l, n ,m+2\bs n)\).
  \end{itemize}
\end{lemma}
The above extensions of the Cayley transform to the cycles space is
linear, however in the parabolic case it is not expressed as a
similarity of matrices (reflections in a cycle). This can be seen, for
example, from the fact that the parabolic Cayley transform does
not preserve the zero-radius cycles represented by matrices with zero
p-determinant. 
Since orbits of all subgroups in \(\SL\) as well as their Cayley
images are cycles in the corresponding metrics we may use
Lem.~\ref{it:cayley-cycle-parab}(p) to prove the following statements
(in addition to Lem.~\ref{le:n-orbits-concentric}):
\begin{corollary}
  \begin{enumerate}
  \item \(A\)-orbits in transforms \(P_e\) and \(P_h\) are segments of
    parabolas with the focal length \(\frac{1}{2}\) passing through
    \((0,-\frac{1}{2})\). Their vertices belong to two parabolas
    \(v=\frac{1}{2}(-x^2-1)\) and \(v=\frac{1}{2}(x^2-1)\)
    correspondingly, which are boundaries of parabolic circles in
    \(P_h\) and \(P_e\) (note the swap!).
  \item\label{it:parabolic-disk5}
    \(K\)-orbits in transform \(P_e\)  are parabolas
    with focal length less than \(\frac{1}{2}\) and in transform
    \(P_h\)---with inverse of focal length bigger than \(-2\).
  \end{enumerate}
\end{corollary}
Since the action of parabolic Cayley transform on cycles does not
preserve zero-radius cycles one shall better use infinitesimal-radius
cycles from \S~\ref{sec:zero-length-cycles} instead. First of all
images of infinitesimal cycles under parabolic Cayley transform are
infinitesimal cycles
again~\cite{Kisil05b}*{\S~\ref{G-sec:cayl-transf-infin}}, secondly
Lemma~\ref{le:infinitesimal-ortho}.\ref{it:infinitesimal-f-ortho} provides a useful expression of
concurrence with infinitesimal cycle focus through f-orthogonality.
Although f-orthogonality is not preserved by the Cayley
transform~\ref{it:cayley-cycle-parab}(p) for generic cycles it did for
the infinitesimal ones,
see~\cite{Kisil05b}*{\S~\ref{G-sec:cayl-transf-infin}}:
\begin{lemma}
  An infinitesimal cycle \(\cycle{a}{\bs}\)~\eqref{eq:inf-cycle} is
  f-orthogonal (in the sense of Lem.~\ref{le:infinitesimal-ortho}.\ref{it:infinitesimal-f-ortho}) to a
  cycle \(\cycle[\tilde]{a}{\bs}\) if and only if the Cayley
  transform~\ref{it:cayley-cycle-parab}(p) of \(\cycle{a}{\bs}\) is
  f-orthogonal to the Cayley transform of \(\cycle[\tilde]{a}{\bs}\).
\end{lemma}
We main observation of this paper is that the potential of the
Erlangen programme is still far from exhausting even for
two-dimensional geometry.

\section*{Acknowledgements}
\label{sec:acknowledgments}
This paper has some overlaps with the paper~\cite{Kisil04b} written in
collaboration with D.~Biswas.  However the present paper essentially
revises many concepts (e.g. lengths, orthogonality, the parabolic
Cayley transform) introduced in~\cite{Kisil04b}, thus it was important
to make it an independent reading to avoid confusion with some earlier
(and na\"{\i}ve!) guesses made in~\cite{Kisil04b}.

The author is grateful to Professors S.~Plaksa, S.~Blyumin and
N.A.~Gromov for useful discussions and comments.  Drs.~I.R.~Porteous,
D.L.~Selinger and J.~Selig carefully read the previous
paper~\cite{Kisil04b} and made numerous comments and remarks helping
to improve this paper as well. I am also grateful to D.~Biswas for
many comments on this paper.

The extensive graphics in this paper were produced with the help of
the \GiNaC~\cites{GiNaC,Kisil04c} computer algebra system. Since this
tool is of separate interest we explain its usage by examples from
this article in the separate paper~\cite{Kisil05b}. The
\NoWEB~\cite{NoWEB} wrapper for \CPP\ source code is included in the
\href{http://arXiv.org}{\texttt{arXiv.org}} files of the
papers~\cite{Kisil05b}.

\small
\bibliographystyle{plain}
\bibliography{arare,abbrevmr,akisil,aclifford,ageometry,analyse,aphysics}

\providecommand{\CPP}{\texttt{C++}} \providecommand{\NoWEB}{\texttt{noweb}}
  \providecommand{\MetaPost}{\texttt{Meta}\-\texttt{Post}}
  \providecommand{\GiNaC}{\textsf{GiNaC}}
  \providecommand{\pyGiNaC}{\textsf{pyGiNaC}}
  \providecommand{\Asymptote}{\texttt{Asymptote}}
  \providecommand{\noopsort}[1]{} \providecommand{\printfirst}[2]{#1}
  \providecommand{\singleletter}[1]{#1} \providecommand{\switchargs}[2]{#2#1}
  \providecommand{\irm}{\textup{I}} \providecommand{\iirm}{\textup{II}}
  \providecommand{\vrm}{\textup{V}} \providecommand{\cprime}{'}
  \providecommand{\eprint}[2]{\texttt{#2}}
  \providecommand{\myeprint}[2]{\texttt{#2}}
  \providecommand{\arXiv}[1]{\myeprint{http://arXiv.org/abs/#1}{arXiv:#1}}
  \providecommand{\doi}[1]{\href{http://dx.doi.org/#1}{doi: #1}}
\begin{bibdiv}
\begin{biblist}

\bib{Arveson76}{book}{
      author={Arveson, William},
       title={An invitation to {C}*-algebras},
   publisher={Springer-Verlag},
        date={{\noopsort{}}1976},
}

\bib{BairdWood09a}{article}{
      author={Baird, Paul},
      author={Wood, John~C.},
       title={Harmonic morphisms from {M}inkowski space and hyperbolic
  numbers},
        date={2009},
        ISSN={1220-3874},
     journal={Bull. Math. Soc. Sci. Math. Roumanie (N.S.)},
      volume={52(100)},
      number={3},
       pages={195\ndash 209},
      review={\MR{2554640 (2010h:53093)}},
}

\bib{BarkerHowe07}{book}{
      author={Barker, William},
      author={Howe, Roger},
       title={Continuous symmetry. {From} {Euclid} to {Klein}},
   publisher={American Mathematical Society},
     address={Providence, RI},
        date={2007},
        ISBN={978-0-8218-3900-3; 0-8218-3900-4},
      review={\MR{MR2362745 (2008k:51001)}},
}

\bib{GiNaC}{misc}{
      author={Bauer, Christian},
      author={Frink, Alexander},
      author={Kreckel, Richard},
      author={Vollinga, Jens},
       title={\textsf{GiNaC} is {N}ot a {CAS}},
        date={2001--},
        note={URL: \url{http://www.ginac.de/}},
}

\bib{Beardon95}{book}{
      author={Beardon, Alan~F.},
       title={The geometry of discrete groups},
      series={Graduate Texts in Mathematics},
   publisher={Springer-Verlag},
     address={New York},
        date={1995},
      volume={91},
        ISBN={0-387-90788-2},
        note={Corrected reprint of the 1983 original},
      review={\MR{MR1393195 (97d:22011)}},
}

\bib{Beardon05a}{book}{
      author={Beardon, Alan~F.},
       title={Algebra and geometry},
   publisher={Cambridge University Press},
     address={Cambridge},
        date={2005},
        ISBN={0-521-89049-7},
      review={\MR{MR2153234 (2006a:00001)}},
}

\bib{BekkaraFrancesZeghib06}{article}{
      author={Bekkara, Esmaa},
      author={Frances, Charles},
      author={Zeghib, Abdelghani},
       title={On lightlike geometry: isometric actions, and rigidity aspects},
        date={2006},
        ISSN={1631-073X},
     journal={C. R. Math. Acad. Sci. Paris},
      volume={343},
      number={5},
       pages={317\ndash 321},
      review={\MR{MR2253050 (2007e:53045)}},
}

\bib{Benz07a}{book}{
      author={Benz, Walter},
       title={Classical geometries in modern contexts. {Geometry} of real inner
  product spaces},
     edition={Second edition},
   publisher={Birkh\"auser Verlag},
     address={Basel},
        date={2007},
        ISBN={978-3-7643-8540-8},
      review={\MR{MR2370626 (2008i:51001)}},
}

\bib{BergerII}{book}{
      author={Berger, Marcel},
       title={Geometry. {II}},
      series={Universitext},
   publisher={Springer-Verlag},
     address={Berlin},
        date={1987},
        ISBN={3-540-17015-4},
        note={Translated from the French by M. Cole and S. Levy},
      review={\MR{882916 (88a:51001b)}},
}

\bib{BocCatoniCannataNichZamp07}{book}{
      author={Boccaletti, Dino},
      author={Catoni, Francesco},
      author={Cannata, Roberto},
      author={Catoni, Vincenzo},
      author={Nichelatti, Enrico},
      author={Zampetti, Paolo},
       title={The mathematics of {Minkowski} space-time and an introduction to
  commutative hypercomplex numbers},
   publisher={Birkh\"auser Verlag},
     address={Basel},
        date={2008},
}

\bib{CatoniCannataNichelatti04}{article}{
      author={Catoni, Francesco},
      author={Cannata, Roberto},
      author={Nichelatti, Enrico},
       title={The parabolic analytic functions and the derivative of real
  functions},
        date={2004},
     journal={Adv. Appl. Clifford Algebras},
      volume={14},
      number={2},
       pages={185\ndash 190},
}

\bib{CatoniCannataZampetti05}{article}{
      author={Catoni, Francesco},
      author={Cannata, Vincenzo, Roberto~Catoni},
      author={Zampetti, Paolo},
       title={{$N$}-dimensional geometries generated by hypercomplex numbers},
        date={2005},
     journal={Adv. Appl. Clifford Algebras},
      volume={15},
      number={1},
}

\bib{CerejeirasKahlerSommen05a}{article}{
      author={Cerejeiras, P.},
      author={K{\"a}hler, U.},
      author={Sommen, F.},
       title={Parabolic {D}irac operators and the {N}avier-{S}tokes equations
  over time-varying domains},
        date={2005},
        ISSN={0170-4214},
     journal={Math. Methods Appl. Sci.},
      volume={28},
      number={14},
       pages={1715\ndash 1724},
      review={\MR{MR2167561}},
}

\bib{Chern96a}{article}{
      author={Chern, Shiing-Shen},
       title={Finsler geometry is just {R}iemannian geometry without the
  quadratic restriction},
        date={1996},
        ISSN={0002-9920},
     journal={Notices Amer. Math. Soc.},
      volume={43},
      number={9},
       pages={959\ndash 963},
      review={\MR{MR1400859 (97e:53129)}},
}

\bib{Cnops94a}{thesis}{
      author={Cnops, Jan},
       title={{Hurwitz} pairs and applications of {M\"obius} transformations},
        type={{Habilitation} Dissertation},
     address={Universiteit Gent, Faculteit van de Wetenschappen},
        date={1994},
        note={See also~\cite{Cnops02a}},
}

\bib{Cnops02a}{book}{
      author={Cnops, Jan},
       title={An introduction to {D}irac operators on manifolds},
      series={Progress in Mathematical Physics},
   publisher={Birkh\"auser Boston Inc.},
     address={Boston, MA},
        date={2002},
      volume={24},
        ISBN={0-8176-4298-6},
      review={\MR{1 917 405}},
}

\bib{CoxeterGreitzer}{book}{
      author={Coxeter, H.S.M.},
      author={Greitzer, S.L.},
       title={{Geometry revisited.}},
   publisher={{Random House}},
     address={New York},
        date={1967},
        note={\Zbl{0166.16402}},
}

\bib{Devis77}{book}{
      author={Davis, Martin},
       title={Applied nonstandard analysis},
   publisher={Wiley-Interscience [John Wiley \& Sons]},
     address={New York},
        date={1977},
        ISBN={0-471-19897-8},
      review={\MR{MR0505473 (58 \#21590)}},
}

\bib{EelbodeSommen04a}{article}{
      author={Eelbode, D.},
      author={Sommen, F.},
       title={Taylor series on the hyperbolic unit ball},
        date={2004},
        ISSN={1370-1444},
     journal={Bull. Belg. Math. Soc. Simon Stevin},
      volume={11},
      number={5},
       pages={719\ndash 737},
      review={\MR{MR2130635 (2005k:30092)}},
}

\bib{EelbodeSommen05a}{article}{
      author={Eelbode, D.},
      author={Sommen, F.},
       title={The fundamental solution of the hyperbolic {D}irac operator on
  {$\mathbb{R}\sp {1,m}$}: a new approach},
        date={2005},
        ISSN={1370-1444},
     journal={Bull. Belg. Math. Soc. Simon Stevin},
      volume={12},
      number={1},
       pages={23\ndash 37},
      review={\MR{MR2134853 (2005k:30093)}},
}

\bib{FillmoreSpringer90a}{article}{
      author={Fillmore, Jay~P.},
      author={Springer, A.},
       title={M\"obius groups over general fields using {C}lifford algebras
  associated with spheres},
        date={1990},
        ISSN={0020-7748},
     journal={Internat. J. Theoret. Phys.},
      volume={29},
      number={3},
       pages={225\ndash 246},
         url={http://dx.doi.org/10.1007/BF00673627},
      review={\MR{1049005 (92a:22016)}},
}

\bib{FjelstadGal01}{article}{
      author={Fjelstad, Paul},
      author={Gal, Sorin~G.},
       title={Two-dimensional geometries, topologies, trigonometries and
  physics generated by complex-type numbers},
        date={2001},
        ISSN={0188-7009},
     journal={Adv. Appl. Clifford Algebras},
      volume={11},
      number={1},
       pages={81\ndash 107 (2002)},
      review={\MR{MR1953513 (2003k:83009)}},
}

\bib{Garasko09a}{book}{
      author={Garas{\cprime}ko, G.I.},
       title={{\cyr Nachala finslerovoi0 geometrii dlya fizikov.} [{E}lements
  of {F}insler geometry for physicists]},
    language={Russian},
   publisher={TETRU},
     address={Moscow},
        date={2009},
        note={268 pp. URL:
  \url{http://hypercomplex.xpsweb.com/articles/487/ru/pdf/00-gbook.pdf}},
}

\bib{Gromov90a}{book}{
      author={Gromov, N.~A.},
       title={{\cyr Kontraktsii i analiticheskie prodolzheniya klassicheskikh
  grupp. {E}dinyi podkhod}. ({R}ussian) [{C}ontractions and analytic extensions
  of classical groups. {U}nified approach]},
   publisher={Akad. Nauk SSSR Ural. Otdel. Komi Nauchn. Tsentr},
     address={Syktyvkar},
        date={1990},
      review={\MR{MR1092760 (91m:81078)}},
}

\bib{GromovKuratov05a}{article}{
      author={Gromov, N.~A.},
      author={Kuratov, V.~V.},
       title={Noncommutative space-time models},
        date={2005},
        ISSN={0011-4626},
     journal={Czechoslovak J. Phys.},
      volume={55},
      number={11},
       pages={1421\ndash 1426},
         url={http://dx.doi.org/10.1007/s10582-006-0020-y},
      review={\MR{MR2223830 (2006k:81197)}},
}

\bib{HerranzSantander02b}{article}{
      author={Herranz, Francisco~J.},
      author={Santander, Mariano},
       title={Conformal compactification of spacetimes},
        date={2002},
        ISSN={0305-4470},
     journal={J. Phys. A},
      volume={35},
      number={31},
       pages={6619\ndash 6629},
        note={\arXiv{math-ph/0110019}},
      review={\MR{MR1928852 (2004b:53123)}},
}

\bib{HerranzSantander02a}{article}{
      author={Herranz, Francisco~J.},
      author={Santander, Mariano},
       title={Conformal symmetries of spacetimes},
        date={2002},
        ISSN={0305-4470},
     journal={J. Phys. A},
      volume={35},
      number={31},
       pages={6601\ndash 6618},
        note={\arXiv{math-ph/0110019}},
      review={\MR{MR1928851 (2004a:53091)}},
}

\bib{HoweTan92}{book}{
      author={Howe, Roger},
      author={Tan, Eng-Chye},
       title={Nonabelian harmonic analysis. {Applications of
  ${{\rm{S}}L}(2,{{\bf{R}}})$}},
   publisher={Springer-Verlag},
     address={New York},
        date={1992},
        ISBN={0-387-97768-6},
      review={\MR{1151617 (93f:22009)}},
}

\bib{Khrennikov05a}{article}{
      author={Khrennikov, A.~Yu.},
       title={Hyperbolic quantum mechanics},
        date={2005},
        ISSN={0869-5652},
     journal={Dokl. Akad. Nauk},
      volume={402},
      number={2},
       pages={170\ndash 172},
      review={\MR{MR2162434 (2006d:81118)}},
}

\bib{KhrennikovSegre07a}{incollection}{
      author={Khrennikov, Andrei},
      author={Segre, Gavriel},
       title={Hyperbolic quantization},
        date={2007},
   booktitle={Quantum probability and infinite dimensional analysis},
      editor={Accardi, L.},
      editor={Freudenberg, W.},
      editor={Sch\"urman, M.},
      series={QP--PQ: Quantum Probab. White Noise Anal.},
      volume={20},
   publisher={World Scientific Publishing, Hackensack, NJ},
       pages={282\ndash 287},
      review={\MR{MR2359402}},
}

\bib{Kirillov76}{book}{
      author={Kirillov, A.~A.},
       title={Elements of the theory of representations},
   publisher={Springer-Verlag},
     address={Berlin},
        date={1976},
        note={Translated from the Russian by Edwin Hewitt, Grundlehren der
  Mathematischen Wissenschaften, Band 220},
      review={\MR{54 \#447}},
}

\bib{Kirillov06}{book}{
      author={Kirillov, A.~A.},
       title={A tale on two fractals},
        date={2010},
        note={URL:
  \url{http://www.math.upenn.edu/~kirillov/MATH480-F07/tf.pdf}, in
  publication},
}

\bib{Kisil08a}{article}{
      author={Kisil, Anastasia~V.},
       title={Isometric action of {${\rm SL}_2(\mathbb{R})$} on homogeneous
  spaces},
        date={2010},
     journal={Adv. App. Clifford Algebras},
      volume={20},
      number={2},
       pages={299\ndash 312},
        note={\arXiv{0810.0368}. \MR{2012b:32019}},
}

\bib{Kisil95d}{article}{
      author={Kisil, Vladimir~V.},
       title={Construction of integral representations for spaces of analytic
  functions},
        date={1996},
        ISSN={0869-5652},
     journal={Dokl. Akad. Nauk},
      volume={350},
      number={4},
       pages={446\ndash 448},
        note={(Russian) \MR{98d:46027}},
}

\bib{Kisil95i}{article}{
      author={Kisil, Vladimir~V.},
       title={M\"obius transformations and monogenic functional calculus},
        date={1996},
        ISSN={1079-6762},
     journal={Electron. Res. Announc. Amer. Math. Soc.},
      volume={2},
      number={1},
       pages={26\ndash 33},
  note={\href{http://www.ams.org/era/1996-02-01/S1079-6762-96-00004-2/}{On-line}},
      review={\MR{MR1405966 (98a:47018)}},
}

\bib{Kisil96d}{incollection}{
      author={Kisil, Vladimir~V.},
       title={Towards to analysis in {${\mathbb{R}}^{pq}$}},
        date={{\noopsort{}}1996},
   booktitle={Proceedings of symposium analytical and numerical methods in
  quaternionic and clifford analysis, june 5--7, 1996, seiffen, germany},
      editor={Spr\"o{\ss}ig, Wolfgang},
      editor={G\"urlebeck, Klaus},
   publisher={TU Bergakademie Freiberg},
     address={Freiberg},
       pages={95\ndash 100},
        note={\Zbl{882.30030}},
}

\bib{Kisil96c}{incollection}{
      author={Kisil, Vladimir~V.},
       title={How many essentially different function theories exist?},
        date={1998},
   booktitle={Clifford algebras and their application in mathematical physics
  ({A}achen, 1996)},
      series={Fund. Theories Phys.},
      volume={94},
   publisher={Kluwer Academic Publishers},
     address={Dordrecht},
       pages={175\ndash 184},
      review={\MR{MR1627084 (99g:30057)}},
}

\bib{Kisil97c}{article}{
      author={Kisil, Vladimir~V.},
       title={Analysis in {$\mathbf{R}\sp {1,1}$} or the principal function
  theory},
        date={1999},
        ISSN={0278-1077},
     journal={Complex Variables Theory Appl.},
      volume={40},
      number={2},
       pages={93\ndash 118},
        note={\arXiv{funct-an/9712003}},
      review={\MR{MR1744876 (2000k:30078)}},
}

\bib{Kisil94e}{article}{
      author={Kisil, Vladimir~V.},
       title={Relative convolutions. {I}. {P}roperties and applications},
        date={1999},
        ISSN={0001-8708},
     journal={Adv. Math.},
      volume={147},
      number={1},
       pages={35\ndash 73},
        note={\arXiv{funct-an/9410001},
  \href{http://www.idealibrary.com/links/doi/10.1006/aima.1999.1833}{On-line}.
  \Zbl{933.43004}},
      review={\MR{MR1725814 (2001h:22012)}},
}

\bib{Kisil97a}{incollection}{
      author={Kisil, Vladimir~V.},
       title={Two approaches to non-commutative geometry},
        date={1999},
   booktitle={Complex methods for partial differential equations ({A}nkara,
  1998)},
      series={Int. Soc. Anal. Appl. Comput.},
      volume={6},
   publisher={Kluwer Acad. Publ.},
     address={Dordrecht},
       pages={215\ndash 244},
        note={\arXiv{funct-an/9703001}},
      review={\MR{MR1744440 (2001a:01002)}},
}

\bib{Kisil98a}{article}{
      author={Kisil, Vladimir~V.},
       title={Wavelets in {B}anach spaces},
        date={1999},
        ISSN={0167-8019},
     journal={Acta Appl. Math.},
      volume={59},
      number={1},
       pages={79\ndash 109},
        note={\arXiv{math/9807141},
  \href{http://dx.doi.org/10.1023/A:1006394832290}{On-line}},
      review={\MR{MR1740458 (2001c:43013)}},
}

\bib{Kisil01a}{misc}{
      author={Kisil, Vladimir~V.},
       title={Spaces of analytical functions and wavelets---{Lecture} notes},
        date={2000--2002},
        note={92 pp. \arXiv{math.CV/0204018}},
}

\bib{Kisil02c}{incollection}{
      author={Kisil, Vladimir~V.},
       title={Meeting {D}escartes and {K}lein somewhere in a noncommutative
  space},
        date={2002},
   booktitle={Highlights of mathematical physics ({L}ondon, 2000)},
      editor={Fokas, A.},
      editor={Halliwell, J.},
      editor={Kibble, T.},
      editor={Zegarlinski, B.},
   publisher={Amer. Math. Soc.},
     address={Providence, RI},
       pages={165\ndash 189},
        note={\arXiv{math-ph/0112059}},
      review={\MR{MR2001578 (2005b:43015)}},
}

\bib{Kisil02a}{inproceedings}{
      author={Kisil, Vladimir~V.},
       title={Spectrum as the support of functional calculus},
        date={2004},
   booktitle={Functional analysis and its applications},
      series={North-Holland Math. Stud.},
      volume={197},
   publisher={Elsevier},
     address={Amsterdam},
       pages={133\ndash 141},
        note={\arXiv{math.FA/0208249}},
      review={\MR{MR2098877}},
}

\bib{Kisil04c}{article}{
      author={Kisil, Vladimir~V.},
       title={An example of {C}lifford algebras calculations with
  \texttt{GiNaC}},
        date={2005},
        ISSN={0188-7009},
     journal={Adv. Appl. Clifford Algebr.},
      volume={15},
      number={2},
       pages={239\ndash 269},
        note={\arXiv{cs.MS/0410044},
  \href{http://www.clifford-algebras.org/v15/v152/kisil152.pdf}{On-line}},
      review={\MR{MR2241255 (2007a:15042)}},
}

\bib{Kisil06a}{article}{
      author={Kisil, Vladimir~V.},
       title={Erlangen program at large--0: Starting with the group {${\rm
  SL}\sb 2({\bf R})$}},
        date={2007},
        ISSN={0002-9920},
     journal={Notices Amer. Math. Soc.},
      volume={54},
      number={11},
       pages={1458\ndash 1465},
        note={\arXiv{math/0607387},
  \href{http://www.ams.org/notices/200711/tx071101458p.pdf}{On-line}},
      review={\MR{MR2361159}},
}

\bib{Kisil05b}{article}{
      author={Kisil, Vladimir~V.},
       title={{S}chwerdtfeger--{F}illmore-{S}pringer-{C}nops construction
  implemented in \texttt{GiNaC}},
        date={2007},
        ISSN={0188-7009},
     journal={Adv. Appl. Clifford Algebr.},
      volume={17},
      number={1},
       pages={59\ndash 70},
        note={\href{http://dx.doi.org/10.1007/s00006-006-0017-4}{On-line}.
  Updated full text, source files, and live ISO image: \arXiv{cs.MS/0512073}},
      review={\MR{MR2303056}},
}

\bib{Kisil07a}{inproceedings}{
      author={Kisil, Vladimir~V.},
       title={Erlangen program at large---2: {Inventing} a wheel. {The}
  parabolic one},
        date={2010},
   booktitle={Trans. {I}nst. {M}ath. of the {NAS} of {U}kraine},
      series={Trans. {I}nst. {M}ath. of the {NAS} of {U}kraine},
      volume={7},
       pages={89\ndash 98},
        note={\arXiv{0707.4024}},
}

\bib{Kisil05a}{article}{
      author={Kisil, Vladimir~V.},
       title={Erlangen program at large--1: Geometry of invariants},
        date={2010},
     journal={SIGMA, Symmetry Integrability Geom. Methods Appl.},
      volume={6},
      number={076},
       pages={45},
        note={\arXiv{math.CV/0512416}. \MR{2011i:30044}},
}

\bib{Kisil10b}{article}{
      author={Kisil, Vladimir~V.},
       title={{E}rlangen {P}rogramme at {L}arge: {A} brief outline},
        date={2010},
        note={\arXiv{1006.2115}},
}

\bib{Kisil09c}{article}{
      author={Kisil, Vladimir~V.},
       title={{E}rlangen program at large---2 1/2: {I}nduced representations
  and hypercomplex numbers},
        date={2011},
     journal={{\cyr Izvestiya Komi nauchnogo centra UrO RAN} [Izvestiya Komi
  nauchnogo centra UrO RAN]},
      volume={1},
      number={5},
       pages={4\ndash 10},
        note={\arXiv{0909.4464}},
}

\bib{Kisil12a}{book}{
      author={Kisil, Vladimir~V.},
       title={Geometry of {M}\"obius transformations: {E}lliptic, parabolic and
  hyperbolic actions of {$\mathrm{SL}_2(\mathbf{R})$}},
   publisher={Imperial College Press},
     address={London},
        date={2012},
        note={Includes a live DVD},
}

\bib{Kisil10a}{article}{
      author={Kisil, Vladimir~V.},
       title={Hypercomplex representations of the {H}eisenberg group and
  mechanics},
        date={2012},
        ISSN={0020-7748},
     journal={Internat. J. Theoret. Phys.},
      volume={51},
      number={3},
       pages={964\ndash 984},
         url={http://dx.doi.org/10.1007/s10773-011-0970-0},
        note={\arXiv{1005.5057}},
      review={\MR{2892069}},
}

\bib{Kisil04b}{inproceedings}{
      author={Kisil, Vladimir~V.},
      author={Biswas, Debapriya},
       title={Elliptic, parabolic and hyperbolic analytic function theory--0:
  Geometry of domains},
        date={2004},
   booktitle={Complex analysis and free boundary flows},
      series={Trans. {I}nst. {M}ath. of the {NAS} of {U}kraine},
      volume={1},
       pages={100\ndash 118},
        note={\arXiv{math.CV/0410399}},
}

\bib{Konovenko09a}{article}{
      author={Konovenko, Nadiia},
       title={Projective structures and algebras of their differential
  invariants},
        date={2010},
     journal={Acta Applicandae Mathematicae},
      volume={109},
      number={1},
       pages={87\ndash 99},
}

\bib{KonovenkoLychagin08a}{article}{
      author={Konovenko, N.G.},
      author={Lychagin, V.V.},
       title={{Differential invariants of nonstandard projective structures.}},
    language={Russian. English summary},
        date={2008},
     journal={Dopov. Nats. Akad. Nauk Ukr., Mat. Pryr. Tekh. Nauky},
      volume={2008},
      number={11},
       pages={10\ndash 13},
        note={\Zbl{1164.53313}},
      review={\MR{MR2529116}},
}

\bib{KuruczWolterZakharyaschev05}{inproceedings}{
      author={Kurucz, A.},
      author={Wolter, F.},
      author={Zakharyaschev, Michael},
       title={Modal logics for metric spaces: Open problems},
        date={2005},
   booktitle={We will show them: Essays in honour of {Dov Gabbay}},
      editor={Artemov, S.},
      editor={Barringer, H.},
      editor={d'Avila Garcez, A.~S.},
      editor={Lamb, L.~C.},
      editor={Woods, J.},
      volume={2},
   publisher={College Publications},
       pages={\href{http://www.dcs.bbk.ac.uk/~michael/dov.pdf}{193\ndash 208}},
}

\bib{Lang85}{book}{
      author={Lang, Serge},
       title={{${\rm SL}\sb 2({\bf R})$}},
      series={Graduate Texts in Mathematics},
   publisher={Springer-Verlag},
     address={New York},
        date={1985},
      volume={105},
        ISBN={0-387-96198-4},
        note={Reprint of the 1975 edition},
      review={\MR{803508 (86j:22018)}},
}

\bib{LavrentShabat77}{book}{
      author={Lavrent{\cprime}ev, M.~A.},
      author={Shabat, B.~V.},
       title={{\cyr Problemy gidrodinamiki i ikh matematicheskie modeli}.
  [{P}roblems of hydrodynamics and their mathematical models]},
    language={Russian},
     edition={Second},
   publisher={Izdat. ``Nauka'', Moscow},
        date={1977},
      review={\MR{56 \#17392}},
}

\bib{McRae07}{article}{
      author={McRae, Alan~S.},
       title={Clifford algebras and possible kinematics},
        date={2007},
        ISSN={1815-0659},
     journal={SIGMA Symmetry Integrability Geom. Methods Appl.},
      volume={3},
       pages={Paper 079, 29 pp. (electronic)},
         url={http://dx.doi.org/10.3842/SIGMA.2007.079},
      review={\MR{MR2322806 (2008d:15074)}},
}

\bib{Mirman01a}{book}{
      author={Mirman, R.},
       title={Quantum field theory, conformal group theory, conformal field
  theory. {M}athematical and conceptual foundations, physical and geometrical
  applications},
   publisher={Nova Science Publishers Inc.},
     address={Huntington, NY},
        date={2001},
        ISBN={1-56072-992-9},
      review={\MR{MR2145706 (2006j:81001)}},
}

\bib{MotterRosa98}{article}{
      author={Motter, A.~E.},
      author={Rosa, M. A.~F.},
       title={Hyperbolic calculus},
        date={1998},
        ISSN={0188-7009},
     journal={Adv. Appl. Clifford Algebras},
      volume={8},
      number={1},
       pages={109\ndash 128},
      review={\MR{MR1648837 (99m:30099)}},
}

\bib{Olver99}{book}{
      author={Olver, Peter~J.},
       title={Classical invariant theory},
      series={London Mathematical Society Student Texts},
   publisher={Cambridge University Press},
     address={Cambridge},
        date={1999},
      volume={44},
        ISBN={0-521-55821-2},
      review={\MR{MR1694364 (2001g:13009)}},
}

\bib{Pavlov06a}{article}{
      author={Pavlov, D.~G.},
       title={Symmetries and geometric invariants},
    language={Russian},
        date={2006},
     journal={{\cyr Giperkompleksny{}e chisla v geometrii i fizike}
  (Hypercomplex Numbers in Geometry and Physics)},
      volume={3},
      number={2},
       pages={21\ndash 32},
        note={URL:
  \url{http://hypercomplex.xpsweb.com/page.php?lang=ru&id=316}},
}

\bib{Pimenov65a}{article}{
      author={Pimenov, R.I.},
       title={Unified axiomatics of spaces with maximal movement group},
    language={Russian},
        date={1965},
     journal={Litov. Mat. Sb.},
      volume={5},
       pages={457\ndash 486},
        note={\Zbl{0139.37806}},
}

\bib{Porteous95}{book}{
      author={Porteous, Ian~R.},
       title={Clifford algebras and the classical groups},
      series={Cambridge Studies in Advanced Mathematics},
   publisher={Cambridge University Press},
     address={Cambridge},
        date={1995},
      volume={50},
        ISBN={0-521-55177-3},
      review={\MR{MR1369094 (97c:15046)}},
}

\bib{NoWEB}{misc}{
      author={Ramsey, Norman},
       title={\href{http://www.eecs.harvard.edu/~nr/noweb/}{Noweb} --- a
  simple, extensible tool for literate programming},
        note={URL: \url{http://www.eecs.harvard.edu/~nr/noweb/}},
}

\bib{RozenfeldZamakhovski03}{book}{
      author={Rozenfel{\cprime}d, B.~A.},
      author={Zamakhovski{\u\i}, M.~P.},
       title={{\cyr Geometriya grupp {L}i. {S}immetricheskie, parabolicheskie i
  periodicheskie prostranstva.} ({R}ussian) [{G}eometry of {L}ie groups.
  {S}ymmetric, parabolic, and periodic spaces]},
   publisher={Moskovski\u\i\ Tsentr Nepreryvnogo Matematicheskogo Obrazovaniya,
  Moscow},
        date={2003},
        ISBN={5-94057-032-1},
      review={\MR{MR2070039 (2005f:53080)}},
}

\bib{Schwerdtfeger79a}{book}{
      author={Schwerdtfeger, Hans},
       title={Geometry of complex numbers},
   publisher={Dover Publications Inc.},
     address={New York},
        date={1979},
        ISBN={0-486-63830-8},
        note={Circle geometry, Moebius transformation, non-Euclidean geometry,
  A corrected reprinting of the 1962 edition, Dover Books on Advanced
  Mathematics},
      review={\MR{620163 (82g:51032)}},
}

\bib{Segal76}{book}{
      author={Segal, Irving~Ezra},
       title={Mathematical cosmology and extragalactic astronomy},
      series={Pure and Applied Mathematics},
   publisher={Academic Press [Harcourt Brace Jovanovich Publishers]},
     address={New York},
        date={1976},
      volume={68},
      review={\MR{58 \#14894}},
}

\bib{Sharpe97}{book}{
      author={Sharpe, R.~W.},
       title={Differential geometry. {C}artan's generalization of {K}lein's
  {E}rlangen program},
      series={Graduate Texts in Mathematics},
   publisher={Springer-Verlag},
     address={New York},
        date={1997},
      volume={166},
        ISBN={0-387-94732-9},
        note={With a foreword by S. S. Chern},
      review={\MR{MR1453120 (98m:53033)}},
}

\bib{MTaylor86}{book}{
      author={Taylor, Michael~E.},
       title={Noncommutative harmonic analysis},
      series={Mathematical Surveys and Monographs},
   publisher={American Mathematical Society},
     address={Providence, RI},
        date={1986},
      volume={22},
        ISBN={0-8218-1523-7},
      review={\MR{88a:22021}},
}

\bib{Uspenskii88}{book}{
      author={Uspenski{\u\i}, V.~A.},
       title={{\cyr {CH}to takoe nestandartnyi0 analiz?} [{W}hat is
  non-standard analysis?]},
    language={Russian},
   publisher={``Nauka''},
     address={Moscow},
        date={1987},
        note={With an appendix by V. G. Kanove\u\i},
      review={\MR{MR913941 (88m:26028)}},
}

\bib{VignauxDuranona35a}{article}{
      author={Vignaux, J.~C.},
      author={Dura\~nona~y Vedia, A.},
       title={{Sobre la teor{\'\i}a de las funciones de una variable compleja
  hiperb\'olica [{On} the theory of functions of a complex hyperbolic
  variable]}},
    language={Spanish},
        date={1935},
     journal={Univ. nac. La Plata. Publ. Fac. Ci. fis. mat.},
      volume={104},
       pages={139\ndash 183},
        note={\Zbl{62.1122.03}},
}

\bib{Wildberger05}{book}{
      author={Wildberger, N~J},
       title={Divine proportions: Rational trigonometry to universal geometry},
   publisher={Wild Egg Books},
        date={2005},
}

\bib{Yaglom79}{book}{
      author={Yaglom, I.~M.},
       title={A simple non-{E}uclidean geometry and its physical basis},
      series={Heidelberg Science Library},
   publisher={Springer-Verlag},
     address={New York},
        date={1979},
        ISBN={0-387-90332-1},
        note={Translated from the Russian by Abe Shenitzer, with the editorial
  assistance of Basil Gordon},
      review={\MR{MR520230 (80c:51007)}},
}

\bib{Zejliger34}{book}{
      author={Zejliger, D.N.},
       title={{\cyr Kompleksnaya linei0chataya geometriya. Poverhnosti i
  kongruencii.} [{C}omplex lined geometry. {S}urfaces and congruency]},
    language={Russian},
   publisher={GTTI},
     address={Leningrad},
        date={1934},
}

\end{biblist}
\end{bibdiv}
\LastPageEnding
\end{document}